# Outer Billiards, Digital Filters and Kicked Hamiltonians

G.H.Hughes

In 1978 Jürgen Moser [MJ4] suggested the outer billiards map (Tangent map) as a discontinuous model of Hamiltonian dynamics.  A decade earlier, J.B. Jackson [JKM] and his colleagues at Bell Labs were trying to understand the source of self-sustaining oscillations in digital filters . Some of the discrete mappings used to describe these filters show a remarkable ability to 'shadow' the Tangent map when the polygon in question is regular.

In [H2] the author outlines a theory of outer billiards dynamics for regular polygons and here we explore the connection between that theory and digital filters.

In 1997 Peter Ashwin [AP] showed that the digital filter map is equivalent to a sawtooth version of the Standard Map. For regular polygons, this provides a link between the global dynamics of the Tangent map and the toral dynamics of the Sawtooth Standard Map.

A.J.Scott,*et al.* [SHM] used a version of the kicked harmonic oscillator from quantum  mechanics, to create a mapping with also 'shadows' the Tangent map for regular polygons, but in a global form. Peter Ashwin noted that the complex valued version of this map is a special case of a Goetz Map [G1]. We have modified this map to obtain a more consistent conjugacy with the Tangent Map.

It is our premise that the geometry revealed by the Tangent map is intrinsic to the polygon itself and there is much evidence that points in that direction - including the contents of this paper which implies that the Tangent map is just one of many ways to illuminate this structure. Mathematica code is available for all mappings either here and at DynamicsOfPolygons.

This is a survey paper which covers a number of different topics so the reader is encouraged to use the table of contents to skip around. We have included a brief introduction to Hamiltonian dynamics in Section 4.

**Section 1: Digital filters and the outer billiards map**

Below is the schematic for a second order Digital Filter with two feedback loops. Circuits such as these are fundamental building blocks for more complex structures such as Analog to Digital Converters. J.B. Jackson and his associates at Bell Labs noted that such circuits can display large oscillations even when no input is present. This is due to non-linearities which can occur when the accumulator overflows.

The circuit shown here consists of three registers with a time delay of one unit between them. The intermediate outputs y(t) and y(t+1) are multiplied by *b* and *a* and fed back in where they are added to the contents of the accumulator.

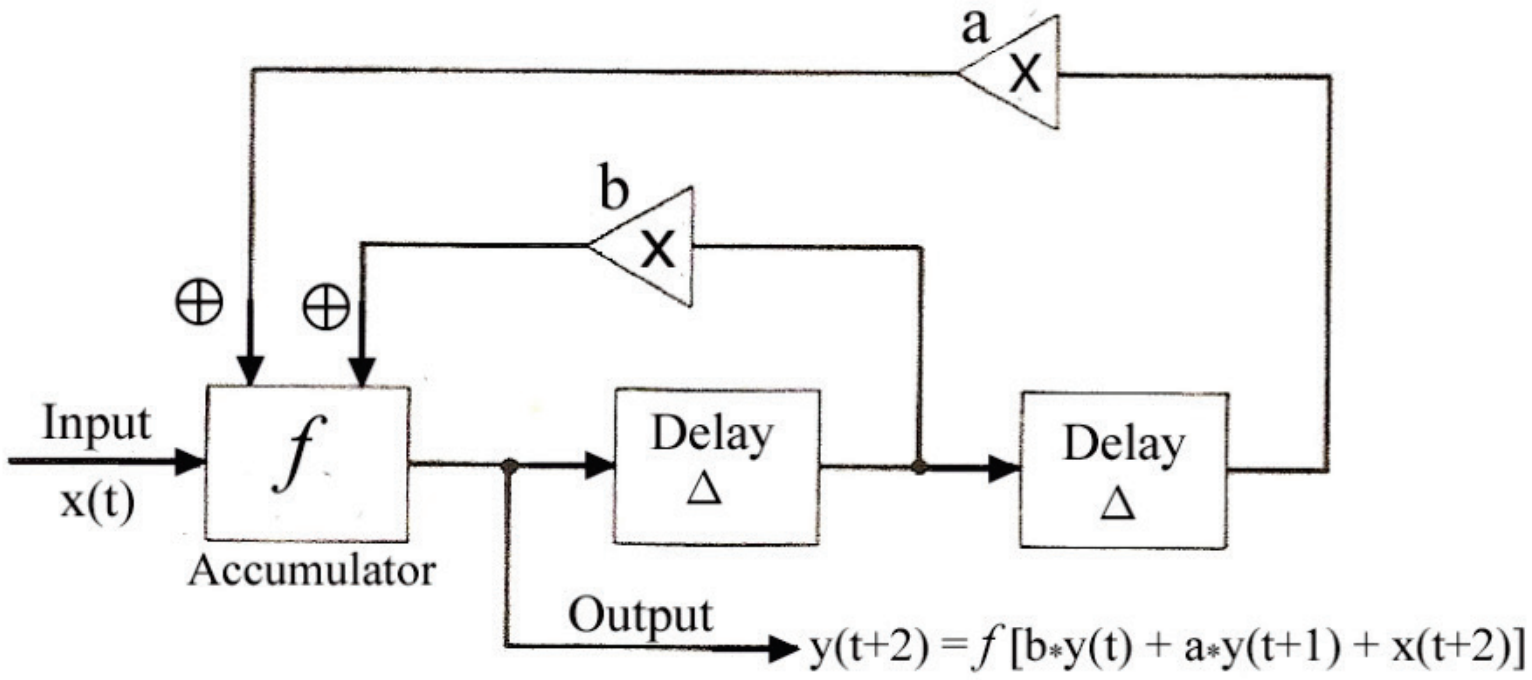

The equation for the output is $y(t+2) = f[\,by(t) + ay(t+1) + x(t+2)]$ which is a second-order difference equation. Under ideal conditions the function *f* would be the identity function but since the registers have finite word-length, there is the issue of possible overflow.

Following studies by Chua & Lin [C1] and Anthony Davies [D], we are interested in the self-sustaining oscillations which can occur even when there is no input, so we will assume $x(t) = 0$ for all t. The corresponding second-order equation can be reduced to two first order equations by setting $x_1 = y(t)$ and $x_2 = y(t+1)$. Then at each time tick, $x_1 \rightarrow x_2$ (by definition) and $x_2 \rightarrow f[bx_1 + ax_2]$. If the accumulator has no overflow, *f* is the identity function and the state equations are

$$\begin{bmatrix} x_1(k+1) \\ x_2(k+1) \end{bmatrix} = \begin{bmatrix} 0 & 1 \\ b & a \end{bmatrix} \begin{bmatrix} x_1(k) \\ x_2(k) \end{bmatrix} \text{ or } X(k+1) = AX(k) \text{ where } A = \begin{bmatrix} 0 & 1 \\ b & a \end{bmatrix} \text{ and } X = \begin{bmatrix} x_1 \\ x_2 \end{bmatrix}$$

In this linear system the fixed point is at the origin and the local behavior is determined by the eigenvalues which are $(a \pm \sqrt{a^2 + 4b})/2$. They are complex when $b < -a^2/4$ and in this case they have the form $\rho e^{\pm i\theta}$ where $\rho^2 = -b$ and $\theta = \arccos(a/2)$ When $|b| < 1$ solutions will damp out in time and when $|b| > 1$, solutions will diverge. When b = -1 the origin is a 'center' and solutions will rotation by θ. This is what circuit engineers call a 'lossless resonator' or a digital oscillator. So the b term controls 'damping' and the *a* term controls the frequency. In a digitally tuned oscillator for a radio, the *a* term would be adjusted to match the frequency of the station so that they would be in 'resonance'.

The registers are assumed to have fixed length so the largest possible value would be 0.11111111 ... $\approx$ 1. Negative number are typically stored in 2's complement form so the smallest possible value is 1.00000000 ... = -1. This means that the accumulator function *f* has the form of a sawtooth: ***f*[z_]:= Mod[z + 1,2] - 1**:

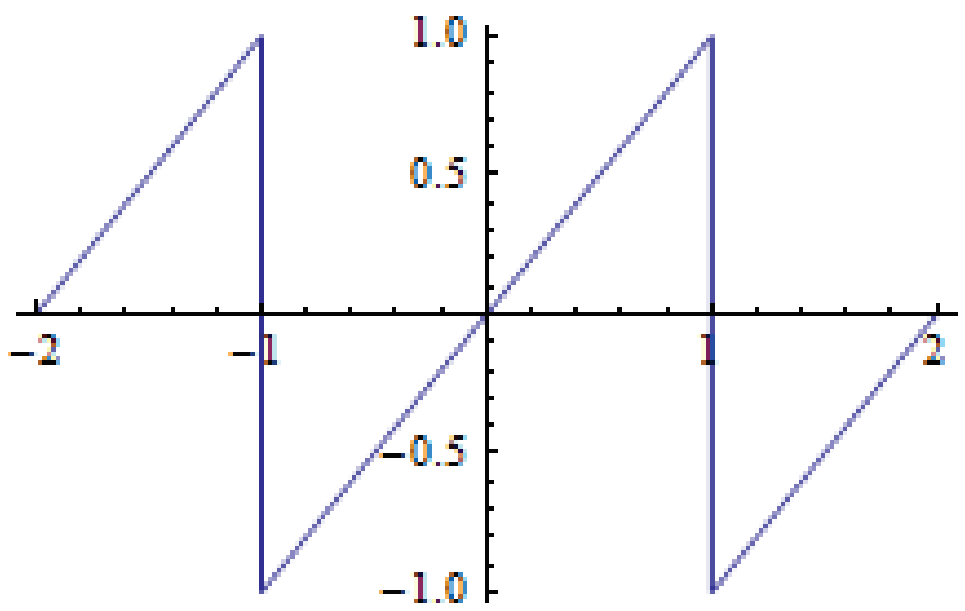

The discontinuity occurs in two possible scenarios: (i) Overflow occurs when adding one bit to 0.11111111 ... $\approx$ 1 and it becomes 1.00000000 ... = -1 . (ii) Underflow occurs when subtracting one bit from 1.00000000 ... = -1 and it becomes 0.11111111 ... $\approx$ 1. In both case the jump is of magnitude 2. The problem is that 1 and -1 are only one bit apart.

The 'ramps' shown here look smooth but actually they have discrete steps - one for each bit. For large register size, *f* is a good approximation. With the assumption that b = -1, then $\theta = \arccos(a/2)$ and the lone parameter can be either *a* or θ. The system equations are:

$$\begin{bmatrix} x_{k+1} \\ y_{k+1} \end{bmatrix} = f\left( A \begin{bmatrix} x_k \\ y_k \end{bmatrix} \right) = \begin{bmatrix} f(y_k) \\ f(-x_k + a y_k) \end{bmatrix} = \begin{bmatrix} y_k \\ f(-x_k + a y_k) \end{bmatrix}$$

where *f* is defined above. Note that $f(y_k) = y_k$ because $y_k$ is in range by assumption.

**Definition**: The Digital Filter map Df: $[-1,1)^2 \rightarrow [-1,1)^2$ is defined as

Df[{x,y}]:={y, *f* (-x + *a*y)} where *f* (z) = Mod[z+1,2]-1

In general the Df map defines a one-parameter family of maps and that parameter could be the winding number (rotation number) ρ or the angle θ = 2πρ.

Since Det[A] = 1, this map preserves area so it is a symplectic map. Symplectic maps on tori have been an area of interest to mathematicians and physicists since Henri Poincare (1854-1912) realized their value in the analysis of conservative (Hamiltonian) systems.

Any linear map of the form A= $\begin{bmatrix} 0 & 1 \\ -1 & a \end{bmatrix}$ will always be symplectic. This implies that the eigenvalues have the form λ and 1/ λ and as we noted above, $\lambda = (a + \sqrt{a^2 - 4})/2$ so for $a \in (0,2)$ the eigenvalues are complex with unit absolute value so $\lambda = e^{2\pi i\theta}$ where θ = arccos(*a*/2).

This implies that A represents a rotation, but it is an 'elliptical' rotation which can be conjugated to a pure rotation. When studying the dynamics of maps based on matrices such as A, if the trace *a* is the solution to a polynomial equation of low degree, there are computational advantages (exact arithmetic) to leaving A in its original form and we will work with A in both its original form and the conjugate form.

The Jordan normal form of A is a rotation matrix R: $\begin{bmatrix} 0 & 1 \\ -1 & 2\cos\theta \end{bmatrix} \sim \begin{bmatrix} \cos\theta & \sin\theta \\ -\sin\theta & \cos\theta \end{bmatrix}$

Note that R is actually a rotation by $-\theta$ which is perfect for the canonical clockwise rotations of the Tangent Map.The conjugating matrix is

$$G = \begin{bmatrix} 1 & 0 \\ \cos\theta & \sin\theta \end{bmatrix} \text{ with } G^{-1} = \begin{bmatrix} 1 & 0 \\ -\cot\theta & 1/\sin\theta \end{bmatrix}$$

G and its inverse will be used to go back and forth between Df space and the Euclidean space of the Tangent map.

For a given *a*, Df is a piecewise isometry with three primary regions (atoms) which can be labeled 1 (overflow), 0 (in bounds), or -1 (underflow).The equations for these atoms are:

$$S[\{x,y\}] = \begin{cases} 1 & \text{if } -x+ay \geq 1 \\ 0 & \text{if } -1 \leq -x+ay < 1 \\ -1 & \text{if } -x+ay < -1 \end{cases}$$

**Example**: Setting $a = 2\cos(2\pi/14)$ corresponds to a 'twist' of 1/14. The three atoms A,B and C are shown below with A being the overflow region and C the underflow. The Df map applies a (clockwise) elliptical rotation of $\theta = 2\pi/14$ to each region and the sawtooth nonlinearity *f* provides the corresponding translation – which is vertical by -2, 0 and + 2 respectively for A,B and C as shown on the right below.

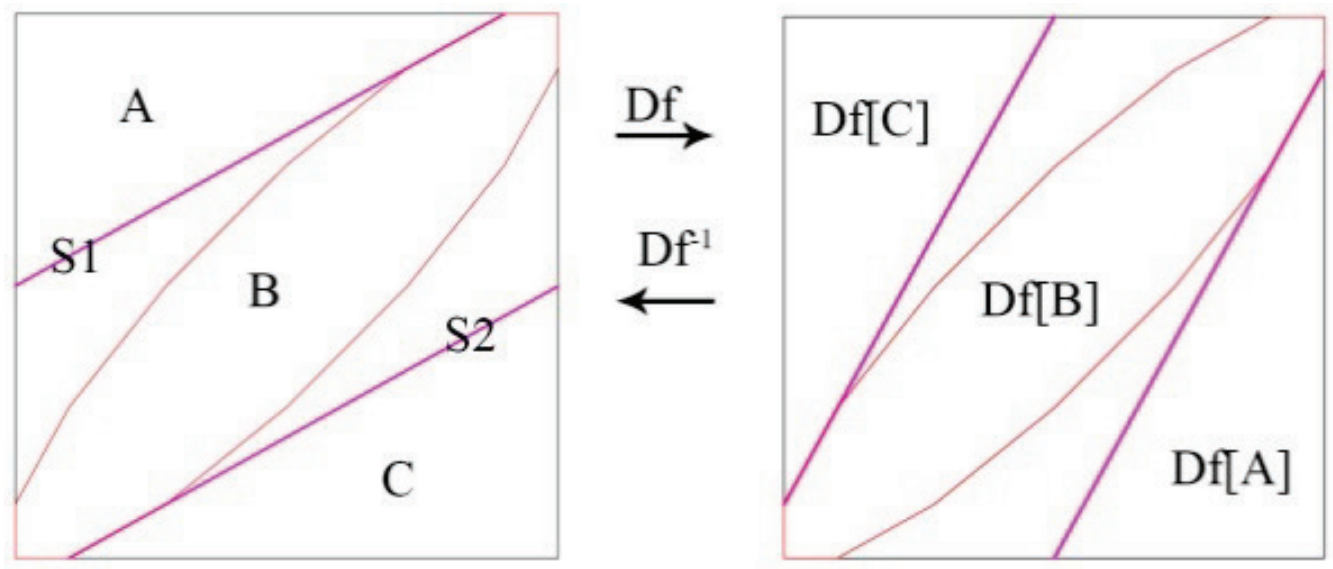

The seperatices S1 and S2 define the maximum extent of the linear center rotation, so they become extended edges of the limiting 14-gon. The bounds of the 2-torus are also extended edges of the limiting 14-gon so S1 and S2 map to these bounds. Any subset of these can serve as the initial level-0 singularity set –depending on what region is of interest. At this time there is no theory that can be used to determine a 'minimal' level-0 set for Df or for τ.

The slope of S1 (or S2) determines the number of 'steps' involved in a rational rotation and of course each 'step' is a rotation by θ. Our convention (described in detail below) is to define the number of steps using the (clockwise) difference between S1 (or S2) and the 'next' bounding edge – so the example above is step-1. This matches the 'twist' ρ which is 1/N.

Applying the map G from above, yields Df in Euclidean space where the elliptical rotation becomes a true rotation and the translations have magnitude 2/sinθ to match the edges of the bounding rhombus. In the rectified plots shown here, we have added a rotation by π/2 to make them match traditional Tangent map space. This means that the displacements of A and C are now horizontal of magnitude +2/sinθ and -2/sinθ respectively.

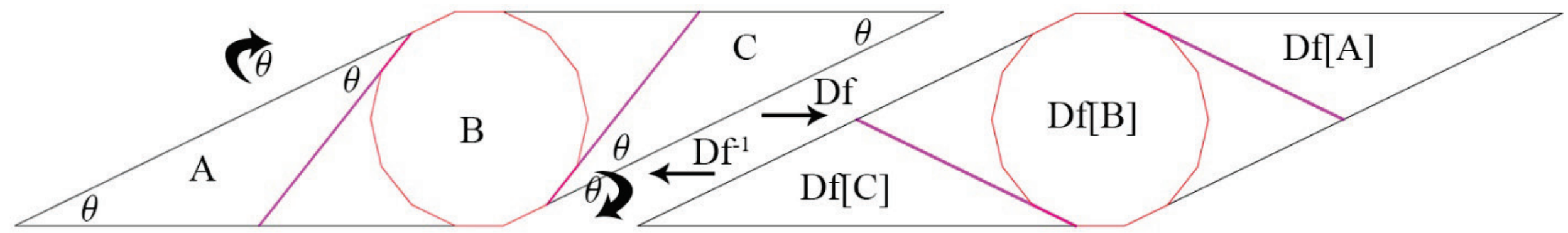

Therefore on each iteration of the web, the seperatices are rotated by θ. In this example, θ was the exterior angle of N = 14, so the resulting web is step-1 which is consistent with the web generation algorithms used for the Tangent map.

**Example**: Below are the first few iterations of the Df web in magenta and the rectified web in blue for $a$ = 2cos(2π/14).This web is obtained by mapping the two magenta seperatices S 1 and S2 under Df. The last iteration is the level 100 web. The Mathematica code for the level-100 web is given below.

**S1=Table[{x,(x+1)/a},{x,-1,a-1,.001}]; S2=Table[{x,(x-1)/a},{x,1-a,1,.001}]; S0=Join[S1,S2];**
**DfWeb = Flatten[Table[NestList[Df,S0[[k]],100],{k,1,Length[S0]}],1];**
**TrWeb = DfToTr[DfWeb] = RotationTransform[Pi/2][IShear/@DfWeb]];** (where IShear is $G^{-1}$ from above, so **IShear[{x_,y_}]:={x,-x *Cot[2Pi/N]+(y /Sin[2Pi/N])};**)

**Graphics[{AbsolutePointSize[1.0], Magenta, Point[DfWeb], Blue, Point[TrWeb]}].**

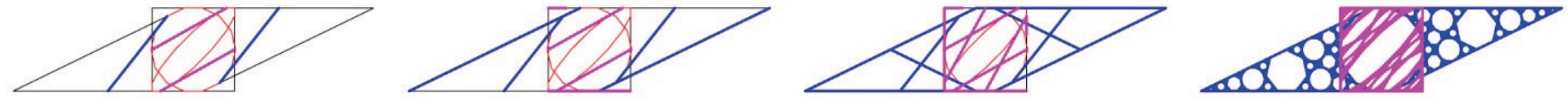

The blue web is a perfect copy of the Tangent map web for N = 14. This is not surprising because in the rectified mapping, the seperatices S1 and S2 become extended 'trailing' edges of N = 14 as shown on the left . The web for N = 14 is formed by iterating these edges under τ and these iterations involve a simple shear and rotation of these edges.

**Definition of the Tangent map**: For an initial point p, each iteration of $\tau$ is a reflection (central symmetry) about one of two possible 'support' vertices, so the formula is $\tau(p) = 2c - p$ where c is the support vertex. We will usually assume a clockwise orientation for the original polygon. The inverse map $\tau^{-1}$ is $\tau$ applied to the polygon with opposite orientation, so for a regular polygon, $\tau$ is essentially its own inverse.

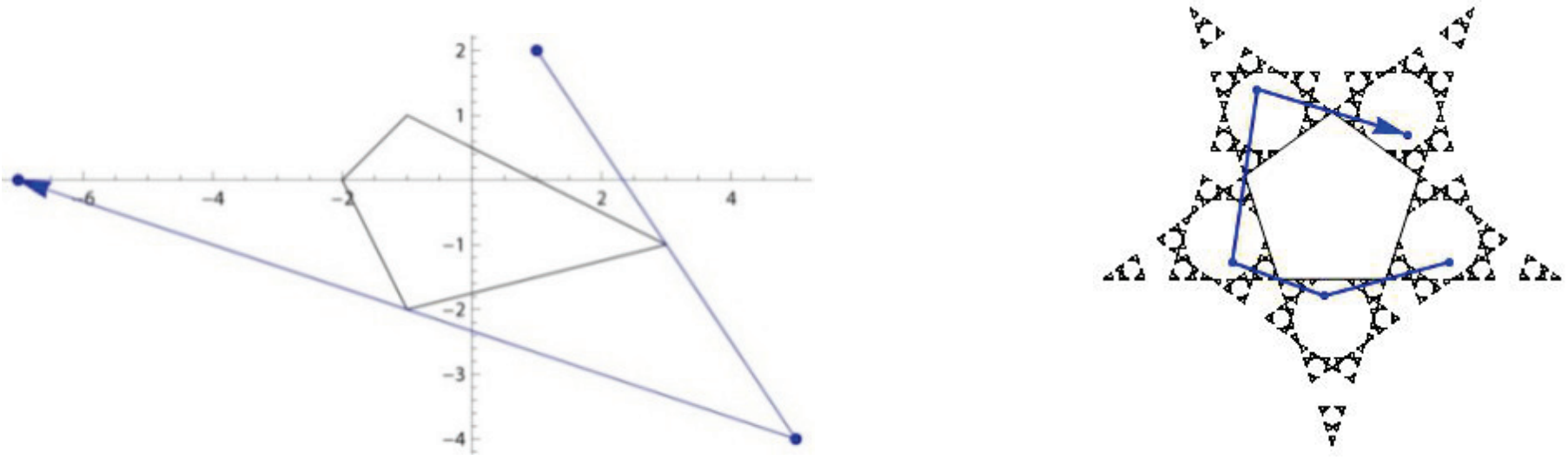

The map $\tau$ is not defined on the extended 'trailing edges' of the polygon and $\tau^{-1}$ is not defined on the extended forward edges, so the union of these extended edges form the 'level-1' singularity set (web). This level-1 web is shown on the left below. If the edges are truncated at the intersections, the web forms a 'star polygon' with Schlafli form {14,6}. This level-1 web provides the bounds of the First Family for N = 14 – as shown in cyan. These 'tiles' are the major 'resonances' of $\tau$ so they arise early in the web generation process. The tiles with the magenta centers are the nucleus of the First Family – they are called S[1] through S[6]. The largest tile, S[6], is also known as $D_{14}$ or simply D when N is understood.

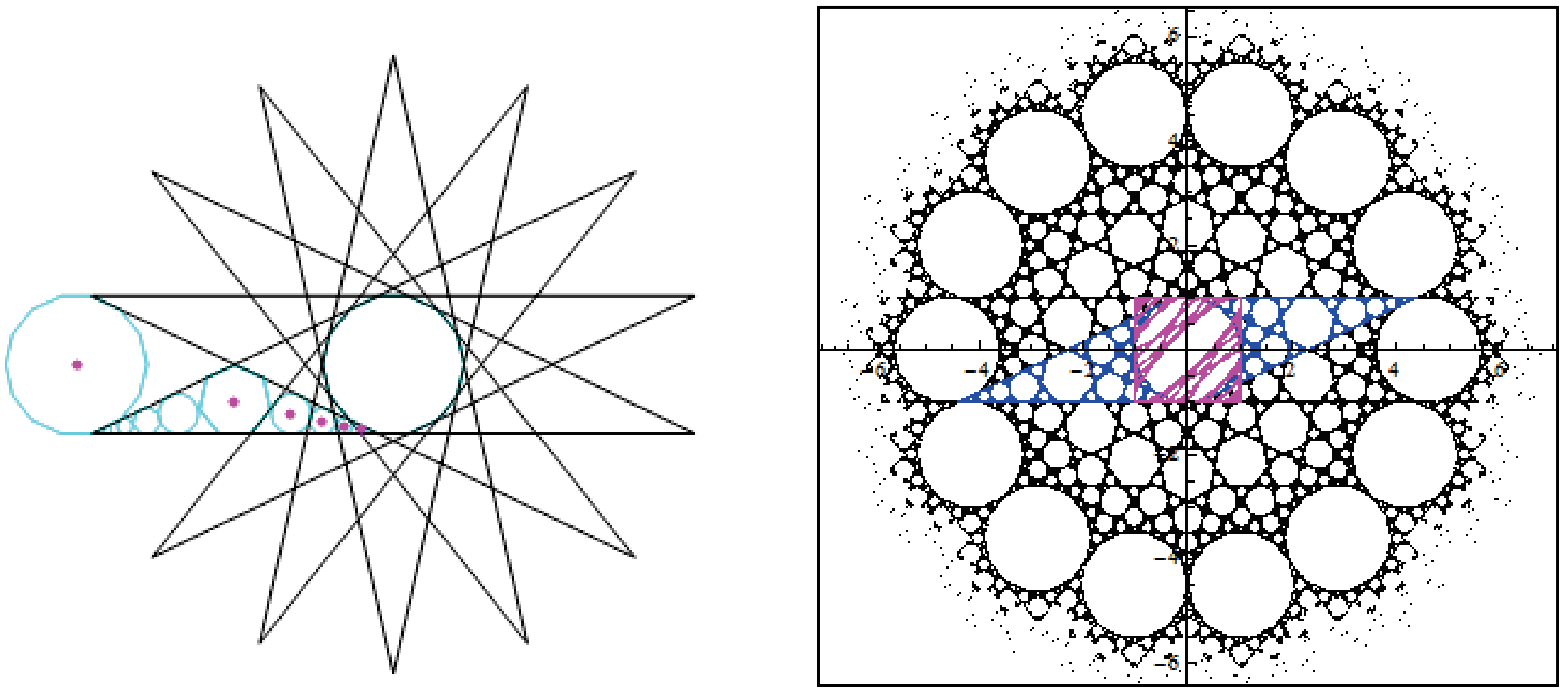

On the right is the level-100 web in black – showing a ring of 14 D tiles. There are endless rings of these maximal D tiles and the region between rings is invariant – so they guarantee that all orbits are bounded [VS]. The region inside the first ring is called the central 'star' region, and it serves as a template for the global dynamics. By symmetry, the blue rhombus can also serve as a template for the global web – so the Df map generates a perfect 'template' for the global web.

Even though the Df web is a perfect match for the Tangent map web, the dynamics are quite different. The example below shows the contrast between the 'toral' dynamics of Df and the 'global' dynamics of $\tau$. In many ways the Df dynamics are more 'natural' because the torus collapses the ring of D tiles down to N and this matches the symmetry of $\tau$ when N is even.

**Example**: (N = 14) The arrow on the left points to a small First Family tile at the foot of D - which is known as DS[1]. This tile has a period 70 orbit shown by the green dots. Note that there are two green dots in the blue rhombus. They evolve symmetrically under Df , so we will follow just the first point. The matching point in Df space is shown at right. It maps to the 'in-bounds' region on the right of N = 14 and therefore it rotates in a step-1 fashion until it experiences an 'underflow' and is reset to the initial point. Therefore the orbit is period 6. These 6 points and the symmetric 6 points are shown in red. (The initial green points have not been recolored.) Therefore the 70 green points are now 12 in Df space.

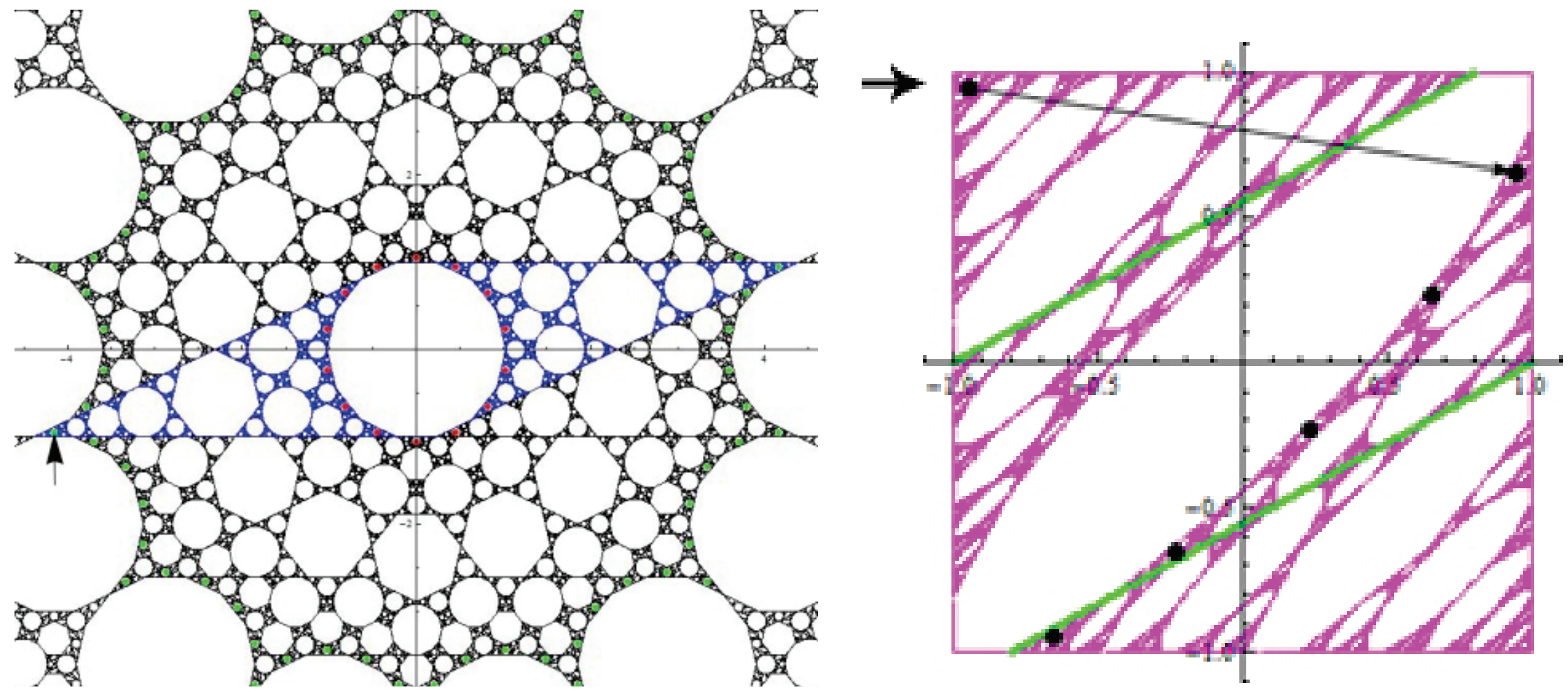

Back in '$\tau$-space', the 14 red points are the orbit of a 'clone' of DS[1]- called S[1]. If these two orbits are united there are now 84 congruent tiles and this yields a perfect 'global'-'toral' correspondence because 84/7 = 12. Note that the blue rhombus contains exactly 1/7 of the total tiles in the 'star' region above. To see this, rotate the blue rhombus 7 times and note that the overlap between rotated versions is identical to the region not covered. This 'global' – 'toral' correspondence is known as the 2kN Lemma.

**The 2kN Lemma**: For the Df map with $\theta = 2\pi/2N$, every periodic tile in Df with period $k > 2$ accounts for 2k congruent tiles in Df space and 2kN congruent tiles in the 'star' region of $\tau$. In this star region, the centers of these tiles may have different periods but the sum of the periods must be 2kN. (If a tile in Df space has period 2, then the center point p must map to $-p$ so there are only 2 congruent tiles in Df space and therefore only 2N congruent tiles in the star region.)

The 2kN Lemma is very general in scope but it can be challenging to obtain τ-periods from tile counts of the 'star' region. One special case where this is feasible is the S[k] and DS[k] First Family tiles. It is easy to find the periods of the S[k] tiles, but the periods of the DS[k] are not as obvious – yet symmetry says that the Df tile count of these two sets must match when N is even.

**Lemma:** For a regular N-gon for N-even, the number of congruent S[k] and DS[k] tiles in the First Family rhombus is N/2-k in both cases – if the central S[N/2 -2] tile is assumed to belong to both sets. Therefore the Df period of any S[k] or DS[k] tile is N/2-k. This is the 'k' in the 2kN Lemma

**Proof**: For any N-gon with N even, the symmetry of the S[k] and DS[k] families dictates that they have the same number of tiles in the First Family rhombus, so it is sufficient to count the S[k] tiles. The S[k] tiles can be arranged in a concentric fashion around S[N/2-2], so #S[k] = N/2-k . Df[cS[k]] will map these N/2-k tiles to themselves and this is true for DF[cDS[k]] also.(Therefore the combined rhombus count for any S[k] and DS[k] can be obtained by doubling the S[k] period, except for the combined S[N/2-2] & DS[N/2-2] which is not doubled.

**Example**: On the left below, the S[k] tiles for N = 14 are colored and the DS[k] are white. If the central S[5] is also counted as DS[5], then these two sets are equal. To count the S[k], perform surgery by slicing the rhombus in half and sliding. This yields the right image which unites the S[k] and the DS[k]. Clearly #S[k] = N/2-k and the DS[k] have the same formula. Under Df, the colored and uncolored tiles map to themselves, so the tile counts are the same as the Df periods. Note that the distinction between N = 14 and D vanishes and this tile could be included in the counts if desired.

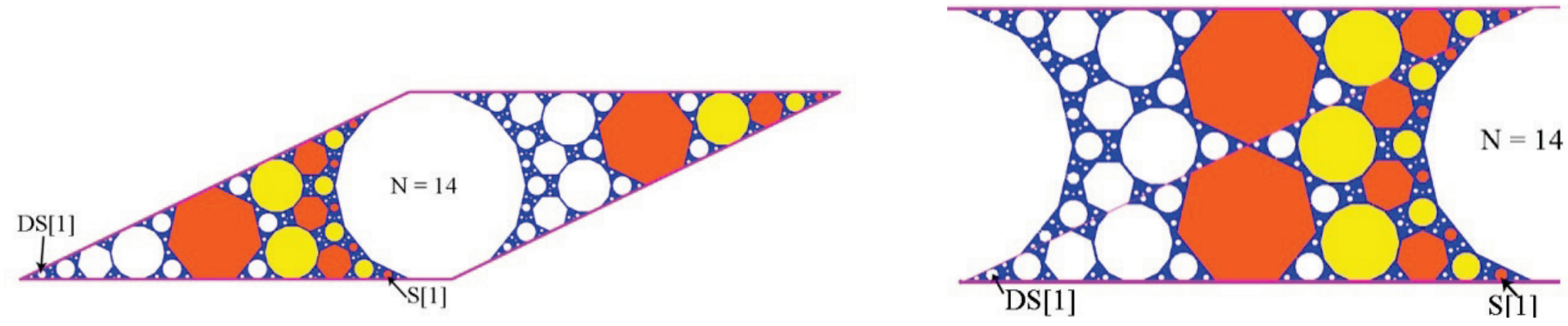

On the previous page we saw that the Df orbit of DS[1] was period 6. DS[1] was mapped to the 5 congruent tiles on the right of N = 14 – which could be called S[1]'s or DS[1]'s, but our convention is to call them DS[1]'s. So Df maps the 6 DS[1] tiles to themselves. The orbit of S[1] will be symmetric and once again, the S[k] will map to themselves as shown in color. Together there are 12 combined DS[1] and S[1] tiles in the Df rhombus, so there are 7*12 such tiles in the star region. (Note that this count could be made without recourse to the Df map.)

The challenge is to reconcile these counts with the τ-periods. This is easy for N twice-prime but not in general. To see the correspondence for N = 14, below is a table showing the periods of the First Family tiles for N = 14 (excluding S[6])

| Tile | S[1] | S[2] | S[3] | S[4] | S[5] | DS[1] | DS[2] | DS[3] | DS[4] | DS[5]=S[5] |
|---|---|---|---|---|---|---|---|---|---|---|
| Period | 14 | 7 | 14 | 7 | 14 | 70 | 56 | 42 | 28 | 14 |

It does not matter if the S[k] orbits decompose or not – the tile count for the inner star region will always be N. This means that the DS[k] count can be determined from the (known) total count, as shown below for N = 14.

| Tiles | S[1]&DS[1] | S[2]&DS[2] | S[3]&DS[3] | S[4]&DS[4] | S[5]=DS[5] |
|---|---|---|---|---|---|
| Total Count | 2*6*7 = 84 | 2*5*7 = 70 | 2*4*7 =56 | 2*3*7 = 42 | 2*1*7 = 14 |
| DS count | 2*5*7 = 70 | 2*4*7 = 56 | 2*3*7 = 42 | 2*2*7 =28 | 2*1*7 = 14 |

To find the DS[k] periods, there is still the issue of possible decomposition, but for N twice-prime, there is no decomposition of the DS[k]. Therefore the DS counts above are also the periods.

**Lemma**: For N twice-prime, the period of DS[k] = N*(N/2-(k+1)). There is a similar formula for N prime: period of DS[k] = N*(N-(k+2)).

These formulas can be used for any N value by making adjustments for decomposition of the DS[k] (and S[k]) when GCD[k,N]>1. For example if the twice-prime formula is applied to N = 18, the predicted periods of the DS[k] are : 126, 108, 90, 72, 54, 36 and 18. All of these are correct except that the periods of S[3] and S[6] are each reduced by a factor of 3. The situation is the same for N = 9, where only the periods of DS[3] and DS[6] are reduced by a factor of 3.

In general the 2kN Lemma yields results which are challenging to reconcile with the τ-map. For example the magenta M[2] tile shown below has τ-period 126 and this orbit includes the matching magenta M[2] tile on the right. The two black M[2]'s also map to each other with period 126, and there are congruent M[2]'s at the foot of N = 14 which form two period 28 orbits – for a tile count of 308. These are all part of a period 22 Df orbit, so 308 = 2*22*7.

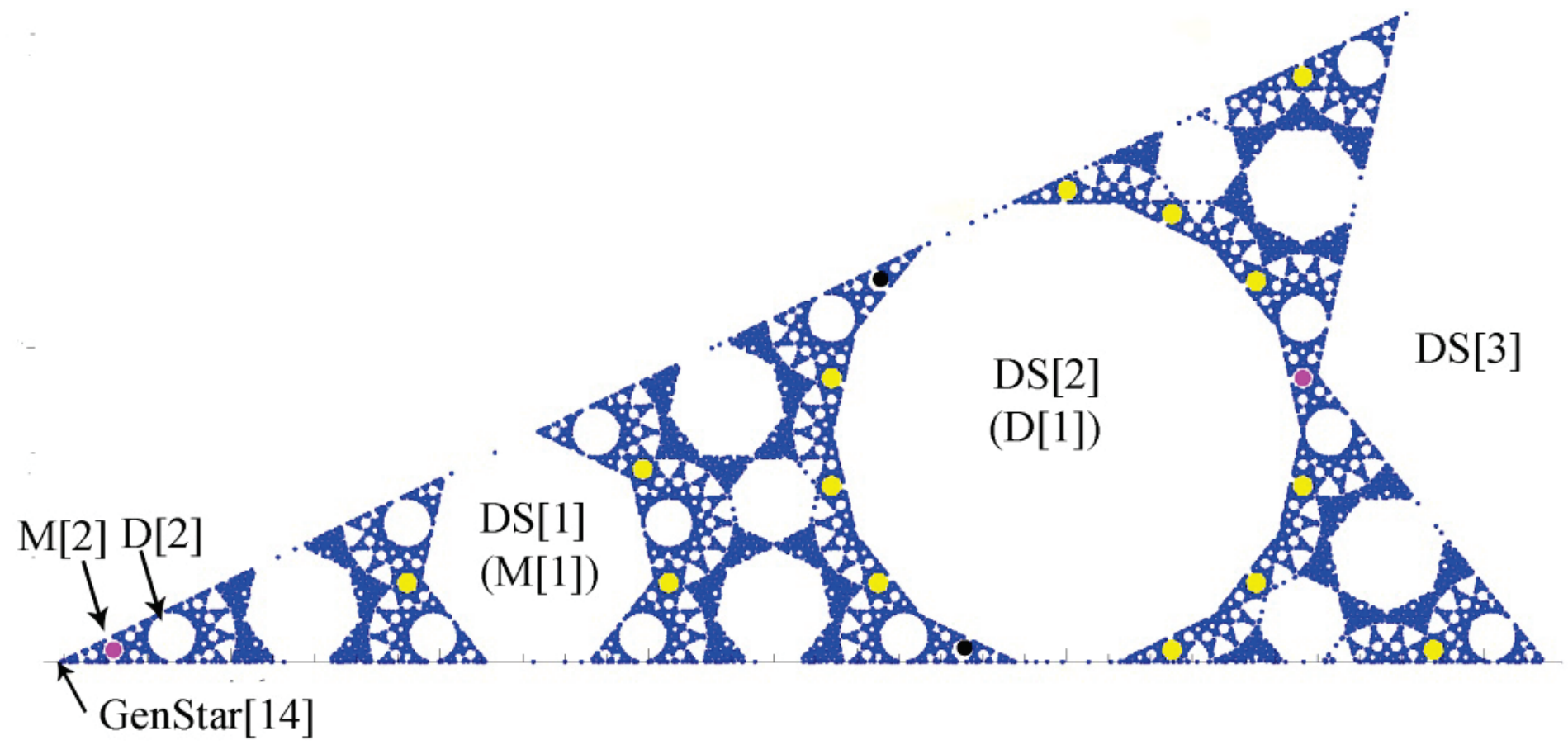

However this Df orbit does not include the 14 yellow M[2]'s shown above. The matching Df orbit is period 86 so there must be 1204 such tiles in τ-space. These 1204 tiles are composed of two outer groups with 476 each and two inner groups with 126 each. When these yellow, magenta and black groups are combined together, they yield 1512 M[2]'s which factors as 108*14. The region shown above is the second generation for N = 14 and it is clearly not self-similar to the first generation – so there is no obvious renormalization which can be applied.

**Extensions of the Df Map**

To mimic the Tangent map for N = 14 required a rotation angle of 2π/14 – which corresponds to a 'twist'of 1/14. This step-1 Df map was a perfect match for the step-1 web of the Tangent map, but a polygon such as N = 14, will also allow step-2 and step-3 Df webs – using θ = 2(2π/14) and θ = 3(2π/14). The former is θ = 2π/7, so it would correspond to N = 7

Below is the step-2 web with θ = 2π/7 - so the corresponding twist is ρ = 2/14, The linear map does indeed generate heptagons, but they are not compatible with the boundary conditions – which require that opposite edges be parallel. Therefore the heptagon are 'chopped up' to match the N = 14 case, but the web is step-2 so it is very different from a normal Tangent map web.

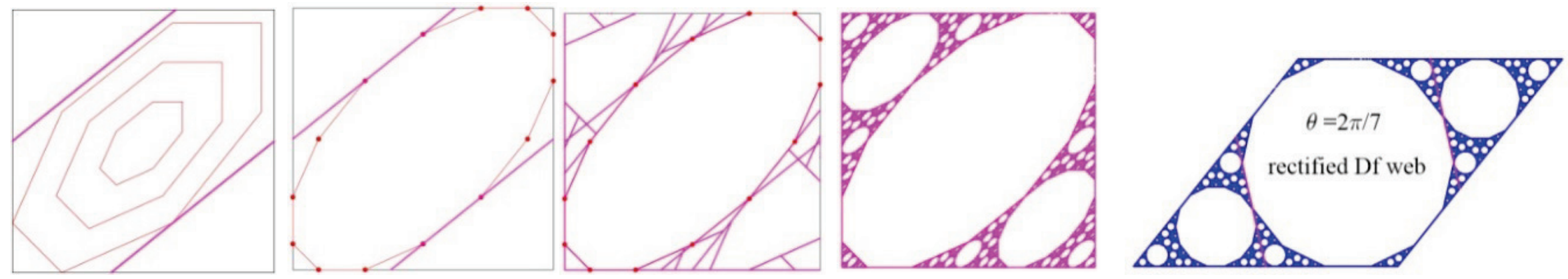

The step-2 relationship between S1 (or S2) and the 'next' edge of the rhombus can be seen in the blue rectified web on the right - where we have reproduced the image of S1 and S2 in magenta. In the case of N = 7 and N = 5, the resulting step-2 webs are recognizable subsets of the traditional Tangent map webs, but for the most part, these attempts at odd rotations are only distantly related to the desired web. As indicated earlier, N = 14 also supports a step-3 web. This is discussed in the Df Theorem below.

In general the Df map defines a one-parameter family of maps and that parameter could be the winding number (rotation number) ρ or the angle θ = 2πρ. We will define a rotation to be 'rational' iff ρ is rational. The following table gives some examples of rational rotations. ('short' families are defined below.)

| ρ (winding #) | 0 | 1/14 | 1/12 | 1/10 | 1/8 | 1/5 | 3/14 | 1/4 | 1/2 |
|---|---|---|---|---|---|---|---|---|---|
| Rotation θ | 0 | 2π/14 | 2π/12 | 2π/10 | 2π/8 | 2π/5 | 2π(3/14) | 2π/4 | 2π/2 |
| 2Cosθ | 2 | 1.80194 | √3 | (√5 +1)/2 | √2 | (√5 -1)/2 | 0.44504 | 0 | -2 |
| polygon | circle | N = 14 | N = 12 | N = 10 (& N = 5) | N = 8 | N = 10 (short) | N = 14 (short) | N = 4 | line |

By default we restrict θ ∊ (0, π/2] so ρ ∊ (0,1/4]. This covers the full range of regular polygons since the odd cases are included in the twice-odd cases . The highlighted region is what we call the 'quadratic' range because it includes the four 'quadratic' regular polygons- N = 5,8,10 & 12.

Below is a series of plots showing the dynamics in this range as *a* = 2cosθ increases from √2 (N = 8) to √3 (N = 12) with N = 10 in between at (√5 +1)/2. These are all 'rational' rotations because the winding numbers are decrements of 1/120- starting with N = 8 at θ = 2π(15/120).

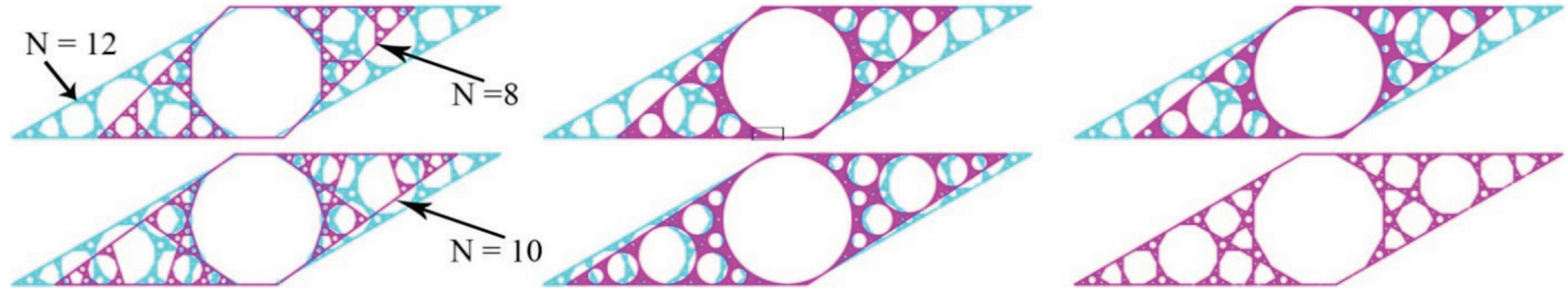

Therefore the 2nd magenta plot above has $\theta = 2\pi(14/120)$. As a twist map it is period 60 so it corresponds to a Df map with N = 60 and a 'step-7' web as shown in the detail below. This detail matches the small rectangle in the plot above.

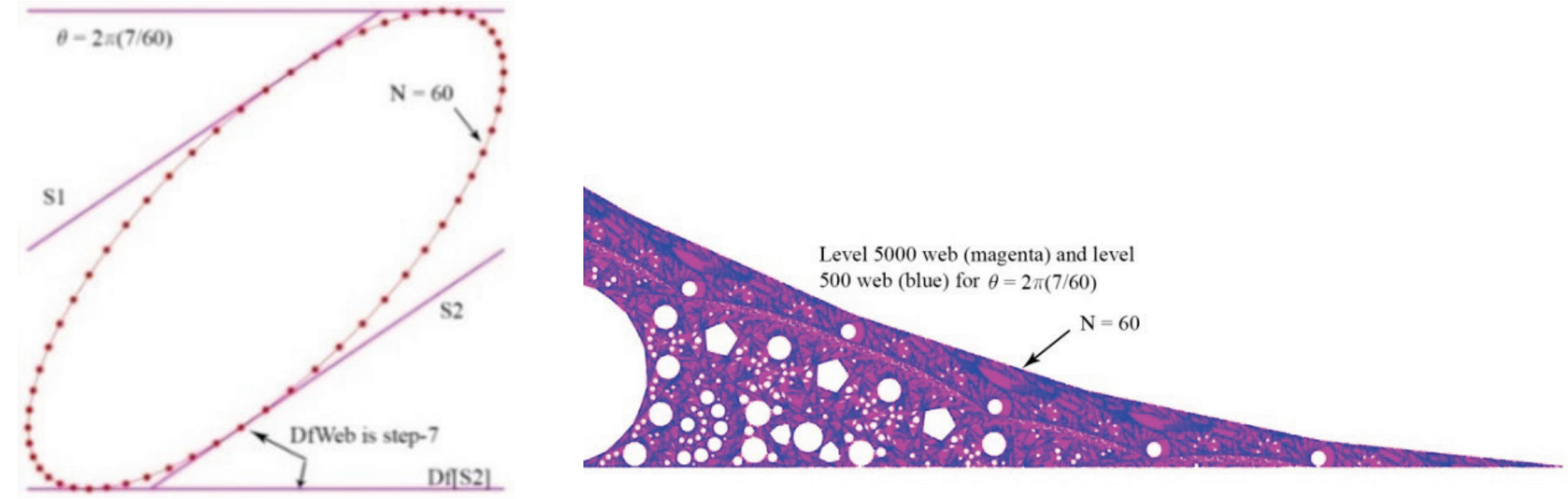

N = 60 has step-k webs ranging from 1 to 15, but the step-15 case reduces to N = 4 with null web, so the practical range extends to 14. It is not clear which of these webs are related to the traditional step-1 web - but the theorem below says that the 'maximal' 14-step web is special in that it reproduces part of the 'in-sito' dynamics of the central S[28] tile of N = 60.

**Lemma:** For the Df map with parameter $\theta = 2\pi\rho$, when the 'twist' $\rho = p/q$ with p and q relatively prime integers and $p/q \leq 1/4$, then the resulting Df web is conjugate to that of a regular q-gon with a step-p web when q is even or a regular 2q-gon with a step-2p web when q is odd.

Since $\rho \leq 1/4$, it follows that the maximal step size for any regular N gon is $\lfloor N/4 \rfloor$ - but as noted above, when N is twice-even this reduces to the null case of N = 4, so we will define the 'maximal' step size to be N/4-1 when N is divisible by 4. Most of these 'rational' webs are unrecognizable but there are two special cases which are described in the Theorem below.

**Theorem (Df map for Regular Polygons):** For a regular polygon N = 2k, the 'effective' range of rotation values for the Df map with $\theta=2\pi\rho$ are $\rho = 1/N$ to $\lfloor k/2 \rfloor/N$ for k odd or $(k/2-1)/N$ for k even. Therefore the Df webs range from step-1 to step-$\lfloor k/2 \rfloor$or step-$(k/2-1)$. These minimal and 'maximal' step values have webs which are described below:

(i) When $\rho = 1/N$ the Df web will be step-1 and this web is locally conjugate to the Tangent map web for N.

(ii) When k is odd the web in (i) will have a central S[k-2] tile which will be a regular N/2-gon. The Twice-Odd Lemma implies that the local web of this polygon will be conjugate to the Tangent map web for N/2, so the web for N includes a scaled copy of the web for N/2. When $\rho = \lfloor k/2 \rfloor/N$, the Df web will have maximal step-$\lfloor k/2 \rfloor$ and the rectified web will yield a shortened version of the step-1 rhombus from (i). This secondary rhombus will be locally identical to the step-1 Tangent map rhombus but contain just the first $\lfloor k/2 \rfloor$ tiles of the First Family of N. This is known as the Inner Family or 'short family'.

(iii) When k is even, the web in (i) will again have a central S[k-2] tile but now this tile will be a regular N-gon instead of a N/2-gon. Since N is a now a multiple of 4, the 'maximal' k/2-1 step will be odd iff N is a multiple of 8. In this case $\rho = (k/2 - 1)/N$ will have a Df web which is conjugate to the local web for the central S[k-2] tile.

To see why the short webs from (ii) and (iii) are faithful to the tangent map web, note that in both cases the gap between consecutive Df boundary edges is $\lfloor k/2 \rfloor+1$. In case (ii) this implies that the maximal step-$\lfloor k/2 \rfloor$ is retrograde step-1 as in the traditional Tangent map web. In case (iii) the 'maximal' k/2-1 step is retrograde step-2 and this matches the evolution of S[k-2] in the Tangent map web when N is twice even. This evolution was described [H2]. S[k-2] is generated by two interwoven step-2 cycles which define the even and odd edges of this N-gon When N is twice-odd, these cycles collapse to a single N/2-gon and yield the Twice-Odd Lemma.

**Examples:** The first few terms in the 'twice-odd' series of ρ are 1/6, 2/10, 3/14, 4/18 so the numerators are the positive integers and the denominators increase by 4. On the resulting rectified web we have included the image of the seperatices S1 and S2 in magenta. Note that these webs are all (retrograde) step-1 on the shortened rhombii.

| N = 10 : $\rho$ = *1/10, a* = 2Cos[2π/10] | $\rho$ = *2/10, a* = 2Cos[2π/5] |
|---|---|
| N = 10 | N = 10 |
| N = 14 : $\rho$ = *1/14, a* = 2Cos[2π/14] | $\rho$ = *3/14, a* = 2Cos[3π/7] |
| N = 14 | N = 14 |
| N = 18 : $\rho$ = *1/18, a* = 2Cos[2π/18] | $\rho$ = *4/18, a* = 2Cos[4π/9] |
| N = 18 | N = 18 |

The first few terms of the ‘twice-even’ series of ρ are 1/8, 3/16, 5/24, 7/32, 9/40 so the numerators are the odd integers and the denominators increase by 8. As above, we have included the images of S1 and S2 in magenta, and here they are (retrograde) step-2.

| N = 16 : $\rho = 1/16$, $a = 2\text{Cos}[2\pi/16]$ | $\rho = 3/16$, $a = 2\text{Cos}[3\pi/8]$ |
|---|---|
| N = 16<br>S[6] | S[6] |
| N = 24 : $\rho = 1/24$, $a = 2\text{Cos}[2\pi/24]$ | $\rho = 5/24$, $a = 2\text{Cos}[5\pi/12]$ |
| N = 24<br>S[10] | S[10] |
| N = 32 : $\rho = 1/32$, $a = 2\text{Cos}[2\pi/32]$ | $\rho = 7/32$, $a = 2\text{Cos}[7\pi/16]$ |
| S[14]<br>N = 32 | S[14] |

Sections 5, 6 and 8 will discuss the connection between these maps and Hamiltonian dynamics. If a matrix such as A was the Jacobian of a Hamiltonian system, the combination of complex eigenvalues and Det[A] = 1, would imply that the origin is a 'marginally stable' fixed point where the stability depends on θ. When the system is perturbed, the rational and irrational rotations would yield different dynamics and the KAM Theorem gives conditions for the survival of a non-zero measure of the irrational rotations.

The full Df map includes a perturbation, but the perturbing function $f$ is discontinuous and this is at odds with the KAM Theorem where any perturbation is assumed to be at least continuous. In fact Jurgen Moser proposed the Tangent Map as an example of a system that violated the continuity conditions of the KAM Theorem - but retained some of the same characteristics. It appears that the Df map is a perfect toral model of the Tangent map, when the polygon in question is regular.

The Df map and related maps with sawtooth nonlinearities are similar to 'kicked' Hamiltonians where the perturbing 'kicks' are periodic. If the kicks are relatively small, they can be used to model continuous perturbations such as planetary interaction. The Standard Map is a kicked Hamiltonian where the kick’s strength can adjusted to see what happens when the KAM Theorem breaks down.

In the Standard Map, the coordinates represent angular position (x) and angular momentum (y) and the perturbing function is naturally of the form $\text{K}\sin(x)$ . For Df, the perturbation frequency depends on the parameter *a* so it remains fixed. This means that the Df map can be 'tuned' and the Standard Map cannot. For the Tangent Map this is critical since the perturbation frequency is fixed by the regular polygon. There is a version of the Standard map which retains this ability to be 'tuned'. Peter Ashwin calls it the Sawtooth Standard Map and he shows that it is conjugate to Df. See Section 5.

**Section 2: Digital Filters in Three Dimensions**

Chua & Lin [C2] describe a third order Df filter that showed dynamics similar to Df but in three dimensions. The circuit diagram is given below.

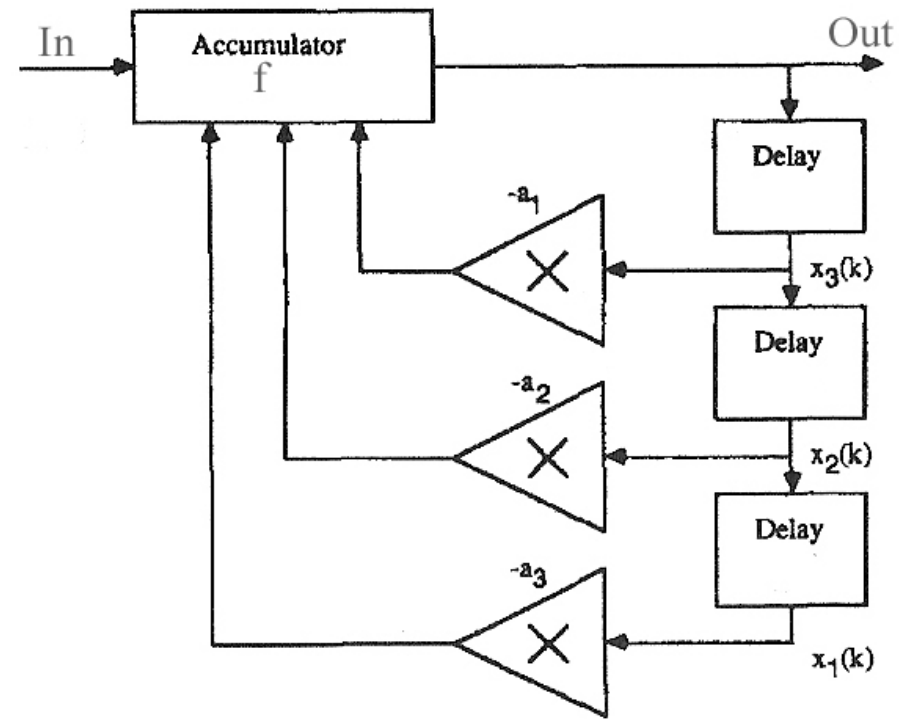

The difference equation for the output is: $y(t+3) = f[\,-a_3y(t) - a_2y(t+1) - a_2y(t+2) + u(t+3)]$ where u(t+3) is the input. If we assume that the accumulator is in 2's complement form, $f$ is unchanged from the second-order case. Since we are once again interested in self-sustaining oscillations, set u(t) = 0 for all t. The output equation can be reduced to three first order equations by setting

$x_1 = y(t)$, $x_2 = y(t+1)$, $x_3 = y(t+2)$. Then at each time tick $x_1 \to x_2$ , $x_2 \to x_3$ and $x_3 \to f[\,-a_3x_1 - a_2x_2 - a_1x_3]$ . The linear system is

$$X(k+1) = AX(k) \text{ where } \mathbf{X} = \begin{bmatrix} x_1 \\ x_2 \\ x_3 \end{bmatrix} \text{ and } \mathbf{A} = \begin{bmatrix} 0 & 1 & 0 \\ 0 & 0 & 1 \\ -a_3 & -a_2 & -a_1 \end{bmatrix}$$

with characteristic equation $\lambda^3 + a_1\lambda^2 + a_2\lambda + a_3 = 0$.

Setting $-(a+c) = a_1$ , $ac - b = a_2$ and $bc = a_3$ this can be factored as $(\lambda - c)(\lambda^2 - a\lambda - b)$ in which case the stable region is $|c| \le 1$, $|b| \le 1$, $|a| \le 1-b$

Choose c and b values from the boundary: c = 1, b = -1, then $|a| < 2$ and *a* is the lone parameter.

The linear evolution matrix is now $\mathbf{A} = \begin{bmatrix} 0 & 1 & 0 \\ 0 & 0 & 1 \\ 1 & -(1+a) & 1+a \end{bmatrix}$

The full non-linear system is D3f[{x_, y_, z_}] := {f[y], f[z], f[x - (1 + a)*y + (1 + a)*z]};

with *f*[x_] := Mod[x + 1, 2] - 1 as with Df. Once again we assume that y and z are 'in range' and the angular variable corresponding to *a* is 2*Cos[w];

**D3f[{x_, y_, z_}] := {y, z, f[x - (1 + a)*y + (1 + a)*z]};**

Just as with Df, the system equations can be written to make the 'kicks' explicit:

$$\begin{bmatrix} x_{k+1} \\ y_{k+1} \\ z_{k+1} \end{bmatrix} = \begin{bmatrix} 0 & 1 & 0 \\ 0 & 0 & 1 \\ 1 & -(1+a) & 1+a \end{bmatrix} \begin{bmatrix} x_k \\ y_k \\ z_k \end{bmatrix} + \begin{bmatrix} 0 \\ 0 \\ 2 \end{bmatrix} s_k$$

For a given initial point $\{x_0,y_0,z_0\}$, the sequence $s_k$ can be found by tracking its orbit with D3f, but now there are 7 possible values of $s_k$: $\{0 ,\pm1, \pm2, \pm3\}$. This is because the orbit of $\{x_0,y_0,z_0\}$ may visit as many as 4 planes. We will give Mathematica code below for finding the s-sequence of any point.

D3f is a map on the unit cube to itself.The dynamics of the linear system lie in an ellipse, as in the Df case. Evaluating the Jacobian, the three eigenvalues 1, $e^{i\theta}$, $e^{-i\theta}$ have unit value except when $|x_1 - (1+a)x_2 + (1+a)x_3|$ is odd (1,3,5,7). The mapping is area preserving except for these discontinuities.

To see the dynamics of the linear map , set **f[x_]:= x** and try **w = 1.0**, so a = 2*Cos[1] ≈ 1.0816.
**Orbit=NestList[D3fx,{1,0,0},1000]; ListPointPlot3D[Orbit]**

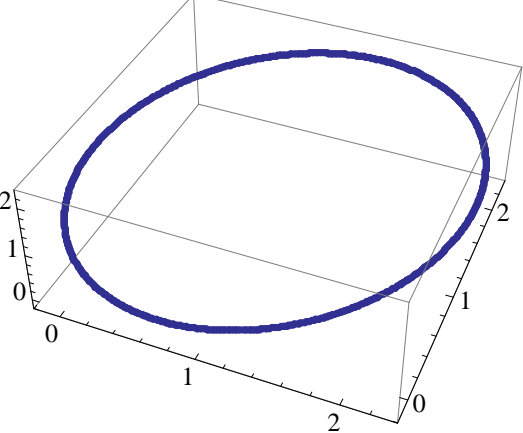

This ellipse can be aligned with the $x_1 = 0$ plane using:T3 = $\begin{bmatrix} 1 & 1 & 0 \\ 0 & Cos[\theta] & Sin[\theta] \\ 1 & Cos[2\theta] & Sin[2\theta] \end{bmatrix}$

**T3[w_] := {{1, 1, 0}, {1, Cos[w], Sin[w]}, {1, Cos[2*w], Sin[2*w]}};**
**D3fToTr[q_] :=Inverse[T3[w]].q ; ROrbit = D3fToTr/@Orbit; ListPointPlot3D[Orbit]**

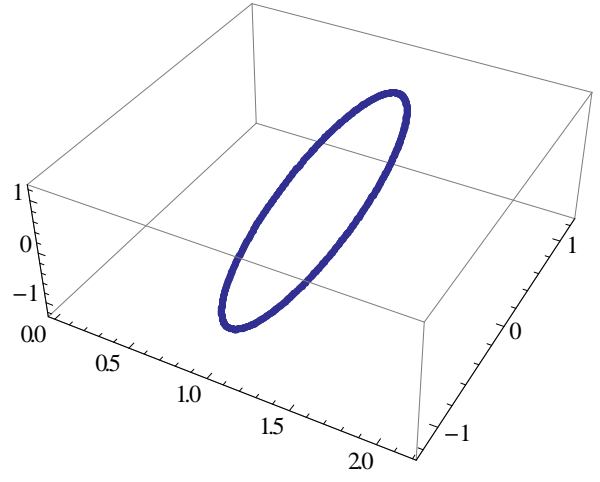

Now we can project this onto the $\{x_2, x_3\}$ plane: **P23=Drop[ROrbit,None,1]** (this drops no rows and the first column). **ListPlot[P23 ,AspectRatio->Automatic]**

Below is the orbit of this same initial point {1,0,0}, but with normal non-linear overflow function *f*[x_]:= Mod[x + 1, 2] - 1. **Orbit = NestList[D3f, {1, 0, 0}, 10000]; ROrbit = D3fToTr/@Orbit; P23=Drop[ROrbit,None,1]; ListPlot[P23]; Graphics3D[Point[Orbit]];**

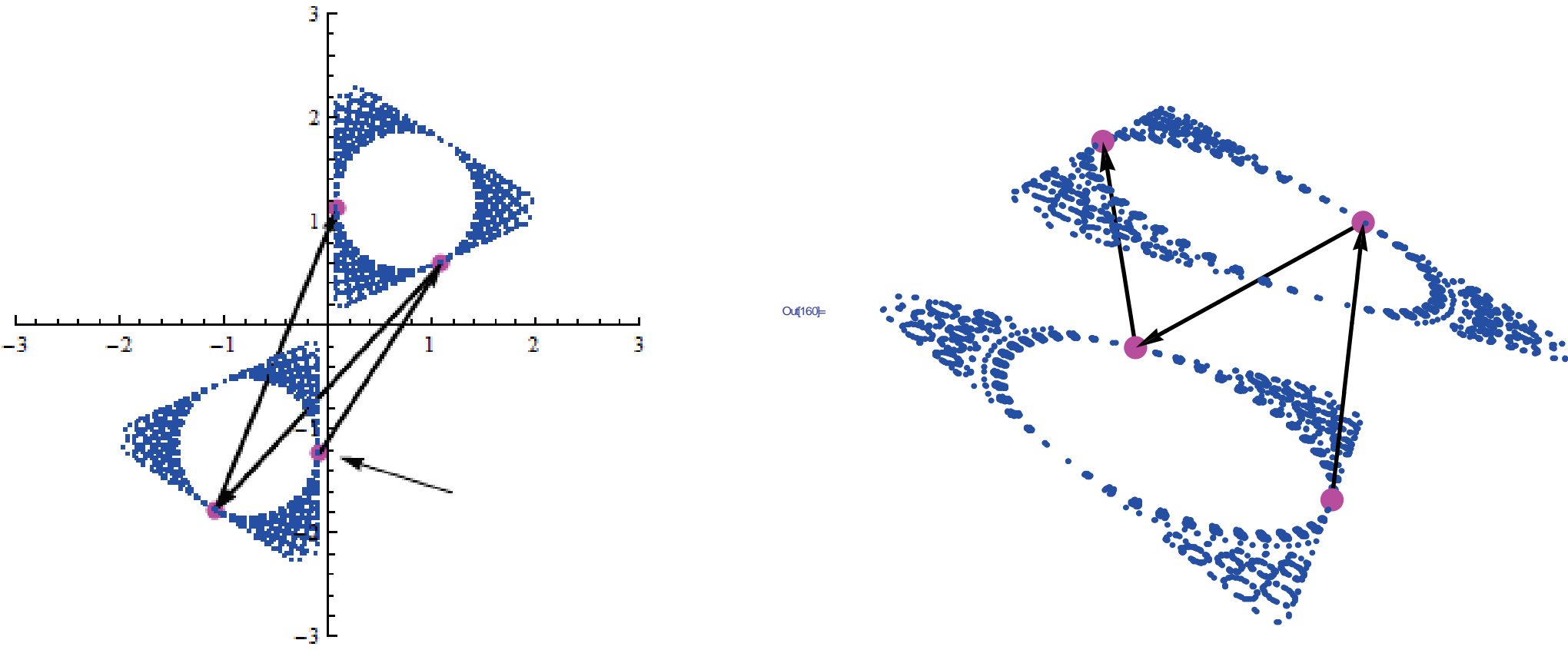

Points jump back and forth between the two planes, but not necessarily in strict alternation as above. The first 4 points in the orbit of {1,0,0} are shown in pink. In the $\{x_2,x_3\}$ projection on the left, the arrow shows the initial point. This is D3fToTr[{1,0,0}]={1.08767, -0.0876713,-1.23629} but with the first co-ordinate dropped. In the 3D plot, the points are the actual orbit, so the initial point shown is {1,0,0} (lower right).

The planes are parallel at distance 2/Sqrt[2+$a^2$ ] ≈ 1.1237. This is the magnitude of a 'kick' for $a$ =2*Cos[1]. So what we are seeing are the kicks from the s-sequence of this orbit. Chua gives the change of co-ordinates that can be used to restrict the map to these planes. In this manner the 3 dimensional dynamics can be analyzed using the same techniques as the 2 dimensional case.

Here is an example, using the orbit above. The new (orthogonal) basis vectors depend on the parameter $a$:
**e1 = (1/Sqrt[2 + a^2])*{1, - a, 1}; e2 = (1/Sqrt[2])*{1, 0, -1}; e3 = (1/Sqrt[4 + 2*a^2])*{a, 2, a};**
These will be the columns of the transformation matrix: **T = Transpose[{e1, e2, e3}];**

$$\textbf{MatrixForm[T]} = \begin{bmatrix} \frac{1}{\sqrt{2+a^2}} & \frac{1}{\sqrt{2}} & \frac{a}{\sqrt{4+2a^2}} \\ -\frac{a}{\sqrt{2+a^2}} & 0 & \frac{2}{\sqrt{4+2a^2}} \\ \frac{1}{\sqrt{2+a^2}} & -\frac{1}{\sqrt{2}} & \frac{a}{\sqrt{4+2a^2}} \end{bmatrix} \quad \textbf{TI = Inverse[T];}$$

The corresponding D3f map in T-space can be decomposed into two maps - one for the transport between planes, and one for the dynamics within the planes. The first co-ordinate of the T-space map tracks the jumps between the planes and since these coordinates are orthogonal, the jump distance is ± 2/Sqrt[2+$a^2$ ]. This means that the first coordinate has the form $u_{k+1} = u_k$ + 2/Sqrt[2+$a^2$ ]*s.

For two planes, the possible values of s are {0,1,-1}. Our current $a$ value is 2*Cos[1] ≈ 1.0816 which is slightly larger than 1, so there may be as many as four planes and this means the possible values for the s-sequence are:{0 ,±1, ±2, ±3}. There actually are 2 more small planes in this example, but they don't show up until after 25,000 iterations. They account for two of the missing vertices of the triangles above. We will plot them shortly, but for now the dynamics are restricted to just the 2 planes shown above.

**Example**: Using {1,0,0} as the initial point, we will use the T-space map to find the first few terms in the s-sequence. The first step is to find the T-space point corresponding to {1,0,0}: TI.{1,0,0} ={**0.561859**, 0.707107, 0.429318} so the first coordinate is u = .561859. Now advance the orbit by one iteration and repeat: D3f[{1,0,0] = {0, 0, -1}, TI.{0,0,-1} = {**-0.561859**, 0.707107, -.429318}. Note that the difference between the first coordinates is the distance between the planes, so $u_1 = u_0$ + 2/Sqrt[2+$a^2$ ]*s where s = -1. D3f[{0,0,-1}={0, -1, -0.0806046} and TI.[%] = {**0.561859**, 0.0569961, -0.829194}, so the s-sequence so far is {-1,1,..}. These are the kicks that we see above.

After 10,000 iterations the orbit is still confined to these two planes, but there are occasional occurrences of 0 in the s-sequence. For example, the 10,000th iterate of {1,0,0} is p0= Last[Orbit] ={-0.100544,-0.0934865,0.999522}, and if we track this point for 10 more iterations we get a mapping from the lower plane to itself as shown below. These become much more common as the complexity increases.

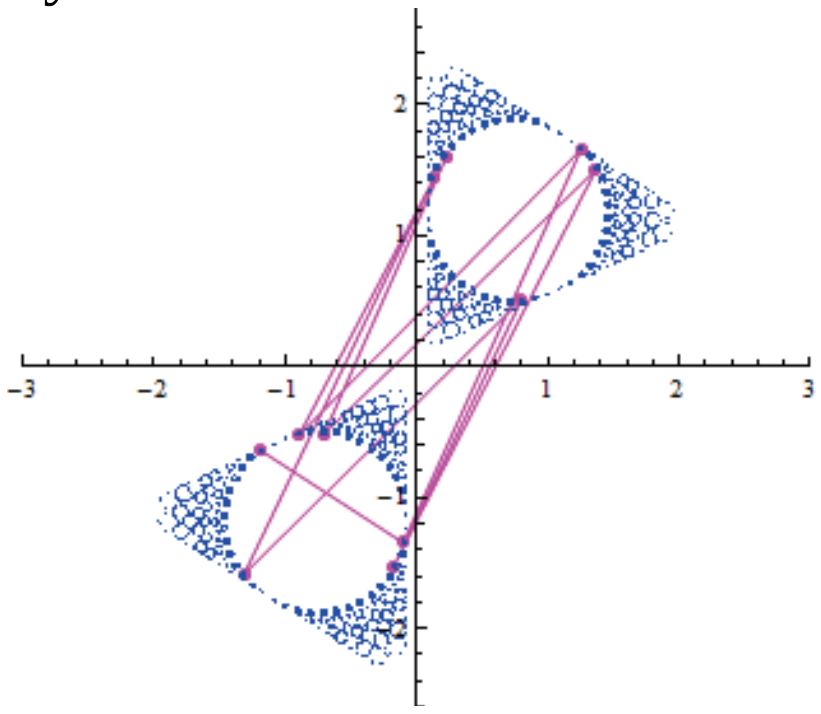

The dynamics of {1,0,0} are very complex and it is likely that this orbit is non-periodic. But the initial dynamics follow a simple pattern of rotation around the ellipses (which are circles above in the rectified plot), and gradually the points move outwards toward the vertices of the triangles. It is there that the 3-dimensional dynamics 'overflow' the boundaries of the unit cube and two more small planes are needed to contain the tips. (The worst case scenario of this overflow and cropping, explains why 4 planes suffice for any dynamics.)

The 'overflow' into planes 3 and 4 is just visible below on the projections and the corresponding 3D plot:

**Orbit=NestList[D3f,{1,0,0},30000];ROrbit=D3fToTr/@Orbit;P23=Drop[ROrbit,None,1]; ListPlot[P23]; Graphics3D[Point[Orbit]]**

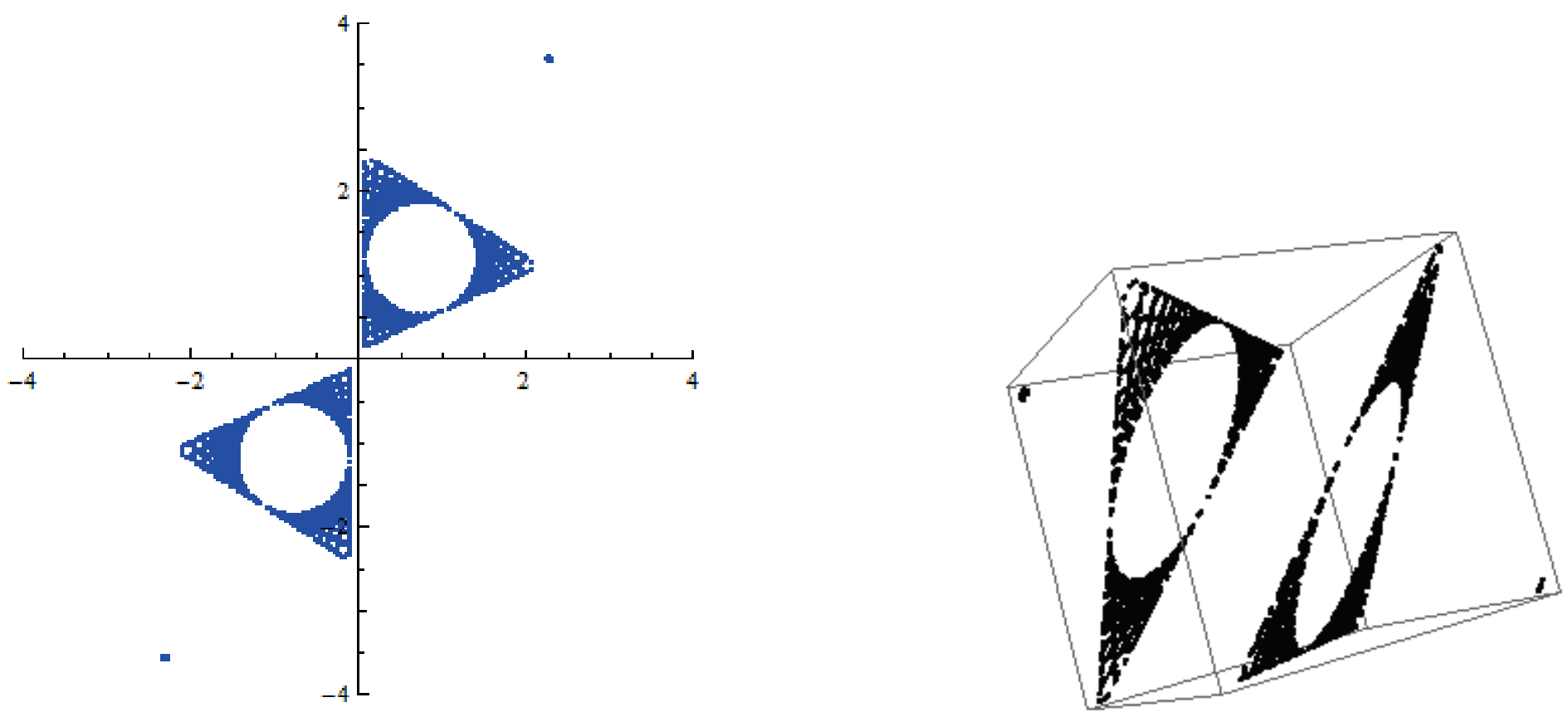

To get some idea of the complexity of the dynamics near the vertices, the plot below is based on 20 million iterations of {1,0,0} which are cropped to the pink region near the origin. Points are colored based on the order in which they were generated

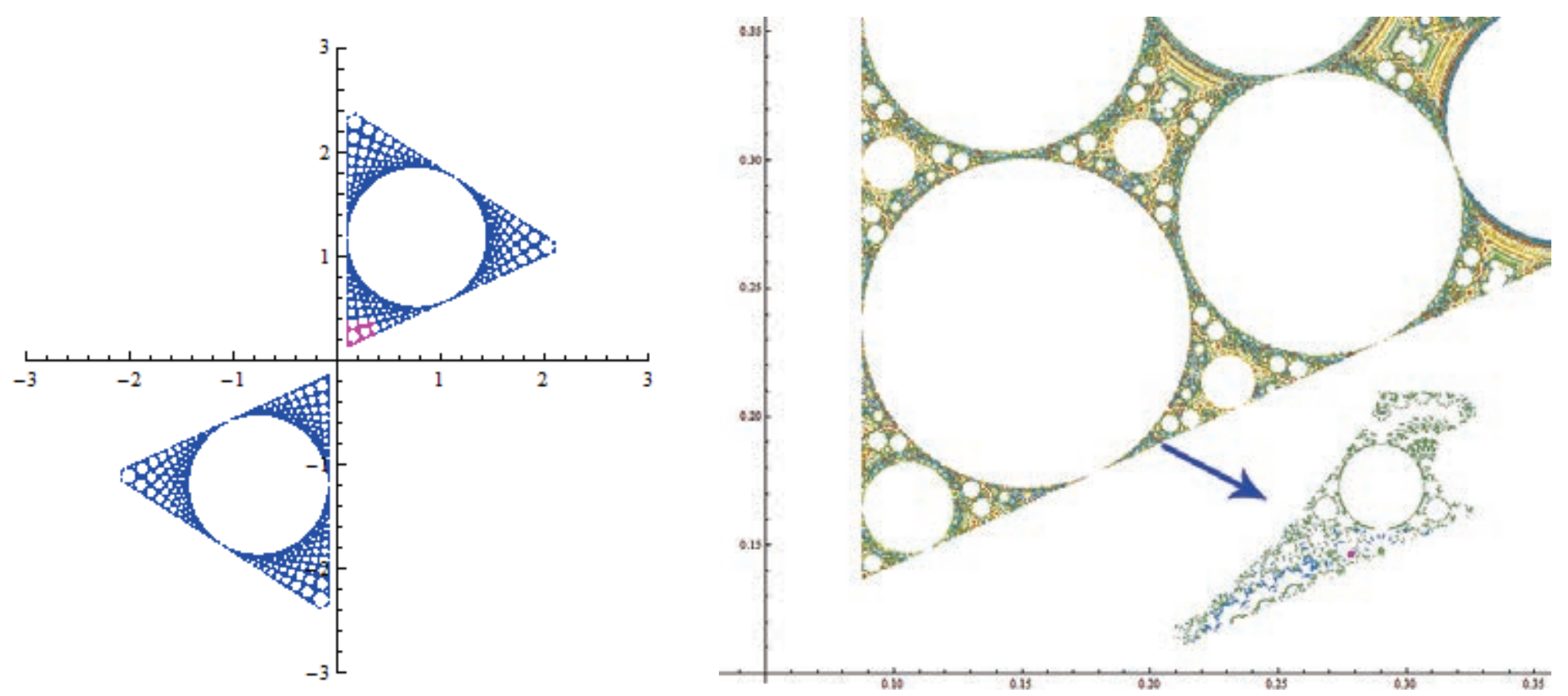

This appears to be a fractal image with structure on all scales.To see what a 'typical' s-sequence looks like for an initial point in this region we will advance the orbit of {1,0,0} to q0 = Orbit[[1098996]]. The red dot above is the projection of D3fToTr[q0]:{0.205138, 0.191607}. We will use q0 as the initial point in a new orbit and track the corresponding s-sequence, first for 20 iterations and then 4000. This is easy to code since we just need the first coordinate of each point in T-space (S1 below). S2 isolates this coordinate and S3 take differences and scales.

**gap = N[2/Sqrt[2+a^2]]; S1 = T1.Transpose[NestList[D3f, q0,20]];S2=Flatten[Drop[S1, -2, None]]; S3=Round[Differences[S2]/gap]** = {1,-1,0,1,-1,1,-1,0,1,-1,1,-1,1,-1,1,0,-1,1,-1,1}

For longer sequences the Tally command will keep track of the jumps. The larger jumps begin after 3234 iterations

**S1 = T1.Transpose[NestList[D3f, q0, 4000]]; S2=Flatten[Drop[S1, -2 ,None]]; S3=Round[Differences[S2]/gap]; Tally[S3]** = {1,1762},{-1,1761},{0,420},{-2,5},{3,25}, {-3,24},{2,3}

Below is the tally for 2 million iterations with initial point {1,0,0}
{{-1,920469},{1,920468},{0,152415},{-2,648},{3,2676},{-3,2676},{2,648}}

This is the type of tally we would expect from a periodic orbit or from a non-periodic orbit which is 'quasi-periodic' in the sense that it is the limit of periodic orbits (as in the regular pentagon).

For the Df map, a given value for the angular parameters *a*, determined the dynamics- which essential take place on 2 'planes'. With D3f, it comes as no surprise that a given value of *a* determines an infinite range of possible dynamics. For a fixed *a* value, every u value will yield a different set of planes. For $|a| < 1$, there are 2 or 3 planes and otherwise 3 or 4. For a given u, Chua defines $I^2_u$ to be the union of these planes. Two sets of planes $I^2_u$ and $I^2_v$ cannot intersect unless u =v. Since the D3f map is invariant on each $I^2_u$ it is possible to use information about 2-dimensional dynamics to analyze the dynamics of D3f. But the diversity of possible dynamics is staggering.

If we were to keep *a* fixed in the example above and change the initial point, even the smallest change could lead to totally different dynamics. There is no simple measure of 'distance' which will tell us which points are 'close to' {1,0,0} in D3f space.

For the Tangent map with a parameter such as w = 2*Pi/14, we are left with an infinite array of possible $I^2_u$ dynamics to consider. Hopefully one of these will mimic the 'perfect' dynamics of the Df map.

**Section 3: Dynamics of an Analog to Digital Converter**

Investigators such as Orla Feely [F]. have found 'qualitatively similar' behavior between Df dynamics with those of a second-order bandpass Delta-Sigma modulator. These modulators are used in Analog to Digital conversion. A typical block diagram for a (first order) ADC is show below

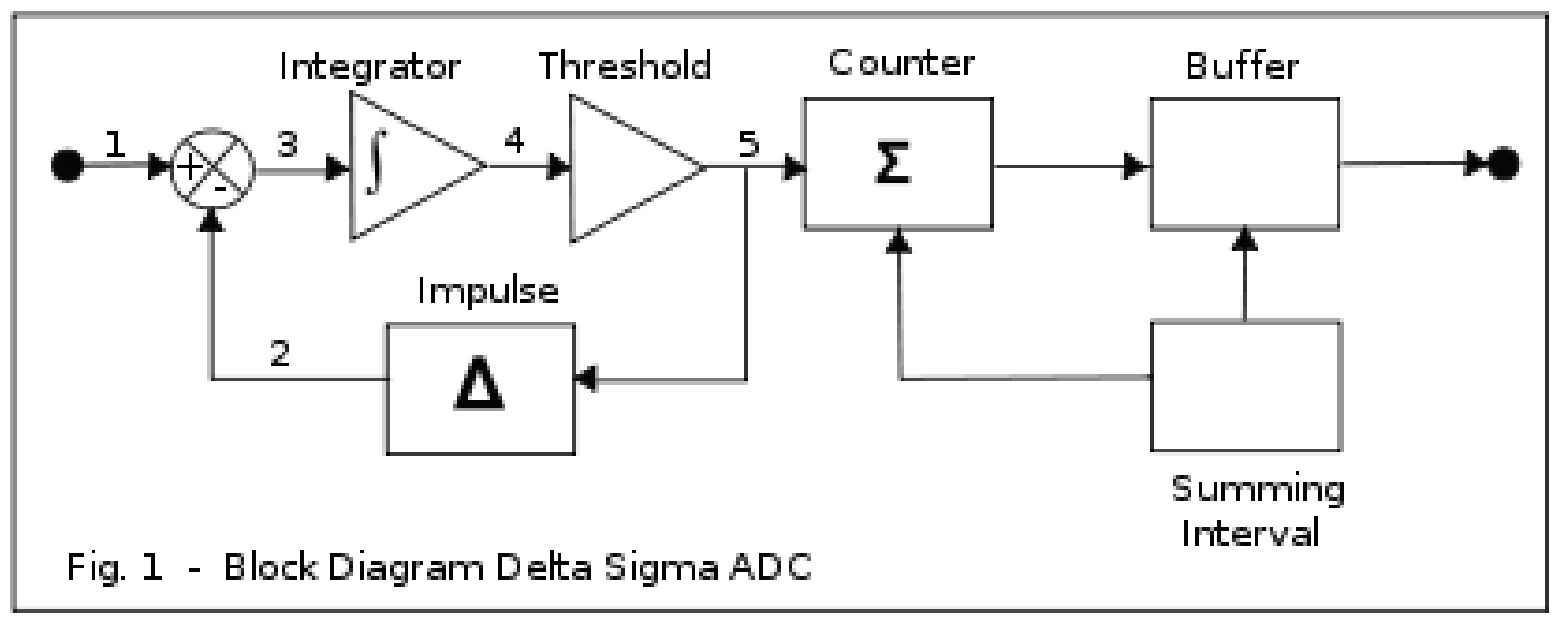

Fig. 1 - Block Diagram Delta Sigma ADC

The input is analog voltage and the output is a digital pulse stream with constant amplitude. The interval between pulses is determined by the primary feedback loop. A low input voltage produces longer intervals between pulses (and hence a smaller count in the Sigma counter). The Delta pulse zeros out the integrator to get ready for the next integration. The Sigma counter, counts the pulses over a given summing interval dt.

The first order ADC has stable behavior for all inputs which are in range, but this is not the case for the second -order ADC shown below. It has two integrators and hence 2 feedback loops.

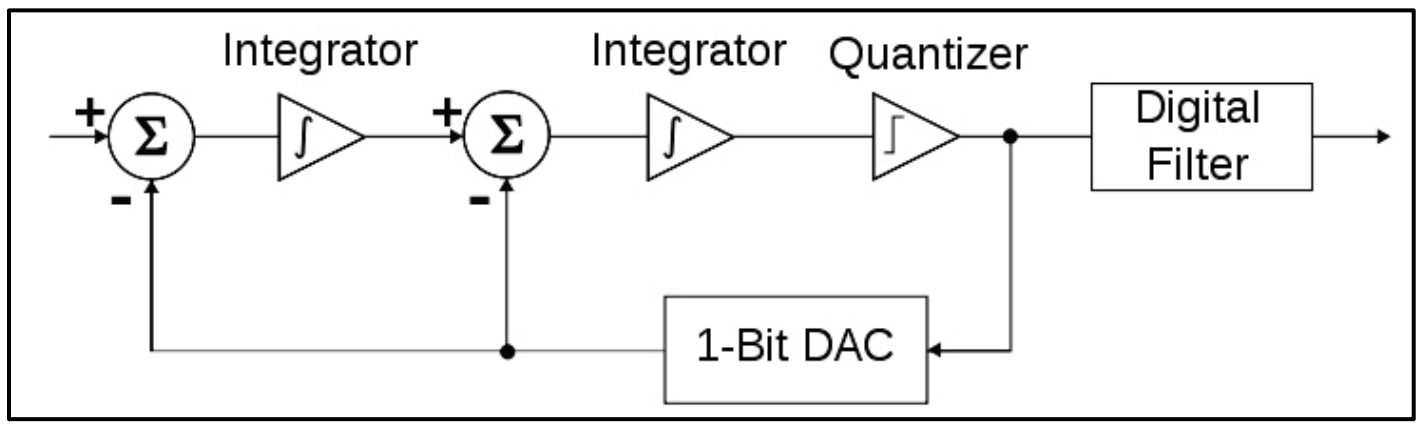

This type of negative feedback is common in biology, economics as well as mechanical systems. In biology it is known as homeostasis and in mechanical systems it may be an attractor or an equilibrium state. In chemistry, it is not unusual for one set of chemical signals to oppose another. Often these opposing forces are non-linear.

The non-linearity in the ADC is due to the one-bit DAC which can be modeled with a sign function:

$$\operatorname{sgn} x = \begin{cases} -1 & if \;\; x < 0 \\ 0 & if \;\; x = 0 \\ 1 & if \;\; x < 0 \end{cases}$$

The input-output equations for the ADC as formulated by Feeley and Fitzgerald, has the form of a second order difference equation:

$$y_{k+2} = 2r\cos\theta\, y_{k+1} - r^2 y_k + 2r\cos\theta(x_{k+1} - \operatorname{sgn} y_{k+1}) - r^2(x_k - \operatorname{sgn} y_k)$$

where x is the input, y is the output, r is the 'gain' and θ is a parameter that is chosen to center the filter at the desired sampling frequency.

In numerical studies we will set the gain r = 1, and to obtain self-sustaining behavior, set x = 0. This reduces the above difference equation to:

$$y_{k+2} = 2\cos\theta(y_{k+1} - \text{sgn } y_{k+1}) - (y_k - \text{sgn } y_k)$$

Set $x_1 = y_k$ and $x_2 = y_{k+1}$ then at each time tick: $x_1 \rightarrow x_2$ and $x_2 \rightarrow 2\cos\theta(x_2 - \text{sgn}x_2) - (x_1 - \text{sgn}x_1)$:
In matrix form:

$$\begin{bmatrix} x_{k+1} \\ y_{k+1} \end{bmatrix} = \begin{bmatrix} 0 & 1 \\ -1 & 2\cos\theta \end{bmatrix} \begin{bmatrix} x_k \\ y_k \end{bmatrix} + \begin{bmatrix} 0 \\ \text{sgn } x_k - 2\cos\theta \, \text{sgn } y_k \end{bmatrix}$$

Note that the 2 by 2 linear matrix A above is identical to the Df matrix with b = -1 and *a* = 2Cosθ. A is conjugate to a simple rotation using DfToTs from Section 1.

In Mathematica**: Adc[{x_,y_}]:={y, 2*Cos[w]*(y - Sign[y]) - (x - Sign[x])};**

**Example**: a = 2Cos[w] = -.8: **Orbit = NestList[Adc, {.001, 0}, 100000];**

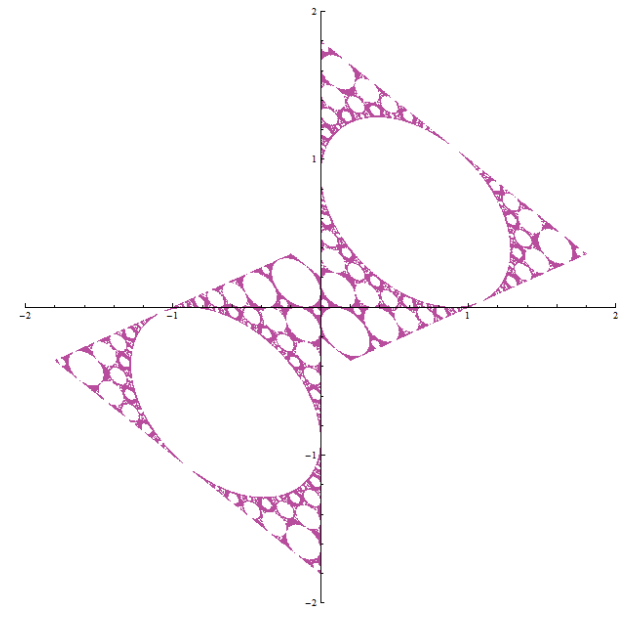

To track the regular pentagon set w = N[2*Pi/5] and again pick an initial point close to the origin: **q1 = {-.001,.001}**; **Orbit = NestList[Adc, q1, 100000]**; Rectify these points using **TrOrbit = DfToTs/@Orbit**

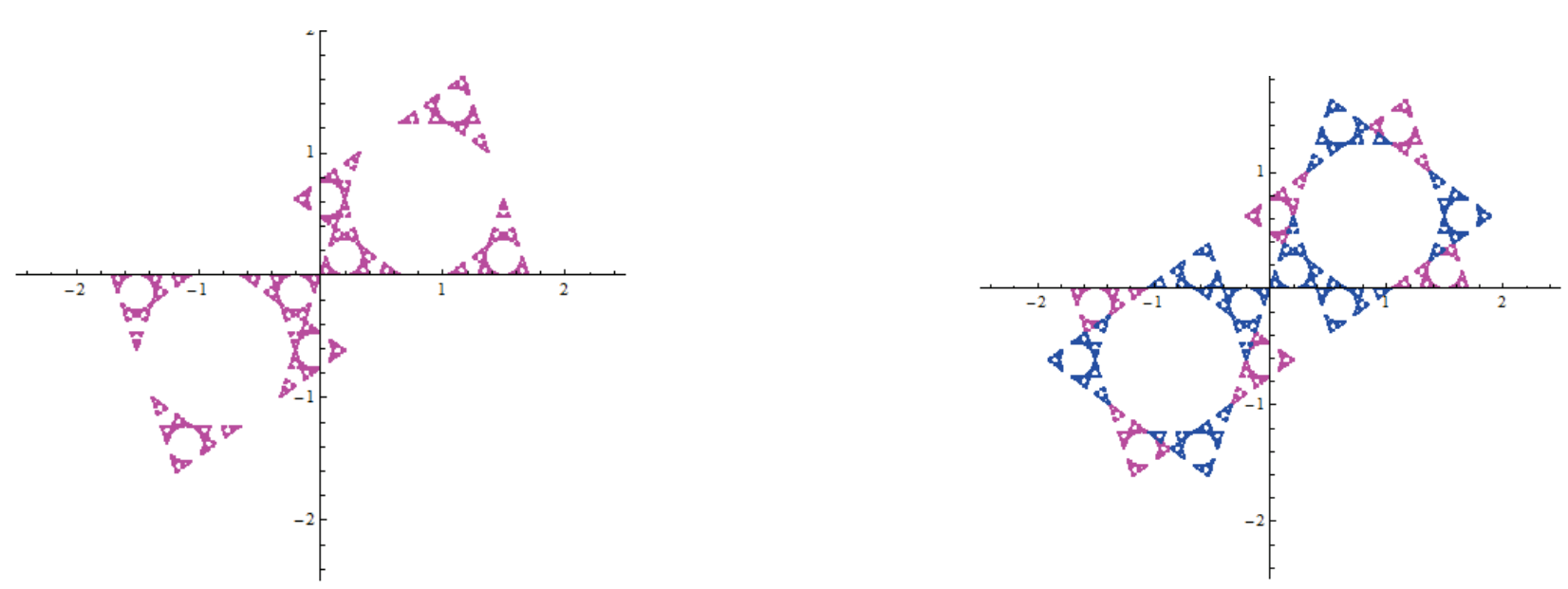

We can fill in the gaps by taking the union of this map with its reflection obtained by swapping x & y. This is shown on the right above in blue overlaid on top of the original map. Even though Adc is globally defined, it is only invertible on two trapezoids which are a function of θ. These trapezoids can be see above for θ= ArcCos[-.4] ≈ 1.98231 and θ = 2π/5, but one condition for invertibility is that θ ≥ π/3, so N = 5 is the last regular polygon that can be modeled. However we will see in Section 6 that a Harmonic Kicked Oscillator yields a map which is locally conjugate to Adc and invertible for almost all θ. This map continues the chain of conjugacies back to the Tangent Map

## Section 4: Hamiltonian Dynamics

From the perspective of Hamiltonian dynamics, the digital filter map and ADC maps are examples of 'kicked' Hamiltonians.

The Hamiltonian of a physical system is a function of (generalized) position $\vec{q}$, and (generalized) momentum $\vec{p}$. In some cases it may change in time, so it is writen as $H(q_i(t), p_i(t), t)$ where i = 1,....N are the degrees of freedom of the physical system. The dynamics of a Hamiltonian system are given by:

$$q_i' = \frac{\partial H}{\partial p_i}, \quad p_i' = -\frac{\partial H}{\partial q_i}$$

So a Hamiltonian system is a system of differential equations which are related by H. Not every system of differential equations has a corresponding Hamiltonian. When H does not depend explicitly on time, the state of the system at any given time is a point ($p_i$,$q_i$) in 2N-dimensional phase space. The motion of this point is determined by the equations above. Starting with an initial set of points, the volume will be preserved over time because the equations above imply that H is conserved: $\frac{dH}{dt} = \frac{\partial H}{\partial p_i}\frac{dp_i}{dt} + \frac{\partial H}{\partial q_i}\frac{dq_i}{dt} = 0$

When N = 2, there are two position coordinates $\vec{q}$ ={$q_1$,$q_2$} and two momentum coordinates $\vec{p}$ = {$q_1$,$q_2$} so the system is 4-dimensional. For example if the differential equations describe the motion of 2 planets around the sun, then the system motion would be periodic in {$q_1$,$p_1$} and {$q_2$,$p_2$} so there would be two frequencies. The Hamiltonian would not be constant, but it would be periodic so these two phase plots would combine together to yield motion on a 2-torus

Below is the case of a 2-year period vs. a 5-year period which is similar to Saturn (29.447 years) and Jupiter (11.861 years), so 5*11.861 ≈ 2*29.447. This is called a 'rational' torus. If the periods were incommensurable, like Earth and Jupiter, the windings would fill the torus and it would be called an irrational torus. The corresponding motion is called a quasi-periodic orbit. These are the most stable orbits because the gravitational interactions have no 'resonances'.

**gr1=ParametricPlot3D[{(2+Cos[u])*Cos[v],(2+Cos[u])*Sin[v],Sin[u]},{v,0,2*Pi},{u,0,2*Pi}, Mesh->False, Boxed->False, Axes->False];** (*this is the solid torus with radii 2 and 1*)

**gr2=ParametricPlot3D[{(2+Cos[5 t]) Cos[2 t],(2+Cos[5 t]) Sin[2 t],Sin[5 t]},{t,0,2 Pi}];**
**Show[gr1,gr2]**

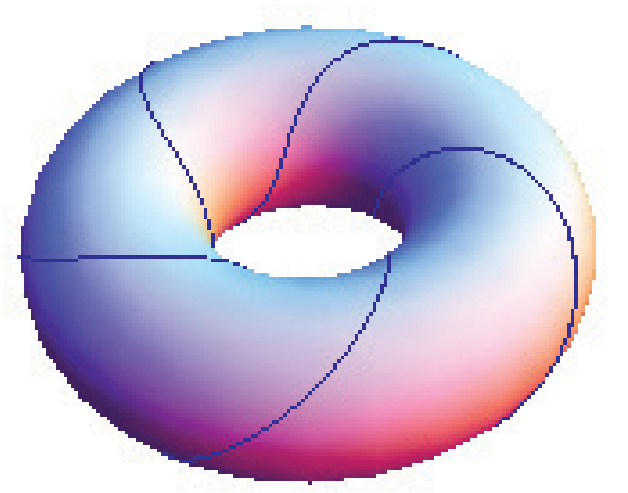

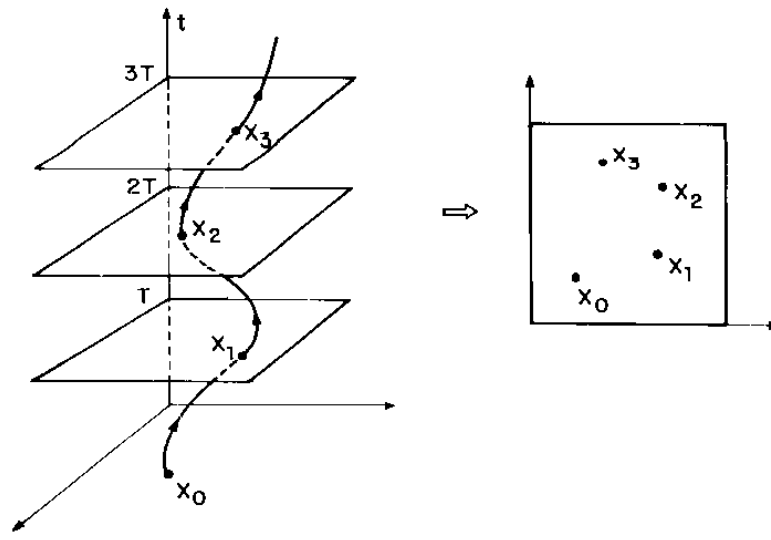

To reduce this motion to a 2-dimensional mapping, choose a reference plane. For example slice the torus at a fixed longitude $q_0$ and 'strobe' the longer 5 year period with the shorter 2-year period by setting $\delta t = 2$ years.

These surfaces of section are also called Poincare return maps. They can be used to help uncover periodic motion or fixed points. Since the system equations are differential equations, continuity would imply that orbits in the vicinity of a fixed point will behave in a predictable fashion based on the Jacobian. Because S is symplectic, the Jacobian at any fixed point will have eigenvalues of the form $\{\lambda, 1/\lambda\}$. There are only three types of behavior possible near a fixed point:

(i) $\lambda$ positive and $\lambda > 1$ ($1/\lambda < 1$) This is an unstable hyperbolic 'saddle' point where the motion expands in one direction and contracts in another.
(ii) $\lambda$ negative and $|\lambda| > 1$, ($1/|\lambda| < 1$). Similar to (i) above but displacements in each direction vary from positive to negative on each iteration.
(iii) $\lambda$ complex - in which case $\lambda$ and $1/\lambda$ must have zero real part so they have the form $\lambda = e^{i\theta}$ and $1/\lambda = e^{-i\theta}$. This is called a 'center' or a neutral fixed point. Nearby points will rotate by $\theta$ on each iteration.

Case (iii) is the only case where stability is not obvious, because this 'marginally stable' behavior may be very sensitive to small perturbations. This is the heart of the KAM Theorem: How will the marginally stable fixed points of a 'well-behaved' Hamiltonian system react when subjected to periodic perturbations ? (i.e. when gravity between Jupiter and Saturn is turned on). Will a finite measure of the quasi-periodic orbits survive ? Is the solar system stable ?

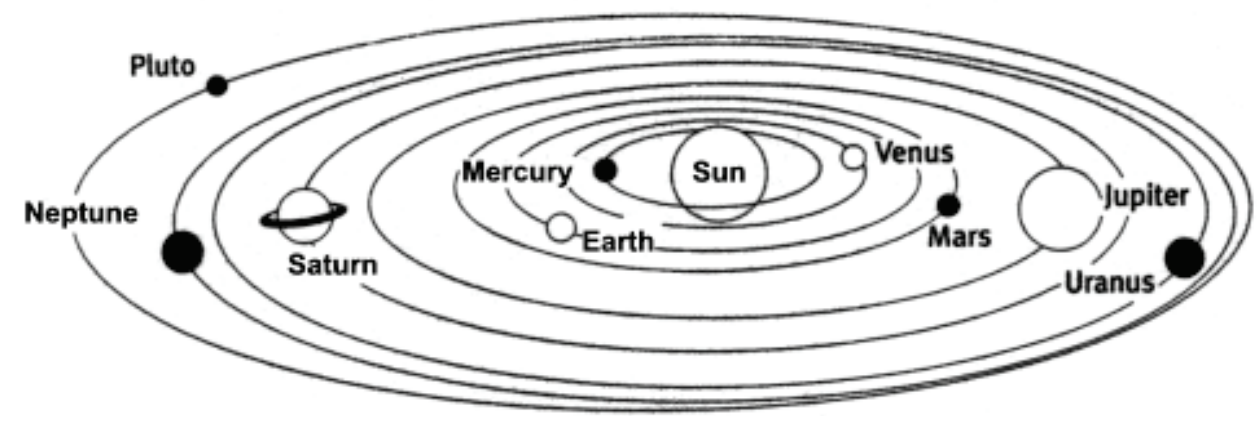

Note that it is not possible for solutions to spiral into or away from a fixed point and also it is not possible to have expansion without contraction. These are consequences of the symplectic property of the Hamiltonian.

A 2-degree of freedom Hamiltonian such as the one above would be classified as 'well-behaved' or integrable, if there are 2 integrals of motion. A function $f(\vec{p}, \vec{q})$ is an Integrable of Motion if it satisfies the 'Poisson bracket' (commutator) relationship.

$$[f, \mathrm{H}] = \frac{\partial f}{\partial p_i}\frac{\partial p_i}{\partial t} - \frac{\partial f}{\partial q_i}\frac{\partial q_i}{\partial t} = \frac{df}{dt} = 0$$

Therefore an integral of motion is also a constant of motion. The system energy $\mathrm{H}(\vec{p}, \vec{q})$ is always a constant of motion. Note that the torus maps shown above have non-Eulidean metric, just like the raw Df maps. We can obtain the integrals of motion by performing a canonical change of variables to convert the map S into 'normal' Euclidean form so that ellipses become circles. The new variables are called 'action-angle' variables: $(\vec{p}, \vec{q}) \rightarrow (\vec{I}, \vec{\theta})$ .In the new surface of section map, θ will represent a simple rotation. The I's are the two constants of motion: $I_j = \oint_{C_j} \vec{p}.d\vec{q}$ where the integral is around the jth axis of the torus.

**Example**: Two coupled oscillators with unit mass m. $H = \frac{1}{2}(p_1^2 + \omega_1^2 q_1^2) + \frac{1}{2}(p_2^2 + \omega_2^2 q_2^2)$ where $\omega_1 = \sqrt{k_1}$ and $\omega_2 = \sqrt{k_2}$ . Transformation to action-angle variables yields $p_1 = \sqrt{2I_1\omega_1}\cos\theta_1$ $q_1 = \sqrt{2I_1\omega_1}\sin\theta_1$ and $p_2 = \sqrt{2I_2\omega_2}\cos\theta_2$ $q_2 = \sqrt{2I_2\omega_2}\sin\theta_2$ so the new Hamiltonian is $H = \omega_1 I_1 + \omega_2 I_2$ and $I_1$ and $I_2$ are the integrals of motion.

The motion is the same as 2 uncoupled oscillations: $\theta_1(t) = \omega_1 t + \theta_1(0)$ and $\theta_2(t) = \omega_2 t + \theta_2(0)$ so the periods are $2\pi/\omega_1$ and $2\pi/\omega_2$ and the motion does not depend on the central spring.

The KAM Theorem starts with an integrable Hamiltonian and asks what happens when it is perturbed. For the N = 2 case of Saturn and Jupiter this amounts to 'turning on' the gravitational attraction between them. In this case the answer to the stability question depends largely on the ratio of the periods, A ratio such as 2/5 could be dangerous over long time periods because it means that Saturn will feel systematic effects which could build up over time. It is not clear whether the canonical change of variables will remain meaningful under perturbations and the new system will generally not be Integrable.

The early solar system was a turbulent place and recent studies show that there may have been an even more dangerous 2:1 resonance between Saturn and Jupiter in the early history of the solar system. This resonance might have caused Neptune and Uranus to migrate outwards and this in turn would have created havoc among the millions of planetesimals in that region.

If there is observable chaotic motion within our solar system currently, it is probably in the orbits of 'irregular' satellites. The first irregular satellite that has been imaged close-up and studied in detail is Saturn's Phoebe. Recently there have been a number of new discoveries of irregular satellites orbiting the four gas giants Jupiter, Saturn, Uranus and Neptune. These discoveries were made possible by attaching new wide-field CCD cameras to large telescopes.

The Earth recently captured a small irregular asteroid. It was called 2006 $RH_{120}$ and it was about 15 meters wide. It only survived for about 4 orbits before returning to its orbit around the sun. The dates were September 2006 to June 2007. Latest news is that there is another companion to the Earth in an irregular orbit. It is known as 2010 $SO_{16}$.

To see what happens to an integrable Hamiltonian when it is perturbed, we will look at the case of one degree of freedom: $H(\theta,p) = \omega p^2$ so the momentum p is constant. The corresponding

Poincare return map is $\theta_{k+1} = \theta_k + p$ so points rotate by a fixed amount p.When the perturbation is of the form $K\sin\theta_k$, this is called the Standard Map.

**The Standard Map**

$$\begin{bmatrix} x_{k+1} \\ y_{k+1} \end{bmatrix} = \begin{bmatrix} x_k + y_k + K\sin x_k \\ y_k + K\sin x_k \end{bmatrix}$$

When x plays the role of the angular variable θ and y plays the role of momentum p. The determinant of the Jacobian is 1, so it is area preserving. Since the perturbation is periodic modulo 2*Pi, it is common to scale x and K by 2*Pi and plot x Mod1. Following normal convention we will plot both x and y Mod 1.

**Std[{x_, y_}] := Mod[{y + (K/(2*Pi))*Sin[2*Pi*x] + x, y + K/((2*Pi))*Sin[2*Pi*x]}, 1];**

Example: **K =0.971635406; gamma = (Sqrt[5]-1)/2;**
**Orbit1= NestList[Std, {0,gamma}, 1000]; Orbit2= NestList[Std, {0,1-gamma}, 1000];**
To get a range of initial conditions passing through the center of the plot:
**VerticalScan = Flatten[Table[NestList[Std, {.5, y}, 1000], {y, 0, 1, .01}], 1];**
**Graphics[{AbsolutePointSize[1.0],Blue,Point[VerticalScan],Magenta, Point[Orbit1], Point[Orbit2]},Axes->True]**

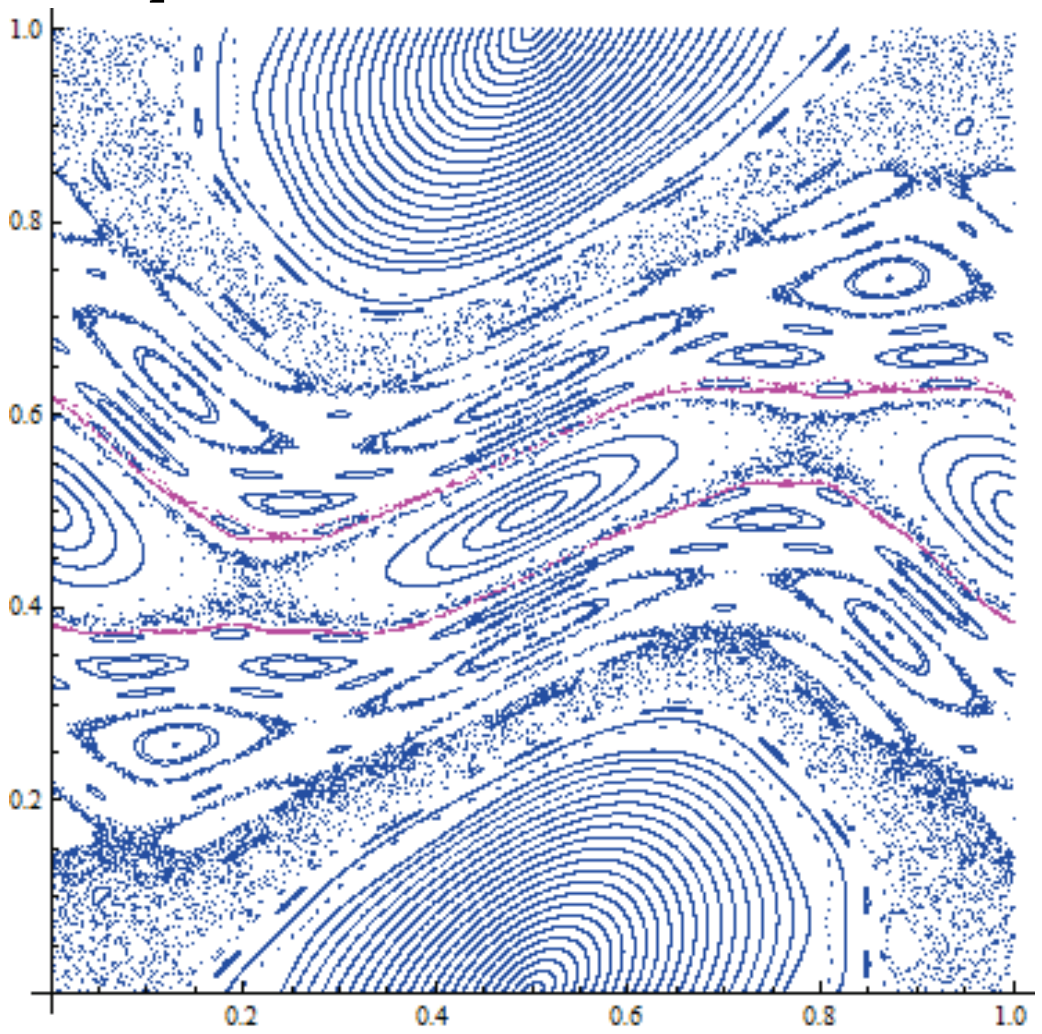

This mapping is a paradigm for Hamiltonian chaos in the sense that it is an area preserving map with divided phase space where 'integrable islands of stability are surround by a chaotic component'. The two magenta orbits are remnants of the last surviving invariant torus - with highly irrational winding number equal to the Golden Mean. (See Section 8). As K increases past 0.971635406 and these curves break down, there is no further impediment to large scale chaos in the central region.

The K = 0 case corresponds to an integrable Hamiltonian such as $H(\theta,p) = \omega p^2$ so the angular momentum p is constant and the Standard Map reduces to a circle map (twist map) $x_{k+1} = x_k + y$ with winding number y. The circles (tori) come in two varieties depending on whether $y/2\pi$ is rational or irrational. When K > 0, The Standard Map simulates what happens to H when it is perturbed by periodic 'kicks' of the form $K \sin x$ . So the Standard Map is a Poincare cross section of perturbed twist map. We will derive the equations below using a kicked rotor.

The fixed points are at {0,0} and {1/2,0}. The Jacobian is $J = \begin{pmatrix} \frac{\partial f_1}{dx} & \frac{\partial f_1}{dy} \\ \frac{\partial f_2}{dx} & \frac{\partial f_2}{dy} \end{pmatrix}$

$= \begin{bmatrix} 1 & 1 \\ K\cos 2\pi x & 1 + K\cos 2\pi x \end{bmatrix}$ At {0,0}, $J = \begin{bmatrix} 1 & 1 \\ K & 1+K \end{bmatrix}$ so this is an unstable saddle point as we can see above. At {1/2,0} $J = \begin{bmatrix} 1 & 1 \\ -K & 1-K \end{bmatrix}$ and this is a center for $K \in (0,4)$

Therefore for K in this range, the motion in the vicinity of {1/2,0} is locally conjugate to a rotation. These rotations are the cross sections of the 'tori' from the Hamiltonian and we know that the curves with irrational winding numbers (like the Golden Mean curve) are most likely to survive as the perturbations are increased. At this resolution it is hard to tell which curves are rational and which are irrational. Below is an enlargement of a central region showing the nesting of rational and irrational curves

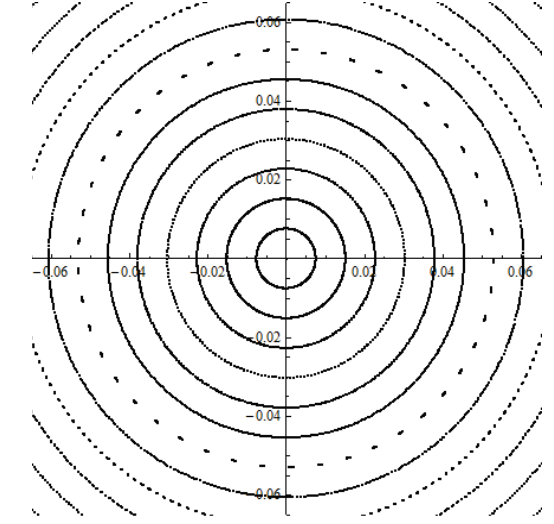

Moving outwards from the center means that the conjugacy may no longer be valid. For small K, the set of surviving curves can be parameterized by a Cantor set, so they form a continuum, but the regions between these irrational curves contains rational rotations with periodic orbits such as the period 46 orbit above. Under perturbations, the KAM Theorem says that tori which are 'sufficiently' irrational will survive, while most 'rational' tori will break up into the resonant islands that we see in the Standard map. This process repeats itself at all scales as shown in the drawing below. The '×'s' are unstable 'saddle points' which form between the islands. Each island represents a locally 'fixed' point because periodic points of a mapping S are fixed points of $S^k$.

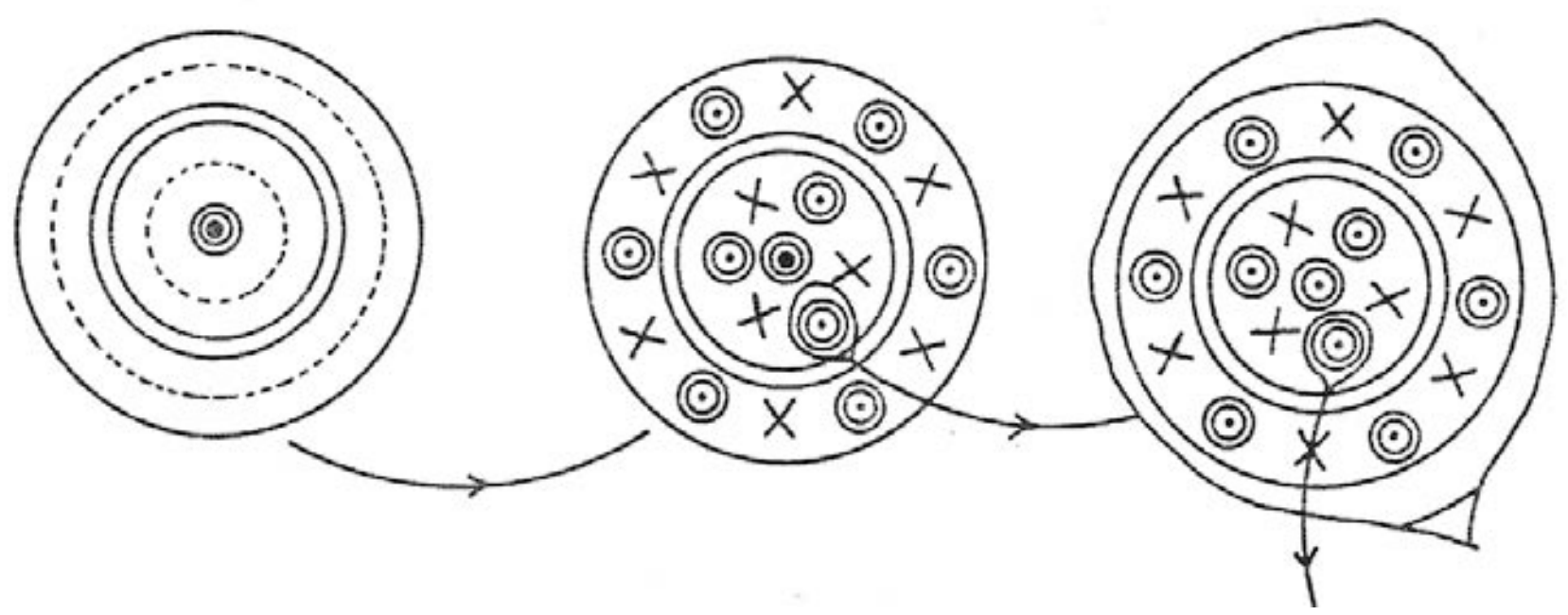

For small perturbations, a finite measure of tori will survive and this guarantees stability. The invariant curves that do survive, partition the phase space and limit the chaotic diffusion. The dark regions in the Standard map plot are the chaotic debris left over from the breakdown of rational torus. If an asteroid had these initial conditions, it might be flung from its orbit.

**Deriving the Standard Map from a Kicked Rotor**

Consider a pendulum allowed to rotate in zero gravity, and subjected to periodic kicks of magnitude K at time intervals $\tau$ as shown here.

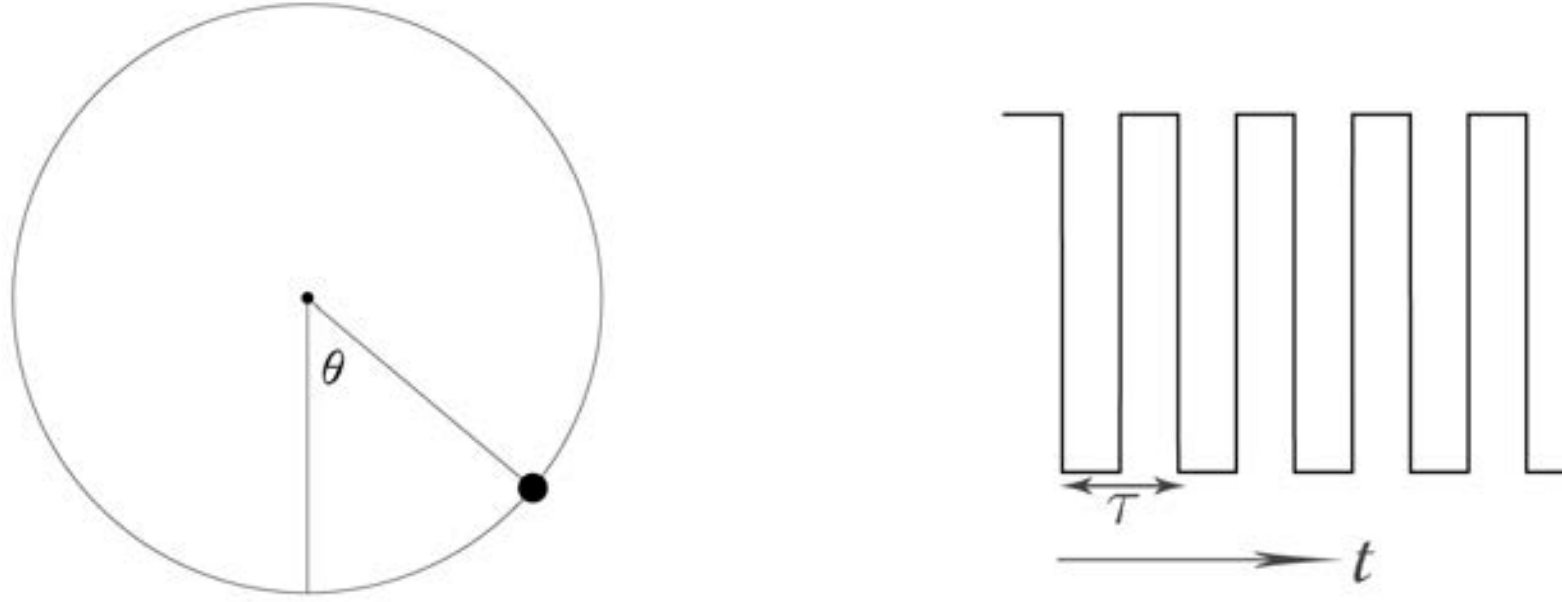

If this were an asteroid in orbit around the sun, the 'kicks' might be gravitational interaction with Jupiter. Another important scenario is a cyclotron where charged particles are accelerated to high speed in a magnetic field created by an alternating electric field.

The Hamiltonian is time dependent: $H(\theta, p, t) = \dfrac{p^2}{2I} + K\cos\theta \sum_n \delta(\dfrac{t}{\tau} - n)$

where I is the moment of inertia of the pendulum, which we will scale to be 1, and $\delta$ is the Dirac Delta function which yields pulses only when the input is 0 and this occurs at time intervals $\tau$ , so

$\delta(t) = K\sum_n \delta(\dfrac{t}{\tau} - n)$ is the periodic pulse stream shown above.

From the Hamiltonian,we can see that the pulse stream alters the Potential Energy at periodic intervals and this shows up in the system equations as a change in Kinetic Energy:

The equations of motion are: $\dfrac{dp}{dt} = K\sin\theta \sum_n \delta(\dfrac{t}{\tau} - n)$ and $\dfrac{d\theta}{dt} = p$

So the rotor receives a periodic torque of magnitude $K\sin\theta$ at time intervals $\Delta t = \tau$. For the discrete equations, we can scale the intervals to get $\tau = 1$. The effect of the kth kick is to update p based on the current value of $\theta$, so $p_{k+1} = p_k + K\sin\theta_k$. Between kicks, the motion is force free so p is constant and the second equation above says that $\theta$ is updated by $p\Delta t = p$ , so $\theta_{k+1} = \theta_k + p_{k+1}$.This gives the following equations which describe the system just after the kth kick.

$$\begin{bmatrix} \theta_{k+1} \\ p_{k+1} \end{bmatrix} = \begin{bmatrix} \theta_k + p_{k+1} \\ p_k + K\sin\theta_k \end{bmatrix} = \begin{bmatrix} \theta_k + p_k + K\sin\theta_k \\ p_k + K\sin\theta_k \end{bmatrix}$$

These equations are identical to the Standard map, and this tells us that the kicked rotor would have very interesting dynamics but zero-gravity is hard to simulate in a lab. We can see this scenario being played out in celestial mechanics, but the planetary time scale is typically very long. On a shorter scale we can the track the orbits of satellites or asteroids where the periodic perturbations may be due to the irregular orbits or irregular shapes. Particle accelerators and fusion reactors live at a much faster pace, and instabilities are a major issue in both cases.

**The Delta-Kicked Rotor in Quantum Physics**

In quantum physics, momentum P and position X are no longer independent. They are related by the Poisson commutation relationship $[X, P] = XP - PX = i\hbar$ where $\hbar = h/2\pi$ is the reduced Plank constant: $\hbar \approx 1.054571 \times 10^{-34}$ joule seconds It represents the proportionality between the momentum and quantum wavelength of a particle. In 1926 Werner Heisenberg realized that this commutation relationship implies that $\Delta X \Delta P \geq \hbar/2$ where $\Delta X$ and $\Delta$P are deviations in the measured values of position and momentum. This is known as the Heisenberg Uncertainty Principle. It says that it is not possible to know both of these quantities with high accuracy at the same time because decreasing one uncertainty, increases the other. This principle underlies all of quantum mechanics and has wide-ranging implications

For large aggregates of atoms, quantum effects can often be ignored and it is possible to apply the 'classical' laws of physics like Newton's laws of motion and the laws of electricity and magnetism. The dividing line between the classical and quantum worlds is very difficult to determine and recently there have been studies showing that quantum affects can be observed all around us if we look carefully. Often these effects are in conflict with the laws of classical physics.

No one doubts that the laws of quantum mechanics govern the behavior of particles on very small scales, but some of these laws seem very bizarre and they do not seem to apply to the macro world around us. The laws of classical physics should be based on these quantum laws, but no one knows how to do this.

Vlatok Vedral describes an example involving ordinary salt. On the left below are random orientation of molecules in a salt crystal.When a magnetic field is applied, the atoms have a tendency to align themselves with the field, but in fact they align themselves much more efficiently than the classical prediction. Physicists believe that this is due to quantum effects, which predict that at a certain scale and energy level, particles have a natural 'entanglement'

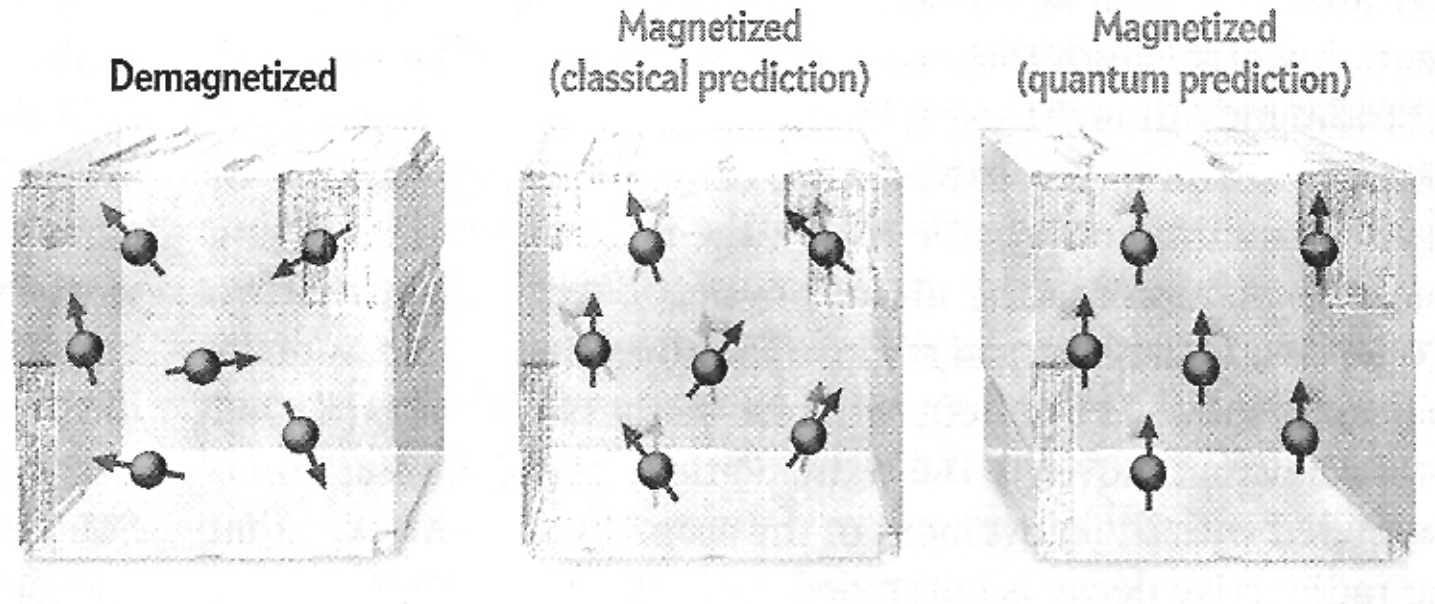

In Moore *et al* a group of five physicists at the University of Texas used ultra-cold salt atoms in an experiment which showed a similar 'entanglement'. The first step was to trap the atoms in a magneto-optical 'trap' so they are ultra-cold with minimal motion. About $10^5$ atoms were trapped and subjected to a modulated laser beam (the 'kicks'). The sodium atoms have a 'ground' state and an 'excited' state and the beam energy causes them to oscillate between states.

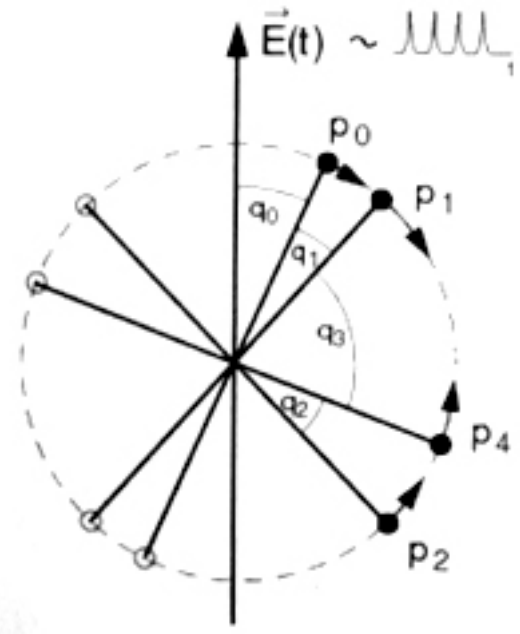

This 'dipole' has the same characteristics as a kicked rotor, but now the kicks are generated by an electric field E(t) that is pulsed with period τ so $E(t) = K\sum_n \delta(\frac{t}{\tau} - n)$.

The Hamiltonian has the form $H = \frac{p^2}{2I} + dE\cos q$ where d is the dipole moment, I is the inertia and q is the angle between the dipole and the field and p is the angular momentum.

If we compare this Hamiltonian with the kicked rotor above, they are virtually identical so the Standard Map can be used to describe the classic behavior of the sodium atoms. The key issue is how closely this 'classical' behavior matches the actual behavior.

Below is a plot the Standard Map in a perturbation range similar to that used for the Texas experiment. This is the 'large-scale chaos' regime with **K= 3.0** The blue points are from a scan to depth 500 and the pink points are 10,000 points in a 'randomly chosen' orbit. This orbit will eventually turn the plot pink except for the invariant curves.

```
Std[{x_,y_}]:=Mod[{y+(K/(2*Pi))*Sin[2*Pi*x]+x,y+K/((2*Pi))*Sin[2*Pi*x]},1];
Scan = Flatten[Table[NestList[Std, {.5, y}, 500], {y, 0, 1, .01}], 1];
```

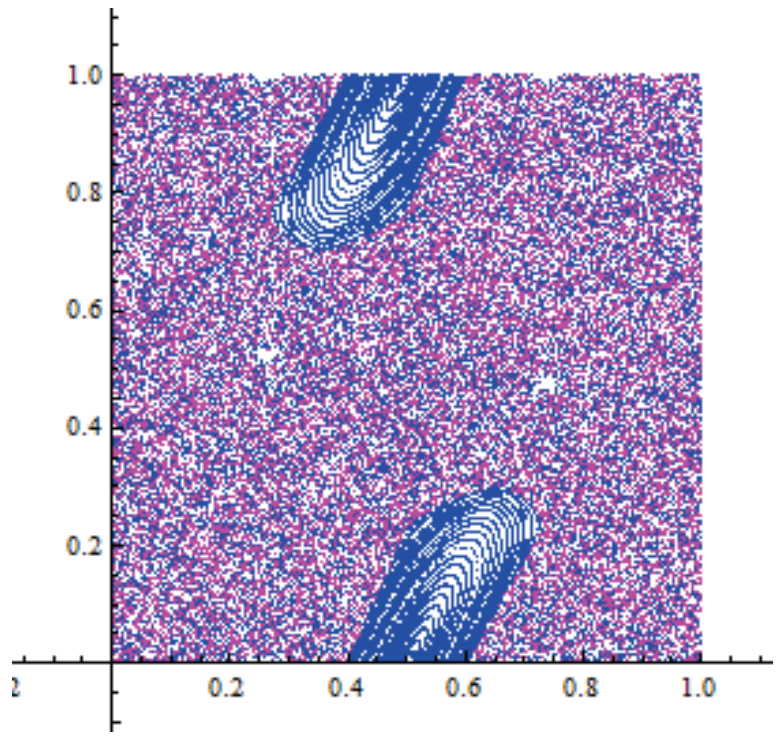

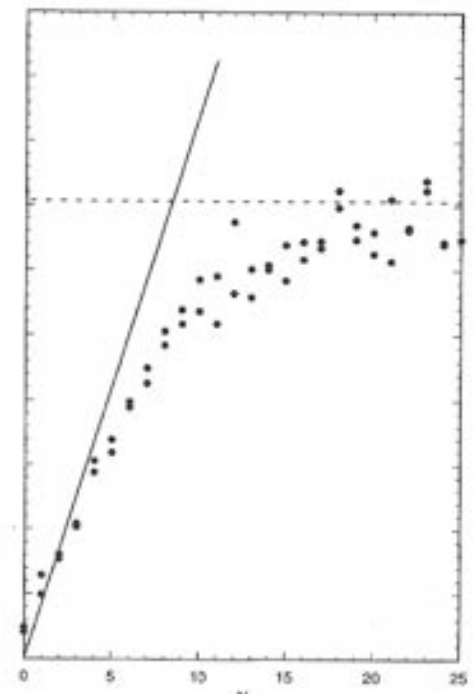

This 'sea-of-dots' actually contains small resonant islands on all scales but the overall system behavior is 'diffusive', meaning that individual orbits tend to fill the entire phase space. In the Texas University experiment, this is the type of 'classical' behavior observed initially. On the right above is a plot of their results. The solid line is the classical prediction of the energy gain as the sodium atoms are pulsed. This corresponds closely to behavior predicted by the delta-kicked rotor and illustrated above by the Standard Map. The 'x-axis' of their plot is the number of cycles of pulsing and 'refreezing'. After about 8 cycles, there was significant drop-off of energy. This drop-off is due to a quantum effect called 'localization', and it had been predicted by theory and by earlier experiments. This behavior appears to be related to the magnetic induced 'entangelment' found in the salt crystals.

As indicated earlier, the classical kicked rotor (and the resulting Standard Map) are regarded by physicists as paradigms for Hamiltonian chaos but this is the first time that this behavior was controlled in a laboratory so that it was possible to observe the transition between classical and quantum dynamics.

For any discrete mapping to be of value in quantum mechanics, it has allow for a consistent 'quantization'. The Standard Map allows researchers to impose this quantization, and the same is true for the kicked Harper Map which we will examine in Section 6.

**Section 5: The Sawtooth Standard Map**

The non-linearities with Df and D3f arose from discontinuities in the overflow function *f*[x_]:= Mod[x + 1,2] - 1 and the non-linearities of the Adc map arose from 'one-bit' discontinuities modeled with x-Sign[x] . These are both examples of sawtooth maps which in Mathematica are called SawtoothWave functions: SawtoothWave[x] is defined to be x-Floor[x]

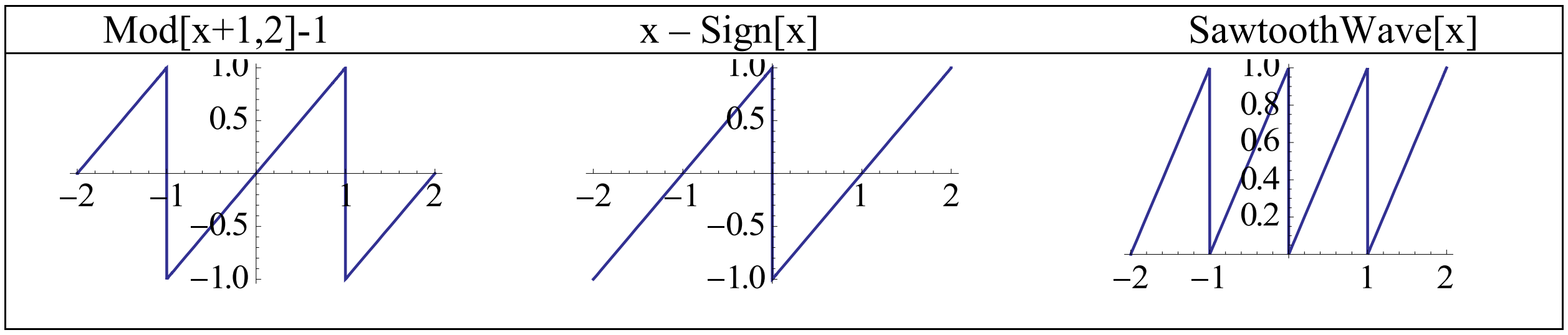

There have been a number of studies done with the Standard Map by replacing the 'kick' function $f(x) = \sin 2\pi x$ with either non-analytic or non-smooth functions. Usually the non-smooth functions are piecewise linear 'sawtooth' maps. For example V.V. Vecheslavov has studied the dynamics of the following piecewise version of $f(x) = \sin 2\pi x$

**Saw[x_]:=Which[ <=.25, 4*x, Abs[.5-x]<=.25, 4*(.5-x), Abs[x-1.0]<=.25, 4*(x-1)];**
**Plot[ Saw[x], {x,0,1}]**

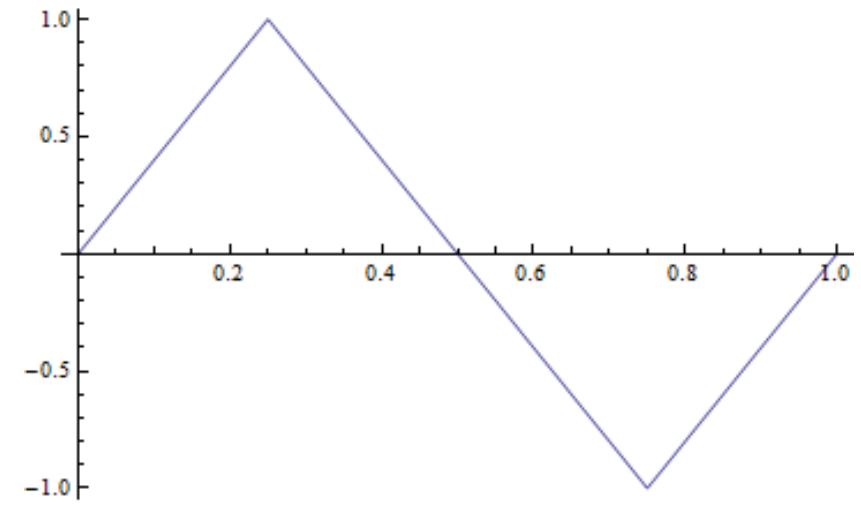

The corresponding 'SawtoothSine' Standard Map is:

$$\begin{bmatrix} p_{k+1} \\ x_{k+1} \end{bmatrix} = \begin{bmatrix} p_k + \mathrm{KSaw}(x_k) \\ x_k + p_{k+1} \end{bmatrix}$$

**SawSine[{x_, y_}] := {Mod[x + K*Saw[y], 1], Mod[x + K*Saw[y] + y, 1]};**

For some K values this map is similar to the Standard Map, but for most values of K the dynamics are quite different - especially near 'resonant' values such as K = 1/3 and K = 1/8. Below are three sample plots using K =2*Cos[2*Pi/14] - 2 , K = .332 and K = 1/2

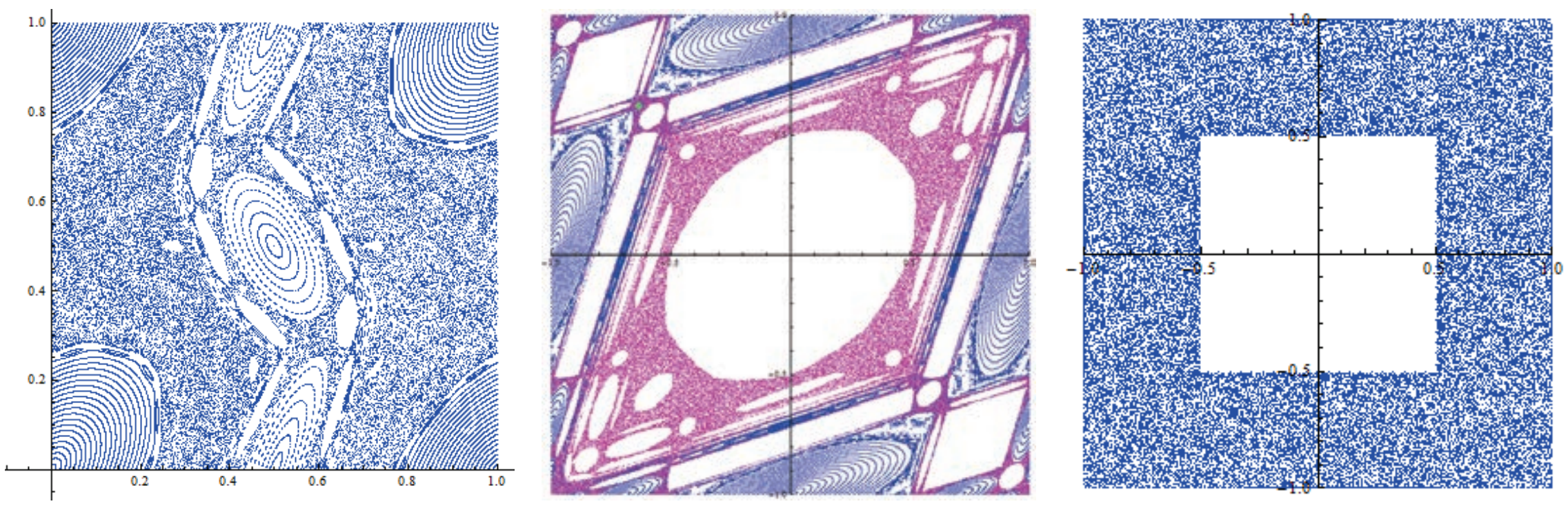

I.C. Percival and others have studied families of 'Sawtooth' maps obtained by linearizing the Standard Map about the stable fixed point at {1/2,0}. These map have the form:

$$\begin{bmatrix} p_{k+1} \\ x_{k+1} \end{bmatrix} = \begin{bmatrix} p_k + K(x_k - 1/2) \\ x_k + p_{k+1} \end{bmatrix}$$

Where only x is taken Mod 1 so it maps the cylinder $[0,1)\times\mathbb{R}$ to itself. The 'kick' function $f(x) = \sin 2\pi x$ of the Standard Map is replaced with $f(x) = (x - 1/2)$ Mod1:
**Plot[Mod[x-1/2,1],{x,0,1},Exclusions->None ]**

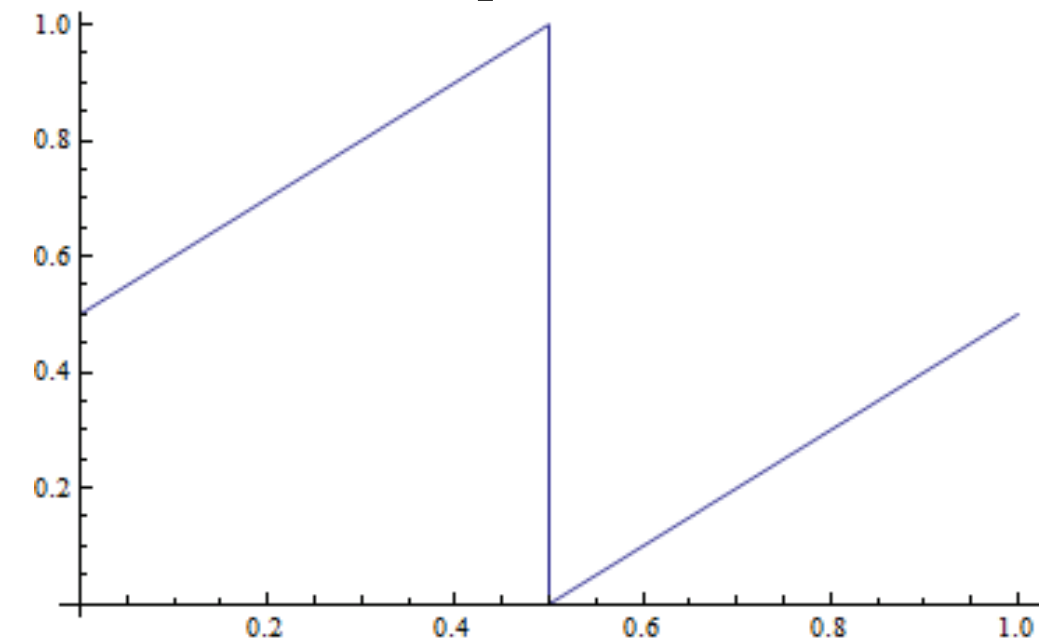

The discontinuity at x = 1/2 destroys the stable fixed point at {1/2,0} (except for integer values of K,-in which case it is a 'cat' map.) yielding a map which is everywhere hyperbolic. The divided phase space of the Standard Map is often not desirable in 'transport' problems which arise for example in fusion reactors. This is an example of an Anasov map. There is a more 'benign' version of this map which preserves the stable fixed point to retain a phase space which contains both regular and chaotic dynamics.This is what Peter Ashwin calls the Sawtooth Standard Map. The corresponding sawtooth map is **Saw[x_]:** $= \boldsymbol{x} - \textbf{Floor}[\boldsymbol{x}] - \mathbf{1/2}$

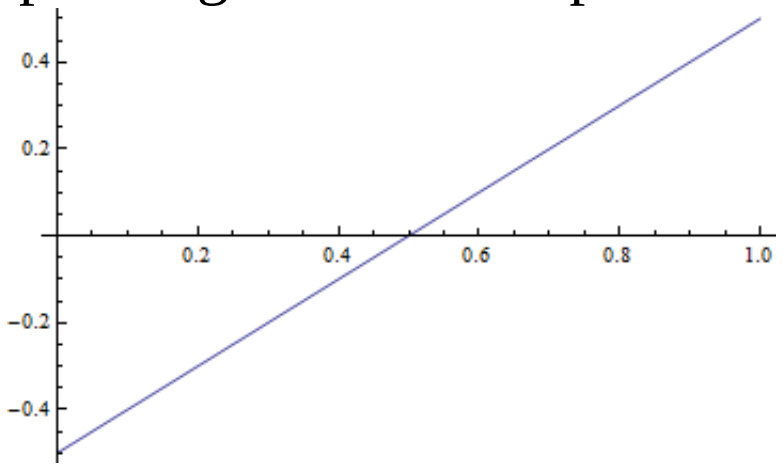

On the interval [0,1], Saw[Saw[x]] = x - 1/2. Below are the (inverted) Standard Map and the Sawtooth Standard Map (both are Mod 1)

$$\begin{bmatrix} x_{k+1} \\ y_{k+1} \end{bmatrix} = \begin{bmatrix} x_k + \dfrac{K}{2\pi}\sin 2\pi y_k \\ y_k + x_k + \dfrac{K}{2\pi}\sin 2\pi y_k \end{bmatrix} \qquad \begin{bmatrix} x_{k+1} \\ y_{k+1} \end{bmatrix} = \begin{bmatrix} x_k + kSaw(y_k) \\ y_k + x_k + kSaw(y_k) \end{bmatrix}$$

For the (inverted) Standard Map, the Jacobian is

$$J = \begin{pmatrix} \dfrac{\partial f_1}{dx} & \dfrac{\partial f_1}{dy} \\ \dfrac{\partial f_2}{dx} & \dfrac{\partial f_2}{dy} \end{pmatrix} = \begin{bmatrix} 1 & K\cos 2\pi y \\ 1 & 1 + K\cos 2\pi y \end{bmatrix} \text{ so at } \{0,1/2\},\ J = \begin{bmatrix} 1 & -K \\ 1 & 1-K \end{bmatrix}$$

The Sawtooth map also has a fixed point at $\{0,1/2\}$ and the Jacobian is $\begin{bmatrix} 1 & k \\ 1 & 1+k \end{bmatrix}$ so k = -K relative to the Standard Map. Note that (unlike the Standard Map) the Jacobian is constant. In both cases Det[J] = 1 so they are area preserving (symplectic). This implies that the eigenvalues come in pairs {$\lambda_1$, $1/\lambda_1$). Here $\lambda_1 = \frac{1}{2}(2 + k - \sqrt{k}\sqrt{4+k})$ . The eigenvalues are complex for $k \in [-4,0]$ . In this range the eigenvalues have the form $e^{\pm i\theta}$ where θ is given by $k = 2\cos\theta - 2$ , so we can use either k or θ as the parameter of the map. The range for θ is [0, π] which is perfect for regular polygons. Using the Saw function defined above, the Sawtooth Standard Map is **StSaw[{x_, y_}] := {Mod[x + k*Saw[y], 1], Mod[x + k*Saw[y] + y, 1]};**

Peter Ashwin [AP] showed that this map is equivalent to the Df map via a linear change of coordinates. We will give details below. As with Df, the web for a regular N-gon can be obtained using θ = 2*Pi/2N not 2*Pi/N.

Example: **k = N[2*Cos[2*Pi/14]-2]:** The following plot generates a web by scanning the x axis from 0 to 1 at depth 1000:

**StSawWeb = Flatten[Table[NestList[StSaw, {x, 0}, 1000], {x, 0, 1, .013}], 1];**

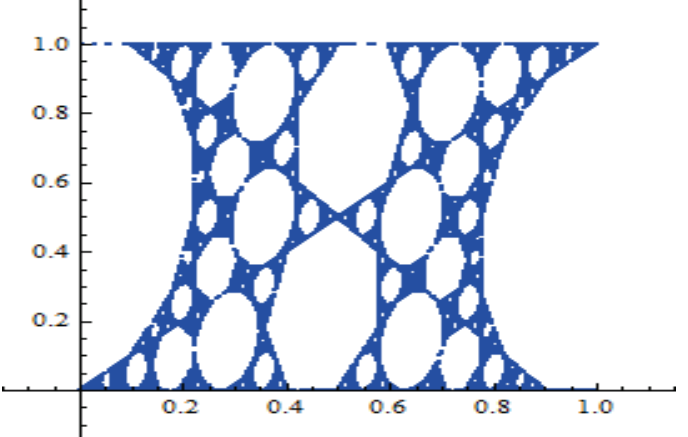

Below are the first 3 iterations of the web overlayed on the above plot.

**XAxis = Table[{x, 0}, {x, 0, 1, .001}]; T1 = XAxis; T2 = StSaw/@T1; T3 = StSaw/@T2;**

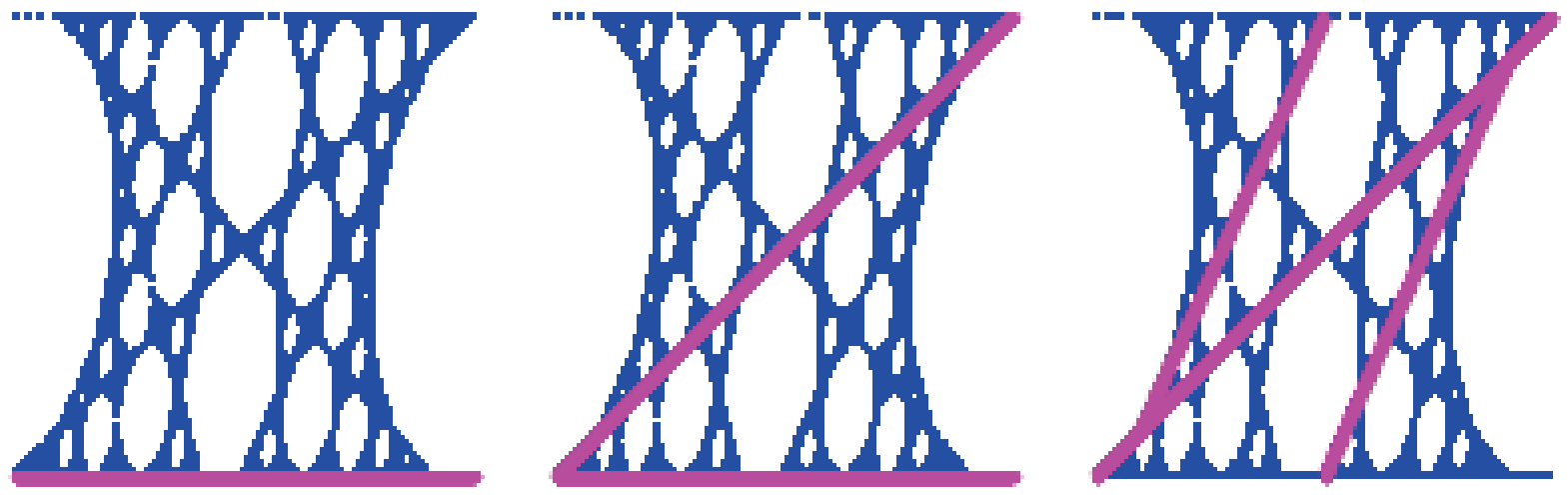

Just like Df, these webs are perfect matches for the Tangent Map webs, but the dynamics are different. To map StSaw space to Ts space for the Tangent Map, we will use Peter Ashwin's linear transformation from StSaw to Df-space and then use DfToTs. Here is an example:

The canonical 'remainder function' is **r[x_] := x - Floor[x]** (which is also SawtoothWave[x]]. Using this function, the change of coordinates from {x,y} in StSaw to {u,v} in Df is:

$$\begin{bmatrix} u \\ v \end{bmatrix} = \begin{bmatrix} 2r(x-y)-1 \\ 1-2y \end{bmatrix}$$

In Mathematica: **StSawToDf[{x_, y_}] := {2*r[x - y] - 1, 1 - 2*y};**

Example: θ= 2*Pi/14 so, **k = N[2*Cos[2*Pi/14] - 2]** as above.

**StSawWeb = Flatten[Table[NestList[StSaw, {x, 0}, 1000], {x, 0, 1, .01}], 1];**
**DfWeb =StSawToDf/@StSawWeb; TsWeb = DfToTs/@DfWeb;**

Below is StSawWeb in pink, DfWeb in blue and TrWeb in black. All these maps are invertible.

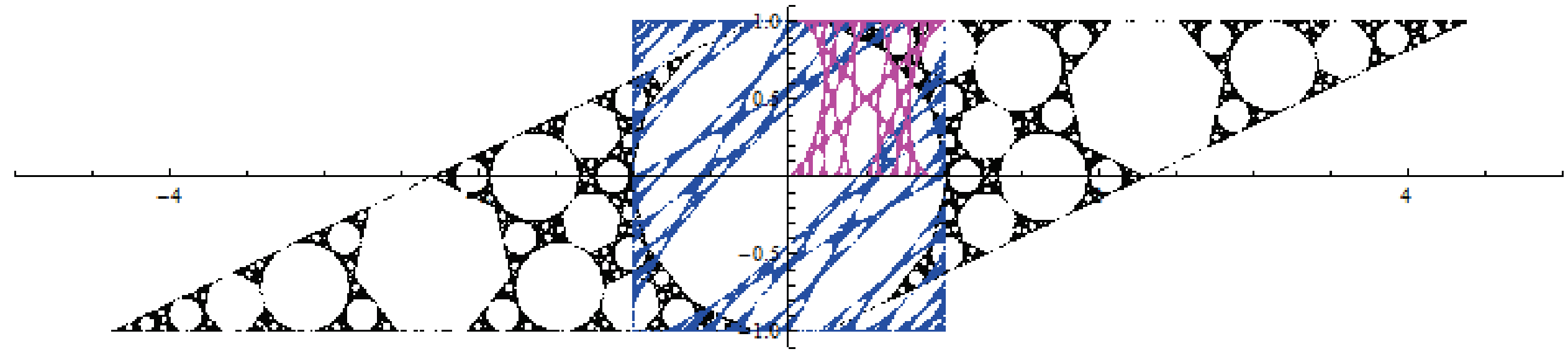

This shows a rather remarkable chain of conjugacies which has its roots in the Standard Map which is a special case of a general Kicked Map where the kicks are applied relative to the 'position' coordinate q, but they alter the 'momentum' coordinate p.

$$\begin{bmatrix} p_{k+1} \\ q_{k+1} \end{bmatrix} = \begin{bmatrix} p_k + Kf(q_k) \\ q_k + p_{k+1} \end{bmatrix}$$

The 'sawtooth' maps replace the Standard Map 'kick' function $f(q) = \sin 2\pi q$ with piecewise linear functions. Jurgen Moser proposed the Tangent Map for the same reason - to see what happens to the KAM Theorem when the 'smoothness' of the perturbation function is violated.

The story of why Moser chose to make 'smoothness' assumptions, rather than the 'analytic' assumptions used by the Russians is interesting - both from a historical perspective and a mathematical perspective. Because of the Iron Curtain, there was very limited exchange of information. In the end, the Russians though Moser's approach was superior (because real analytic functions are a proper subset of smooth functions), but Moser felt that Kolmogorov's insights were the key element in the KAM Theorem. In 1971, Moser [MJ2] wrote a textbook on the subject, explaining both approaches.

In his original version of the 'Twist Theorem' Moser needed to assume a function with 333 continuous derivatives ($C^{333}$) but this was later refined by Russman to a $C^5$ curve and Moser conjectured that it may even may hold for $C^2$. Counterexamples exist for the $C^1$ case (F.Takens - 1971). The Tangent Map and piecewise linear maps are neither analytic nor $C^1$, so they can be used as testing grounds to see what happens when the KAM Theorem breaks down.

**Section 6: The Kicked Harmonic Oscillator (KHO)**

The Standard Map is based on a kicked rotor in zero gravity, so there is no natural frequency of rotation for the unperturbed system - it is a free rotor whose frequency is determined by the initial conditions. The KHO is based on a harmonic oscillator, so it mimics the linear pendulum with gravity turned on. Now there is a natural frequency of oscillation ($\omega_0$) which depends on the length for the pendulum or the spring constant for the oscillator shown below:

x = 0

K.E. = $\frac{1}{2}mv^2 = \frac{(mv)^2}{2m} = \frac{p^2}{2m}$ where v is x'. The P.E depends on the 'spring constant' k. Hooks law says that (for small vibrations and a 'soft' spring), the restoring force on the spring is proportional to the displacement, so F = $k\omega_0 x$ for a constant k.

The Hamiltonian is $\frac{p^2}{2m} + k\omega_0 x^2$ and m can be scaled to get $\omega_0\left(\frac{p^2}{2} + \frac{x^2}{2}\right)$

The kicked system has the form: H = $\omega_0\left(\frac{p^2}{2}+\frac{x^2}{2}\right)+V(x)\sum_n \delta\left(\frac{t}{\tau}-n\right)$ where the kicking potental V(x) is usually of the form $\mu\cos kx$ or $\mu$sinkx ('even' or 'odd').
In quantum mechanics, KHO's have been used to model the behavior of charges in a uniform magnetic field subject to periodic 'kicks', electronic transport in semiconductor lattices, or cyclotron resonance where a resonant particle ( $\omega_0 = n\frac{2\pi}{\tau}$ ) will experience continuous acceleration. Typically the 'even' and 'odd' cases have very different dynamics. Below is an example of an even potential.

Using V(x) = $-\frac{\omega_0 K}{\tau}\cos x$ the Hamiltonian is: $\frac{1}{2}(p^2+\omega_0{}^2x^2)-\frac{w_0K}{\tau}\cos x\sum_n \delta\left(\frac{t}{\tau}-n\right)$
which reduces to the Standard map when $\omega_0\to 0$, $K\to\infty$ and $\omega_0 K$ = constant.

The discrete version is sometimes called the (stochastic) WebMap:

$$\begin{bmatrix} u_{k+1} \\ v_{k+1} \end{bmatrix} = \begin{bmatrix} \cos\alpha & \sin\alpha \\ -\sin\alpha & \cos\alpha \end{bmatrix}\begin{bmatrix} u_k + K\sin v_k \\ v_k \end{bmatrix}$$

where u = $x'/\omega_0$ , v = $-x$ and $\alpha = \tau\omega_0$. The linear form of this map is the unperturbed harmonic oscillator which yields simple rotation in phase space, so between kicks the particle rotates clockwise by α. For certain parameter values the Web Map models particle motion in a constant magnetic field.

In 'generic' {x,y} form, the equations are:
**KHO[{x_, y_}] := {(x + K*Sin[y])*Cos[w] + y*Sin[w], -(x + K*Sin[y])*Sin[w] + y*Cos[w]};**

Example**: K = 1.5; w = 2*Pi/4; H2 = Table[{x, x}, {x, -6, 6, .0374}];**
**T1 = Flatten[Table[NestList[KHO, H2[[k]], 1000], {k, 1, Length[H2]}], 1];**

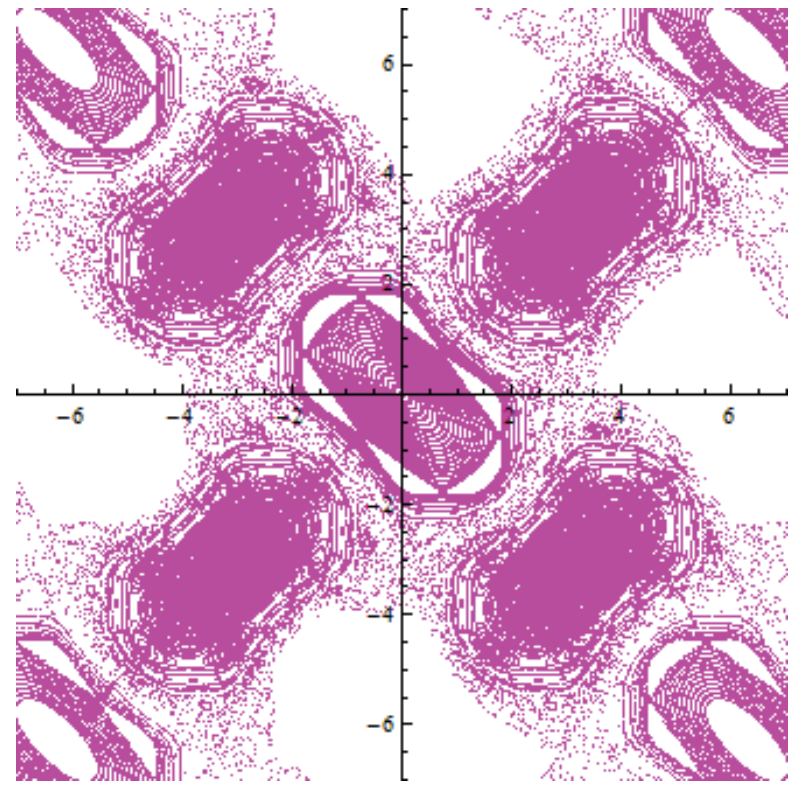

Since w = α = $\tau\omega_0 = 2\pi/4$, the plot above represents a 1:4 resonance between the kick frequency ($\omega_k = 2\pi/\tau$) and the resonant frequency $\omega_0$. It is the ratio of these frequencies which determines the global behavior of the system. For the Tangent Map, the case α = 2π/4 corresponds to the square N = 4. In this case the web tiles the plane with congruent squares (or rectangles) as we see above. This is called 4-fold symmetry. The cases for α = 2π/3 ( N= 3) and α = 2π/6 (N = 6) also tile the plane in a regular fashion. Actually these two tilings are congruent.

The KHO map for the regular pentagon, N = 5 ( α = 2π/5, K = 1.5) is shown below. This is a quasi-crystalline structure which shows less robust diffusion than N=3,4,6. The right side is the corresponding Tangent Map for N = 5 showing the rings of D's which confine the dynamics.

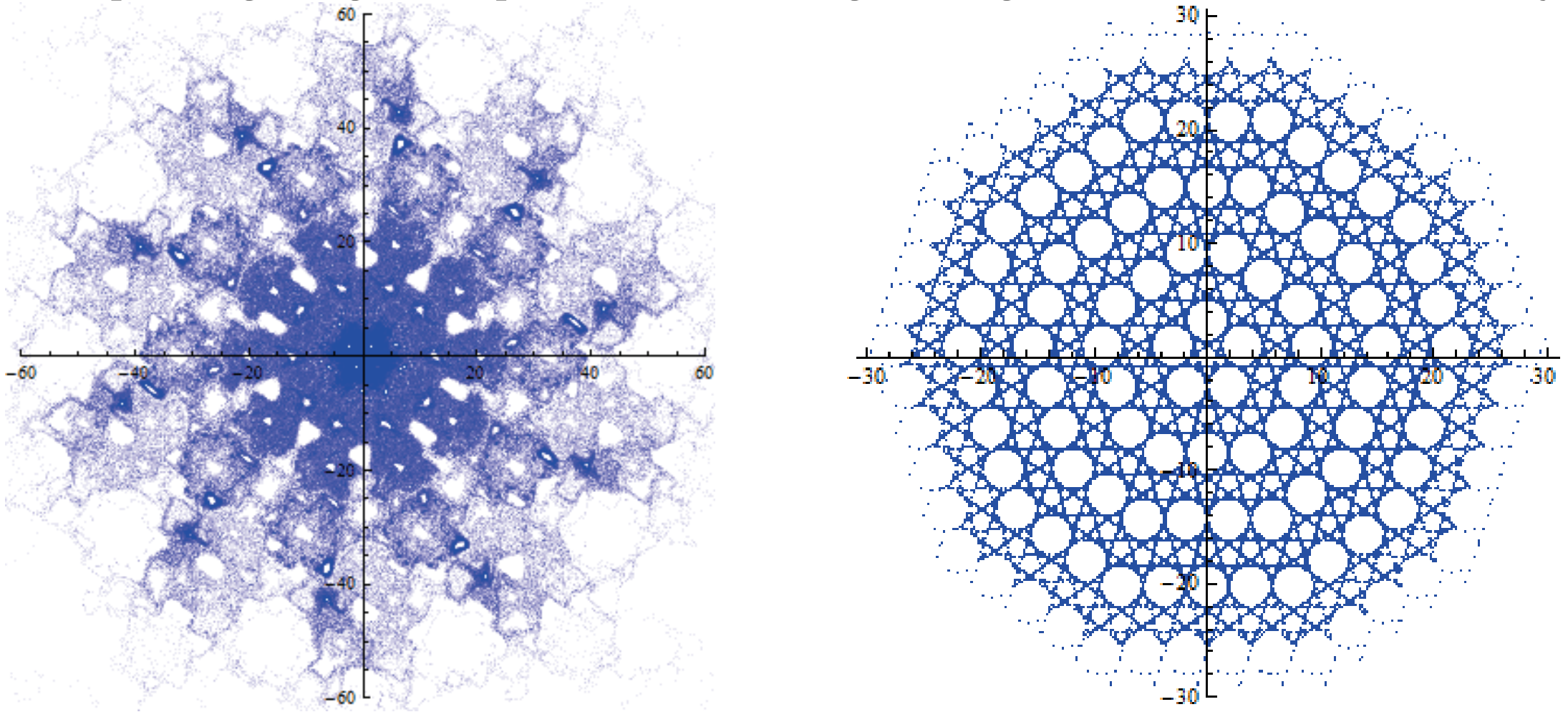

The N = 1 and N = 2 cases are degenerate lattices which correspond to operating regimes of the early cyclotrons. For other integer N values the structure is quasi crystalline and the corresponding quantum evolution show initial suppression of diffusion - until K is large enough to allow tunneling. For irrational ratios, all traces of crystalline structure disappear, but there is still diffusion at high energy levels. This case corresponds to a Tangent Map on a circle or ellipse and all orbits are also circles or ellipses with no intervening channels.

In quantum mechanics the most important cases are N = 3, 4 and 6, because the corresponding tiles represent extended quantum eigenstates and the channels allow diffusion over time, so they are testing grounds for quantum delocalization- of the type observed in the cold-atom experiments. The classical model of diffusion over time is nearly linear as shown below, while the quantum model shows sharp drop-off due to delocalization.

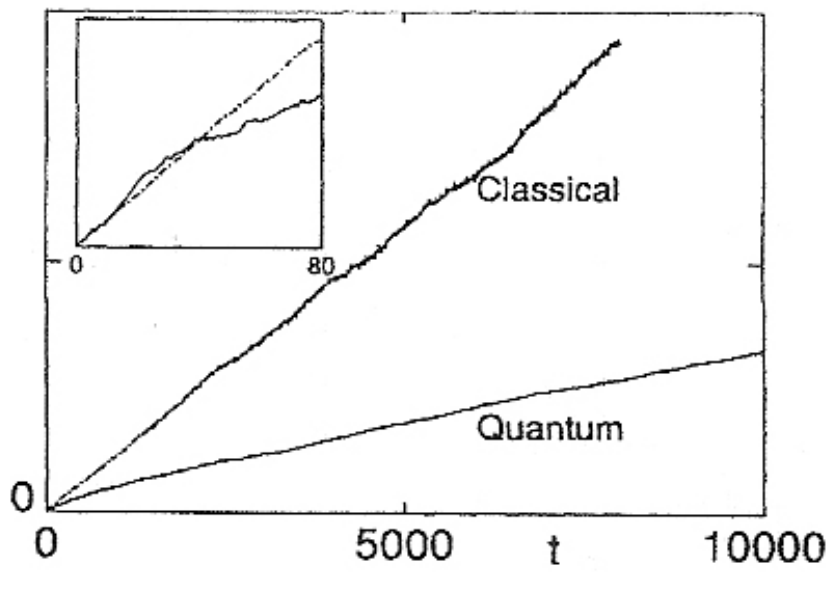

Once again we note that the maps shown here only simulate the classical predictions. To track the quantum evolution of the system it is necessary to generate time-dependent waveforms - usually with the aid of Fractional Fourier Transforms.

Over the past 15 years many laboratories world-wide have replicated the Texas University cold-atom experiments showing this drop-off at the classical-quantum transition. Experiments have advanced from cold-atoms to more robust Boise-Einstein condensates. But these experiments are based on the kicked-rotor model and the corresponding Standard Map. By contrast, the kicked harmonic oscillator models particle flow in lattices - and this much more difficult to observe. The most promising material for observing this diffusion is lateral surface superlattices constructed at nano scales as shown here.

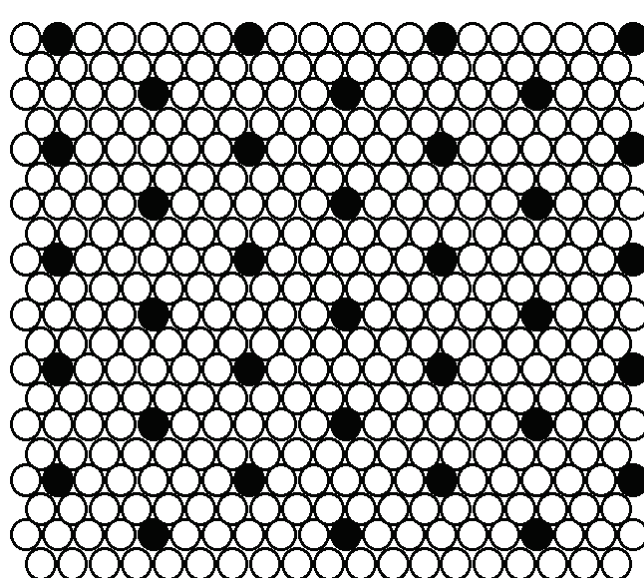

**The Harper Kicked Map**

The N = 4 case is known as the Harper Kicked Map (HKM). The kicked rotor and the Harper kicked map are the two most widely studied models in quantum mechanics. Lowenstein et al have studied this case in the context of 'sticky' orbits where iterates seem to be attracted to the neighborhood of stable 'islands' for long periods of time before being ejected.

HKM is based on Harper's 1955 equation - which is a stationary Schrodinger equation $\hat{H}\psi = E\psi$ .The Hamiltonian has the form: $\hat{H} = 2\cos\hat{p} + 2\cos\hat{x}$ where $\hat{p} = -i\hbar\partial/\partial x$ is the effective momentum and $\hbar = 2\pi\sigma$ where σ is proportional to the magnetic field and has the units of flux quanta per unit cell area. The kicked Harper Hamiltonian has the form:

$$\hat{H} = L\cos\hat{p} + K\cos\hat{x}\ \delta_1(t)$$

where $\delta_1(t)$ is a period 1 delta function. The discrete form is the Harper kicked Map:

$$\begin{bmatrix} x_{k+1} \\ p_{k+1} \end{bmatrix} = \begin{bmatrix} x_k \\ p_k \end{bmatrix} + \begin{bmatrix} -K1\sin p_{k+1} \\ K2\sin x_k \end{bmatrix}$$

In Mathematica: **HKM[{x_, p_}]:= {x - K1*Sin[p +K2*Sin[x]], p+ K2*Sin[x]}**

**Example**: Below is an unaltered overlay of HKO in blue (K1 =K2 = 1.5) and KHO from above in red ( K = 1.5 and w = 2*Pi/4). Using the same diagonal scan parameters as above. They are a perfect match but KHO has radial symmetry - which we can match with a dual scan of HKO.

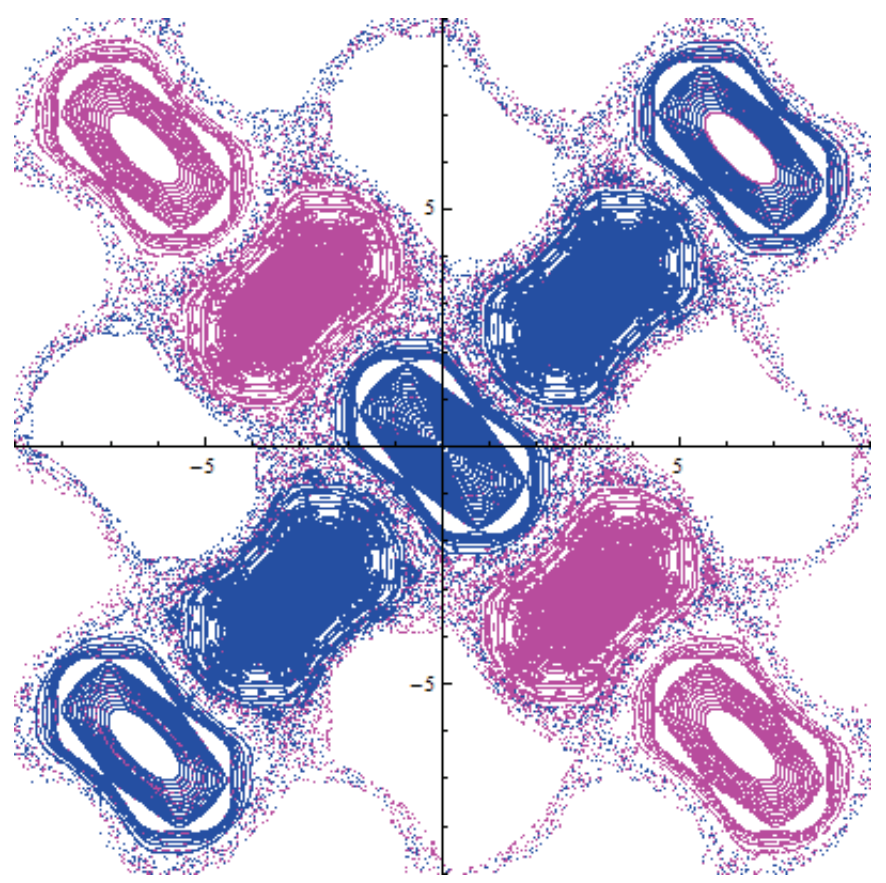

When K1 = K2 = k →0, this becomes the unkicked Harper map, but as soon as k > 0 there are stochastic layers of the type we see in the Standard Map. As k increases these layers widen and the resonant islands break down allowing diffusive motion in x and p. By k = 5 there are no visible islands and the map may be ergotic.

The quantum form of this motion is obtained by solving the time-dependent Schrodinger equation. Since the kicks are periodic, the Floquet Theorem says that there are well-defined quasienergy eigenstates . As the flux, σ, is varied these energy spectrum yield beautiful Hofstadter butterflies - one of which is shown below. The quasienergies can be obtained by taking the Fourier Transform of the time evolution of the wave packet - starting with a random initial packet. For example when $k/\hbar = 1$, σ = 2πk and the plot below shows the quasienergy levels in the range of -2 to 2 on the y axis as σ is varied from 0 to 1.

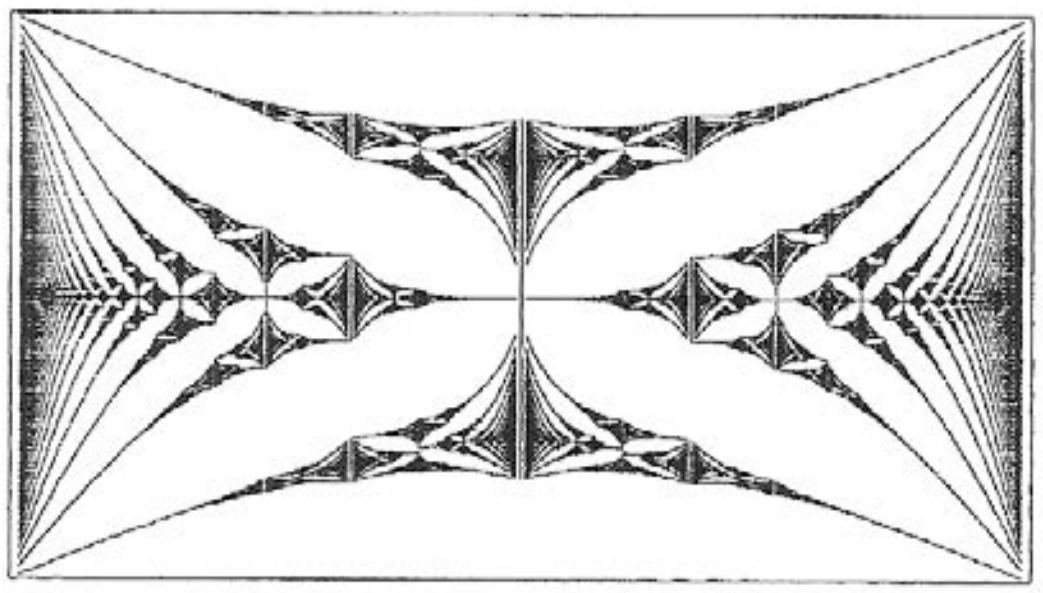

To see the divergence in the HKO stochastic layer we will map the small red grid centered at {-3.125,0}: **K1 =1.5; K2= 1.5**

**M[n_]:= Flatten[Table [NestList [HKM, grid [[k]],n],{k, 1, Length [ grid]}] , 1];**
Shown below are M[0], M[5] and M[10]

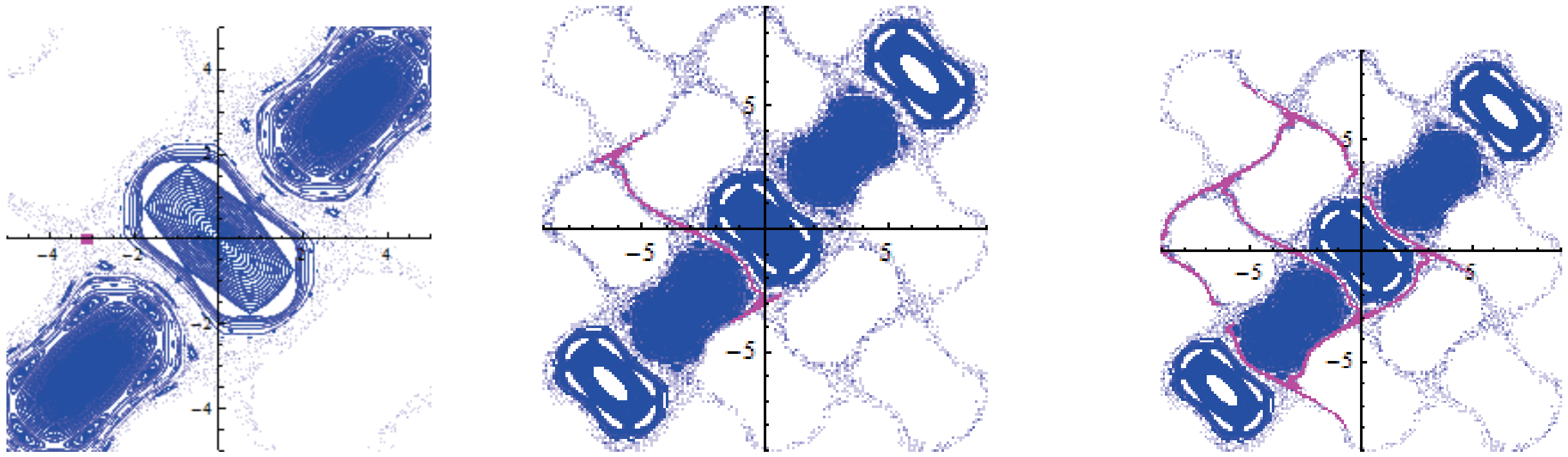

The Standard Map also has chaotic transport. It is easier to see if we consider it as a map on a cylinder and drop the Mod 1.

```
Std[{x_, y_}] :={y + (K/(2*Pi))*Sin[2*Pi*x] + x, y + K/((2*Pi))* Sin[2*Pi*x]} ; K= .97165
```

The initial grid is centered at {.5,.5}.The first 12 iterations are shown below. It appears that these iterates will eventually be 'dense' in the stochastic layer but we know that there will be embedded islands of stability at all scales - so there will be smaller stochastic layers like this at all scales.

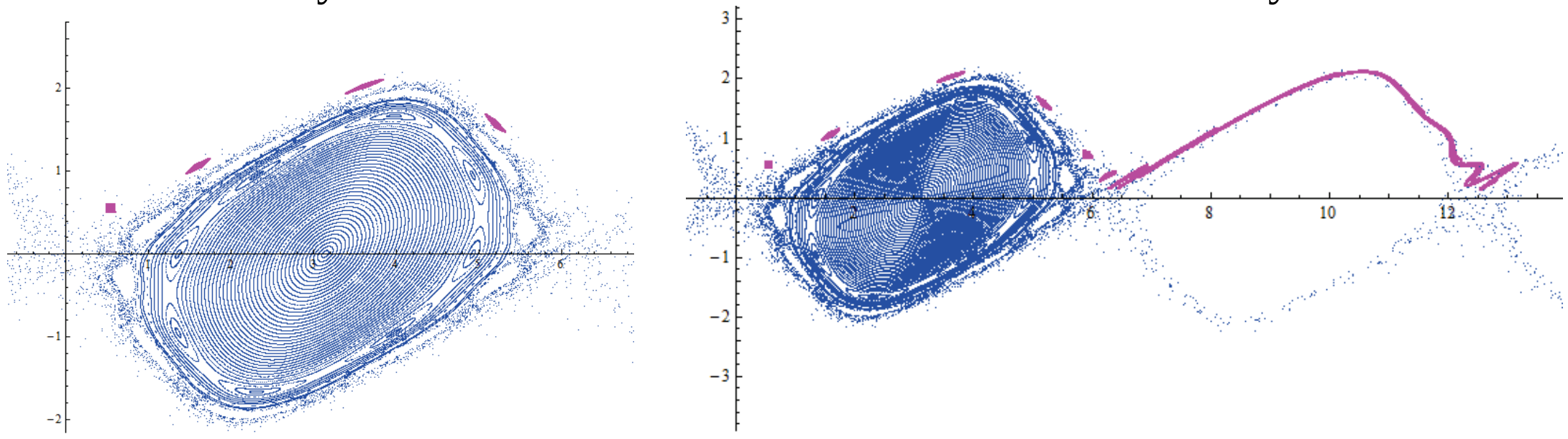

Since these maps are both area-preserving, the apparent increase in red area is due to the fact that the rectangle is stretched horizontally and contracts vertically. The rate of expansion and contraction are determined by the eigenvalues and for an area preserving map they are related by $\lambda_1 = 1/\lambda_2$ In the Harper Map and the Standard Map this expansion/contraction is accelerated when the iterates pass close to an unstable saddle point. This is known as the 'ballistic' effect.

The Poincare Recurrence theorem says that for an area preserving map T and any set U of positive measure such as the rectangle above, then the points that do not return to U form a set of measure 0. So almost all points in U return to U. This can be seen with the Harper map and Standard Map- if we are patient enough. Below is iteration 70 of the Standard Map it shows some points returning close to the initial small rectangle.

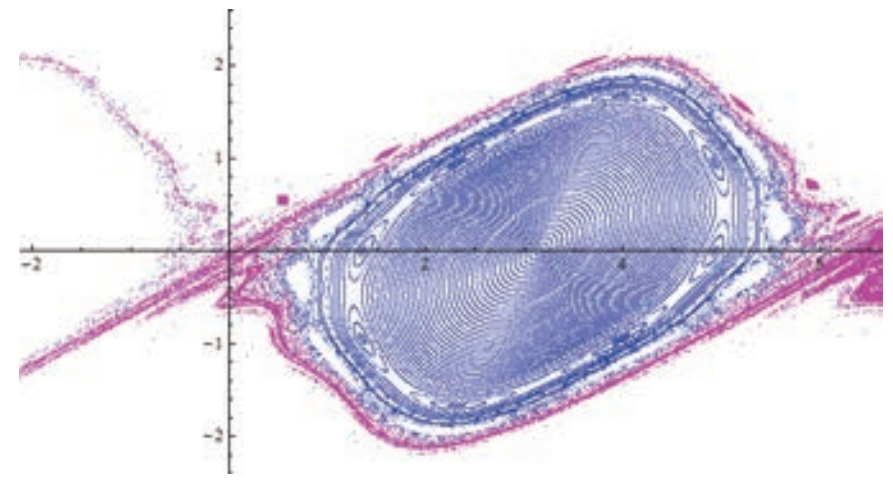

Because these points are confined to the outer regions, they will not visit the whole space, so for these K values, the Harper map and the Standard Map are not ergotic. A map T on a set S is ergotic if for any set U of positive measure, $T^k(U)$ will eventually cover S, except possibly for a set of measure 0. For example a circle map such as $f(\theta) = 2\theta$ is ergotic even though there are an infinite number of periodic orbits. These are exceptional orbits with measure 0 and 'almost all' orbits will be dense in the circle.

The Standard Map and the Harper map are not ergotic for small K values because they are based on integrable Hamiltonians and the KAM theorem says for small K values there will be invariant 'islands' which survive. For larger K values these islands appear to break down and there is evidence that both of these may be ergotic for large enough kick strength K. In these maps it appears that no islands of stability remain for large K values, but this is very difficult to prove.

What makes the kicked Haper Map interesting is that it exhibits diffusion in both the classical and quantum versions - although the connection between the two is poorly understood. In earlier 1-dimensional models the quantum diffusion was quenched by localization.One of the major differences between the classical model as shown here and the actual quantum behavior is that in the quantum model there is diffusion for all initial conditions.This is certainly not the case with the Standard Map or the Harper kicked map. For example, the origin {0,0} is a fixed point of the HKM, but the energy of this quantum state grows with time due to tunneling.

For the Tangen map with a regular polygon, the web structure consists of concentric rings of 'Dads' which serve as bounds for the dynamics. In the non-regular case, these rings may break down and allow points to diverge. The only proven unbounded case is the Penrose Kite but numerical simulations indicate that such orbits are pervasive for large classes of polygons. This means that many polygons support a form of ‘unbounded transport'. There is a class of polygons formed from nested regular polygons. Some of these woven polygons show a form of transport which is similar to the gravitational webs which pervade the solar system and are sometimes used as gravitational assists for space probes.

Note:We saw above that the Harper Map is a special case of the Web Map (KHO) so it has 'normal' form $\begin{bmatrix} y_{k+1} \\ x_{k+1} \end{bmatrix} = \begin{bmatrix} \cos\alpha & \sin\alpha \\ -\sin\alpha & \cos\alpha \end{bmatrix} \begin{bmatrix} y_k + f(x_k) \\ x_k \end{bmatrix}$ where ω= 2π/4 and *f* is the kick function

so $\begin{bmatrix} y_{k+1} \\ x_{k+1} \end{bmatrix} = \begin{bmatrix} 0 & 1 \\ -1 & 0 \end{bmatrix} \begin{bmatrix} y_k + f(x_k) \\ x_k \end{bmatrix}$ This means that the map has the form F(y,x) = {x,- y- *f*(x)} and the kick function $f$ can be taken to be the sawtooth map *f*(x) = *a*x mod1. In this way the global dynamics of the Harper Map can be reduced to a Df type map on a torus:[0,1) → [0,1). This is possible because of the four-fold symmetry of the Harper Map. In the next section we will look at another sawtooth version of the KHO.

**Section 7: A Dissipative Kicked Harmonic Oscillator**

A dissipative form of the KHO can be written as: $\omega_0\left(\frac{p^2}{2}+\frac{x^2}{2}\right)+\mu\,|\,x\,|\sum_n \delta\left(\frac{t}{\tau}-n\right)$

A.J.Scott,*et al.* **[SHM]** used this Hamiltonian to study circle packing. Setting μ = 1, their discrete mapping is:

$$\begin{bmatrix} x_{k+1} \\ p_{k+1} \end{bmatrix} = F\begin{bmatrix} x_k \\ p_k \end{bmatrix} = \begin{bmatrix} \cos\omega & \sin\omega \\ -\sin\omega & \cos\omega \end{bmatrix}\begin{bmatrix} x_k \\ p_k - \operatorname{sgn} x_k \end{bmatrix}$$

The absolute value function models the bi-directional form of the dissipation. As usual, the 'kicks' are position dependent, but they are applied to the momentum. They create an s-sequence on the symbols {1,0,-1} where $s_k = \operatorname{sgn}(x_k)$.

$$\begin{bmatrix} x_{k+1} \\ p_{k+1} \end{bmatrix} = \begin{bmatrix} \cos\omega & \sin\omega \\ -\sin\omega & \cos\omega \end{bmatrix}\left(\begin{bmatrix} x_k \\ p_k \end{bmatrix} + \begin{bmatrix} 0 \\ -\operatorname{sgn} x_k \end{bmatrix}\right)$$ Setting $z_k = x_k + ip_k$ $\mathrm{F}(z_k) = e^{-i\omega}(z_k - is_k)$

Comparing this to the WebMap above for the Kicked Harmonic Oscillator, the only difference is that the Sine function is now a Sign function - so DKHO is a sawtooth version of KHO. Like the Sawtooth Standard Map, this is a fundamental change because it allows the model to be 'tuned' to the frequency ω without the 'noise' of the Sine function. The 1-bit DAC in the Adc map was also modeled with a Sign function, but by necessity it was applied symmetrically to x and y.

In the Mathematica code below, we swap the order of the variables, so that the line of discontinuity is the x-axis, rather than the y-axis. This creates a scenario similar to the Tangent Map, where we scan the x-axis to construct webs.

```
DKHO[{x_, y_}] := {-y*Sin[w] + (x - Sign[y])*Cos[w], y*Cos[w] + (x - Sign[y])*Sin[w]};
```

Like KHO, DKHO is defined globally, but symmetry is broken and the effects are most pronounced in the vicinity of the origin. The 'kicks' have diminishing effect as we move outwards but the small scale dynamics local to the origin are preserved at all distances. This is similar to the Tangent Map for regular polygons where the (invariant) dynamics inside the first ring of Dads is a template for the global dynamics.

For regular polygons under the Tangent map, the key issue is to understand how the web evolves on each iteration. This 'cobwebbing' mimics the methods used to track the dynamics of functions of one variable. For the regular pentagon the cobwebbing is self-similar and hence can be described analytically. For most polygons there is no simple analytic description , but it is still possible to use the geometry of the web to prove basic facts - such as invariance of the inner star region. We will briefly show how the N = 5 web evolves under the DKHO Map and then compare this with N = 7 web.

Example: **w= 2*Pi/5; Xaxis=Table[{x,0},{x,-10, 10,.01}];**
**Web=Flatten[Table[NestList[DKHO,Xaxis[[k]], 500],{k,1,Length[Xaxis]}],1];**

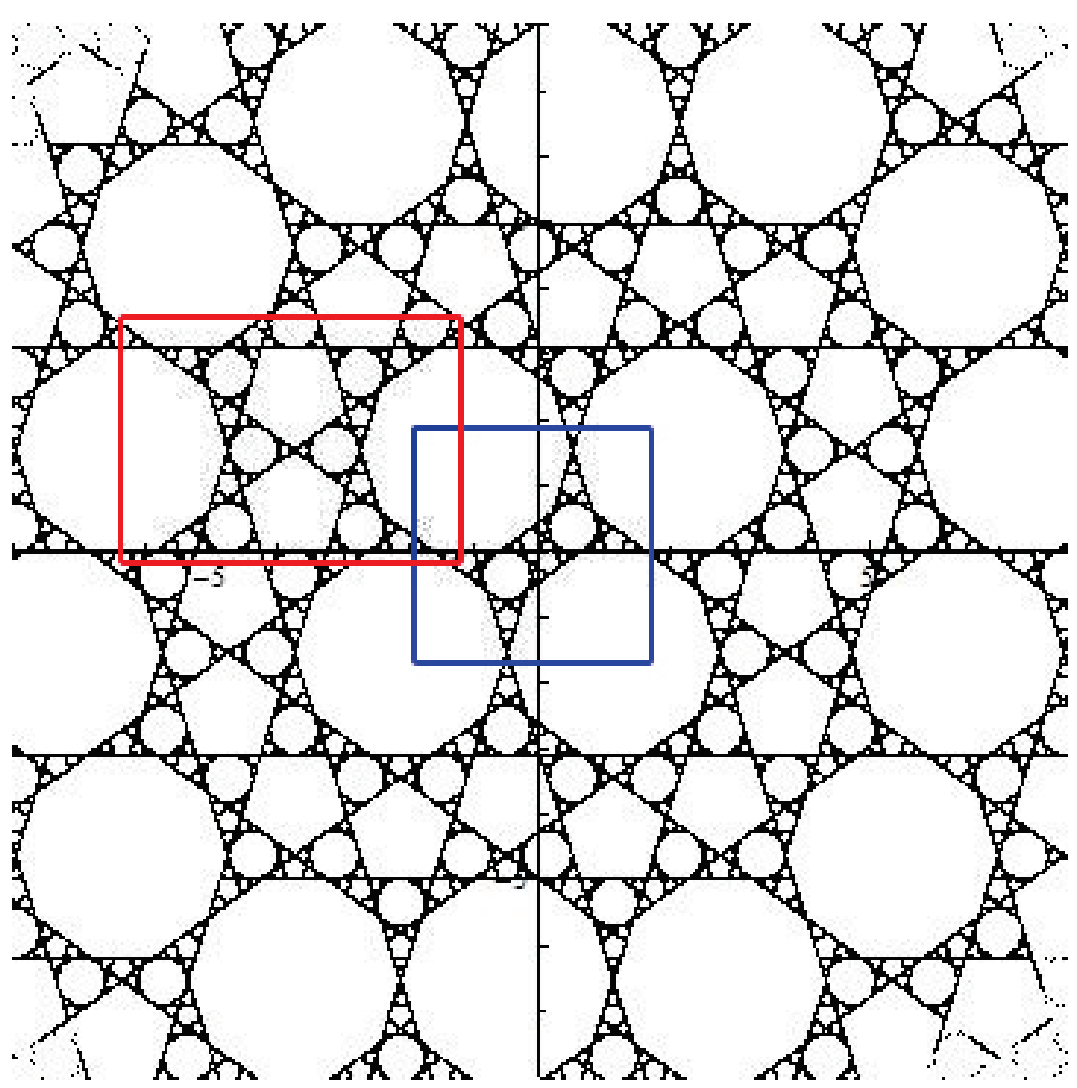

This is the level-500 web.The web outside the first ring of 4 D’s is locally a perfect reproduction of the N = 5 Tangent Map web. The red rectangle outlines a typical family with M, S1 (alias D[1]), D and a small M[1] on the edge of D. The blue rectangle is the invariant 'inner star' region which is bounded by the 4 D’s. If this was a canonical N = 5 inner star region, it would be bounded by a ring of 5 D’s and contain a full generation. Here M is missing, but the 4 D’s are fostering normal families. The two S1 progeny are shared by two D’s and four of the M[1]’s are also shared

Below are the level 1, level 2 and level 3 webs for this inner star region. It is easy to find the exact parameters of the D’s and we have drawn one for reference.

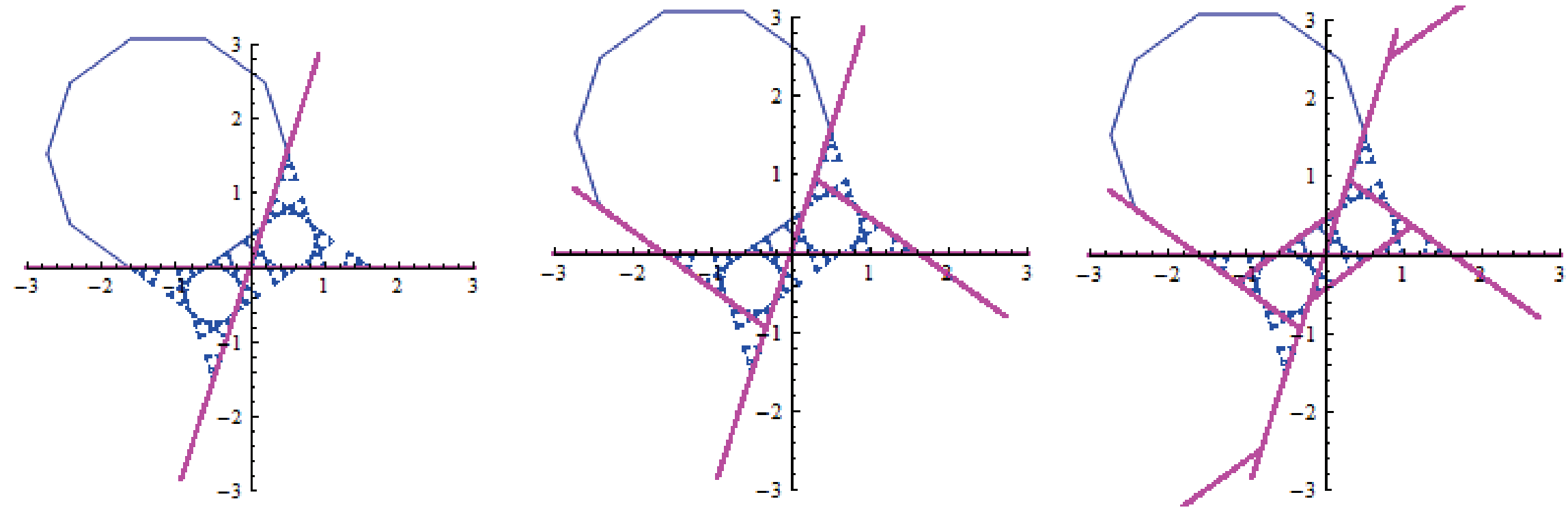

The level 1 web on the  left above is a rotation of the x-axis about the origin. The angle of rotation is 2*Pi/5. There are no 'kicks', because y = 0  On the next iteration, these points have non-zero y value (except the origin) so they will be displaced by {x,y} → {x-Sign[y],y}and then rotated by w. Therefore in any neighborhood of the origin, the points from the first generation will be displaced backwards (or forwards) by magnitude 1 along the direction of the first iteration (and then rotated by w). This forms the level 2 web shown in the middle above. It also determines the ultimate length of a side of D which will be 1.( Note that the second edge of D is a step 2 from first edge. This same relationship will exist between the level 3 and level 4 edges).On the level 3 web, the two displaced 'origins' are now kicked forward or backwards by magnitude 1 and rotated by w to form two parallel edges of the opposing D's.

Because the limiting geometry of the N= 5 web is self-similar, it is easy to find points with non-periodic orbits. In the enlargement below M[1] is drawn in red. It has a vertex at the origin.

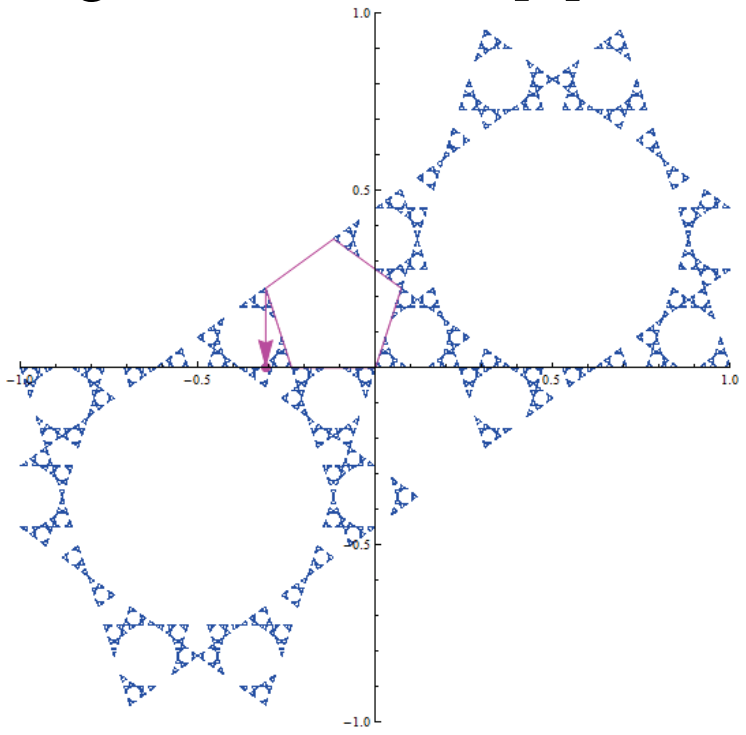

This is a second generation M, so the radius is rM*GenerationScale[5] and both of these are known. M has the same side as D, so rM= RadiusFromSide[1,5] ≈ 0.85065080 and GenerationScale[5]= (1-Cos[Pi/5])/Cos[Pi/5]. This allows us to generate M[1] (and D) to any desired accuracy. A corresponding non-periodic point is at the tip of the arrow above: s1 = {M1[[3]][[2]], 0}≈{-0.30901708364323477, 0} (or the symmetric points). The blue plot above was generated from 50,000 points in the orbit of s1. **Orbit = NestList[DKHO, s1, 50000];**
Note that it lacks the tell-tale edges which are artifacts of web scans.

The plot below is an overlay of the Adc plot from Section 3 in pink on top of the map above: **AdcOrbit = NestList[Adc, p1, 30000]; TrOrbit = DfToTs /@ AdcOrbit.** The generating point p1 is TsToDf[s2] where s2 is the symmetric limit point on the right side of M[1]. Since it is on the y-axis in Adc space, it would be fixed unless it is displaced by an arbitrarily small amount. We used a displacement of $\{\varepsilon,0\}$ where $\varepsilon < 10^{-12}$. The Adc plot is scaled by Cos[3w/4] to obtain a perfect overlay.

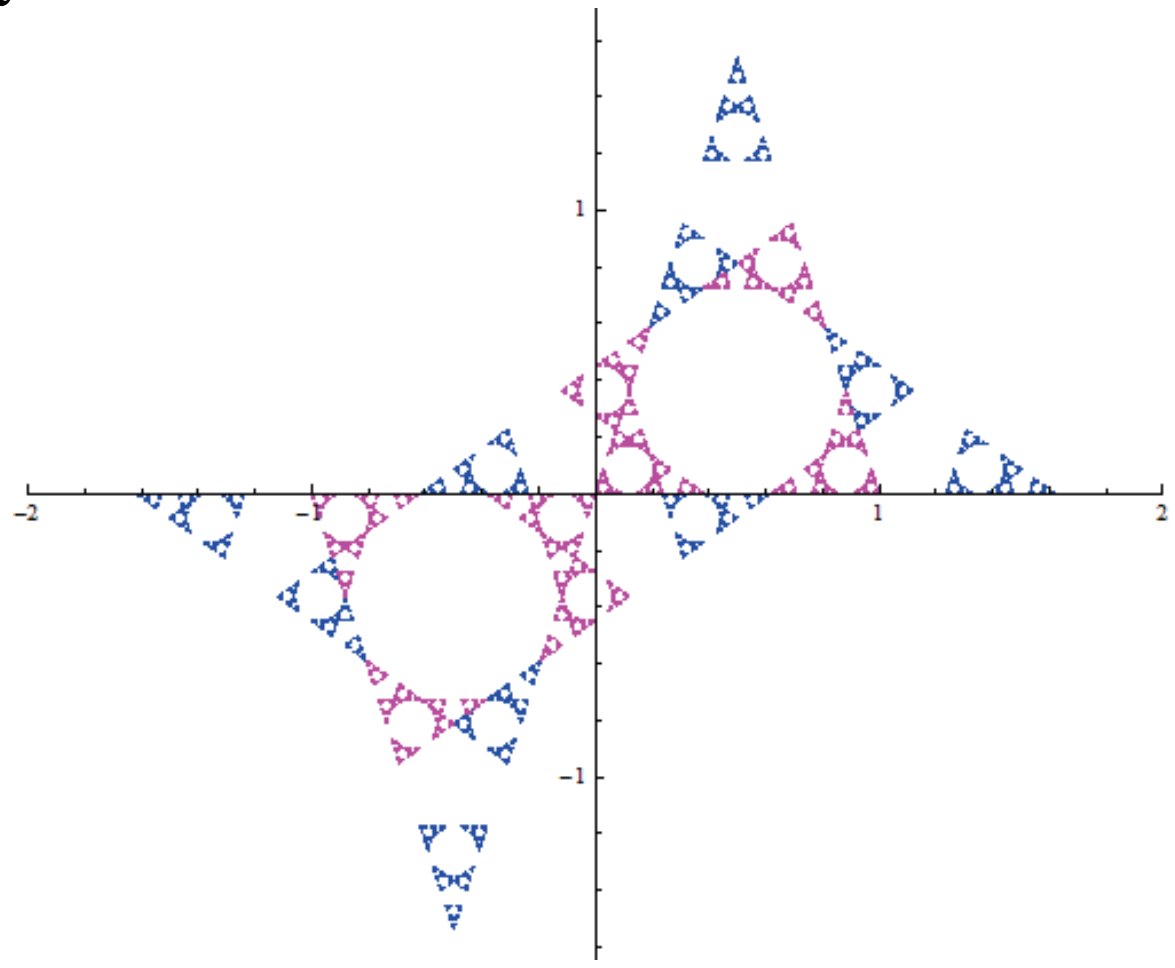

N = 5 is an exceptional case and N = 7 is more typical of what to expect for lager N values - namely very complex webs on the small scale and very predictable webs on a large scale. For all regular polygons with the DKHO map, it appears safe to use the regions outside the first ring of D's. Below is such a region for N = 7 using w = 2Pi/7. The central region is clearly devoid of M's but the shaded region is a perfect reproduction of the N = 7 web.

**w = N[2*Pi/7]; Xaxis = Table[{x, 0}, {x, 0, 4, .0174}];**
**T1 = Flatten[Table[NestList[DKHO, Xaxis[[k]], 3000], {k, 1, Length[Xaxis]}], 1];**

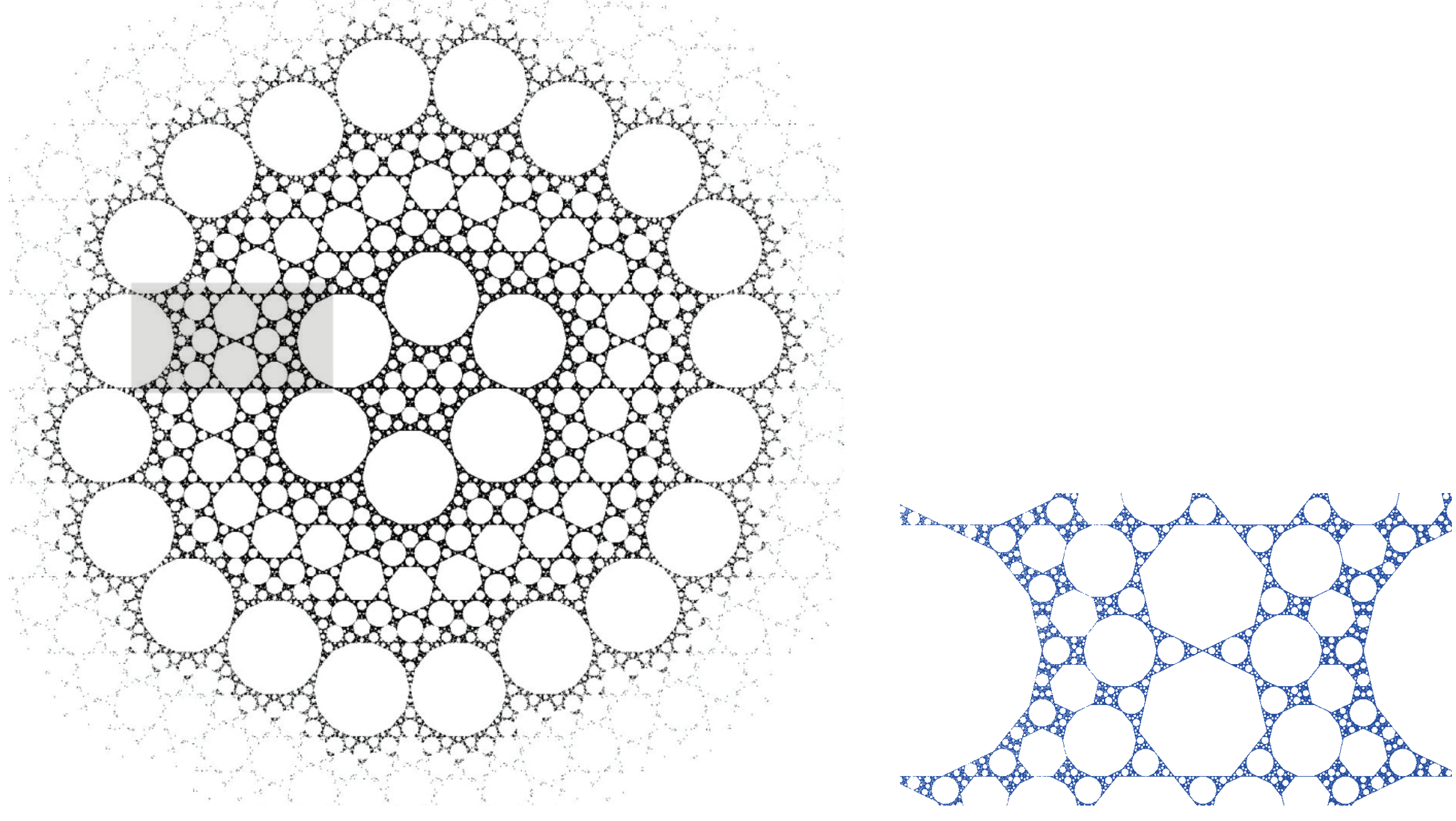

In Section 8 we will return to this map and show how to reclaim the origin by aligning the kicks with the edges of the N-gon. This will place one of the vertices of M at the origin and generate a locally perfect web.

**Section 8: Complex Valued Mappings**

In the early part of the 20th century, Henri Poincare and George Birkoff pioneered the qualitative study of planar maps because of their applications to physical process such as celestial mechanics and quantum mechanics. Their work laid the foundation for topology and the study of manifolds.

Because of this framework, the KAM authors could prove results about Hamiltonian systems in the context of symplectic mappings of the plane. The key issue for the KAM Theorem was what Moser called the 'function theoretic center problem' - the stability of neutral fixed points (the center) when an integrable Hamiltonian system was perturbed periodically. Poincare knew that in these cases, the linearization process might fail because of small divisors involving powers of the eigenvalues. In his famous 1885 paper on the n-body problem for the King Oscar II prize, he predicted that these small divisors could lead to chaotic behavior as shown below for the 3-body problem.

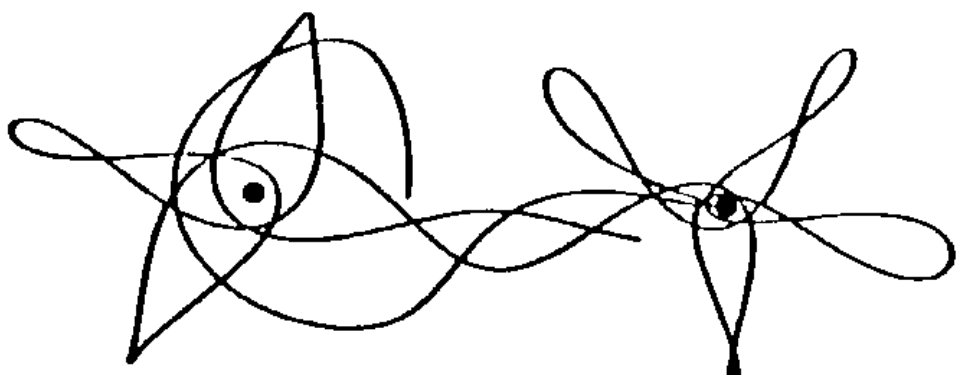

The KAM authors found clever convergence schemes which could be used to bypass the small divisor problem.We remarked earlier that the Russians assumed a mapping which was (real) analytic and Moser assumed that it was smooth.

Poincare was equally comfortable with real or complex analysis and he made fundamental contributions in both areas. The new science of topology provided a common framework. He formulated an 'idealized' center problem in complex analysis as follow:

*Given a holomorphic (complex analytic) function such as $f(z) = \rho z + z^2$, with $\rho = e^{i2\pi\theta}$, for what values of $\theta$ will $f$ be linearizible near the fixed point $z = 0$ ?*

A function such as *f* is 'linearizible' near a fixed point if it is conjugate to a rotation via a holomorphic change of coordinates. The corresponding linear map will have the form $L(z) = ze^{i2\pi\theta}$ so it will be a rotation by θ.

Note that in the example above, $f'(0) = \rho$, and $|\rho| = 1$, so $z = 0$ is a 'center' or indifferent fixed point. Poincare knew that when θ was rational, $|\rho - \rho^k|$ will vanish for some k and there was little hope of convergence, so the only issue is the case when θ is irrational. In this case $z = 0$ is called an irrationally indifferent fixed point (or holomorphic germ)

This problem has a long history and it was only recently been settled by A. Brjuno and J. Yoccoz using a class of irrationals of Brjuno type. In 1917 G. Pfeifer [P] gave a counterexample which showed that some irrational values of θ are not linearizible. Gaston Julia [J] (1893-1978) and Pierre Fatou (1878-1929) were working on this problem from a different perspective. For a given rational complex valued function *f*, they wanted to know which points are 'tame' under iteration and which points are 'exceptional' (chaotic). The latter points make up what we now call the Julia set of *f* ,and the former are the Fatou set. Neither Fatou or Julia knew how to classify the irrationally indifferent points and in 1919 Julia gave an incorrect proof that these points were never linearizable. In 1927 H. Cremer settle the case of rational θ by showing that they are not linearizable. We will see an example of this below.

The big breakthrough came in 1942 when Carl Siegel [SC] showed that if θ satisfied a certain Diophantine condition then it was linearizable. The condition guaranteed that θ could not be approximated closely by rationals.

Jurgen Moser was a friend and colleague of Siegel and in the early 1950's they both followed the developments in Russia as Kolmogorov claimed to have a proof of the stability problem. When they realized that Kolmogorov only had an outline of a proof, Moser worked on his own proof. He did not follow the path of Siegel's 1942 proof, because he knew that complex analytic functions were too 'well-behaved' to be applicable to Hamiltonian problems. His Twist Theorem assumed a real-valued 'smooth' function, while V.I. Arnold [A] (1937-2010) assumed a real analytic function in his proof of Kolmogorov's 'conjecture'. (In the complex plane the distinction between analytic and smooth disappears.)

Any complex polynomial can be used to illustrate Siegel's result. Suppose $f(z) = z^2 + c$. Then the fixed point is no longer $z = 0$ , but we can find it easily by solving $f(z) = z$. The solution depends on c. Suppose c is chosen to yield an indifferent fixed point $z_0$ with $|\rho| = |f'(z_0)| = 1$. We can assume that ρ has the form $e^{i2\pi\alpha}$ and normally this would imply that the motion is locally conjugate to a rotation by α, but when α is close to being rational this conjugacy might fail in any neighborhood of $z_0$ because the power series expansion has terms in the denominator of the form $|\rho^k - \rho|$. If $\alpha \approx m/n$ , this term could get arbitrarily small and this could happen even if α was irrational.

Siegel developed a Diophantine condition on α that would guarantee that the power series converged: There exists $\varepsilon > 0$ and $\mu > 0$ such that $\left|\alpha - \frac{m}{n}\right| > \frac{\varepsilon}{n^\mu}$ for all integers m,n with n positive. By this measure of irrationality, the safest possible number is one with a continued fraction expansion of 0,1,1,1,.. and this number is the Golden Mean $\alpha = \gamma = \left(\sqrt{5} - 1\right)/2$

It is no coincidence that the last surviving curve of the Standard Map has winding number γ. It appears that any irrational with a continued fraction expansion which ends in 1's ( a 'noble' prime) is a local 'survivor' while its neighbors are destroyed. This makes the Golden Mean the ultimate survivor .

When $f(z) = z^2 + c$, $f'(z) = 2z$ so we can find a Golden Mean center $z_0$ by setting $2z_0 = \rho = e^{i2\pi\alpha}$ with $\alpha = \gamma$ . In Mathematica: **gamma = (Sqrt[5]-1)/2;**
**z0 = (1/2)*Exp[2*Pi*I*gamma]** ≈ -0.36868443903915993 - 0.3377451471307618*I
Since this is a fixed point, $z_0^2 + c0 = z_0$ so c0 ≈ -0.3905408702184 - 0.5867879073469687*I

Mathematica does not care whether a function is real or complex, so define: **f[z_]:= z^2+ c0** and iterate *f* in the neighborhood of $z_0$:

**ObitScan =Flatten[Table[NestList[f, z0+ x, 1000],{x, 0, .55, .1}] ,1];**
Two special orbits: **Orbit1 =NestList[f, 0, 3000]; Orbit2 =NestList[f, -z0, 1];**
To plot a complex valued point: **CPoint[z_]:=Point[{Re[z],Im[z]}];**
**Graphics[{ CPoint/@OrbitScan, Blue, CPoint/@Orbit1, AbsolutePointSize[5.0], Magenta, CPoint/@Orbit2}, Axes->True]** (*on the left below*)

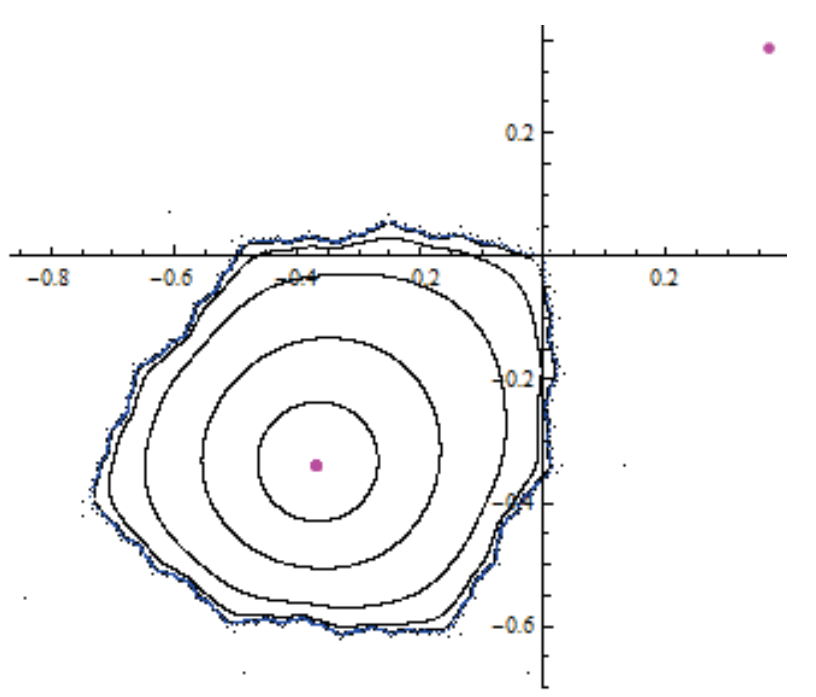

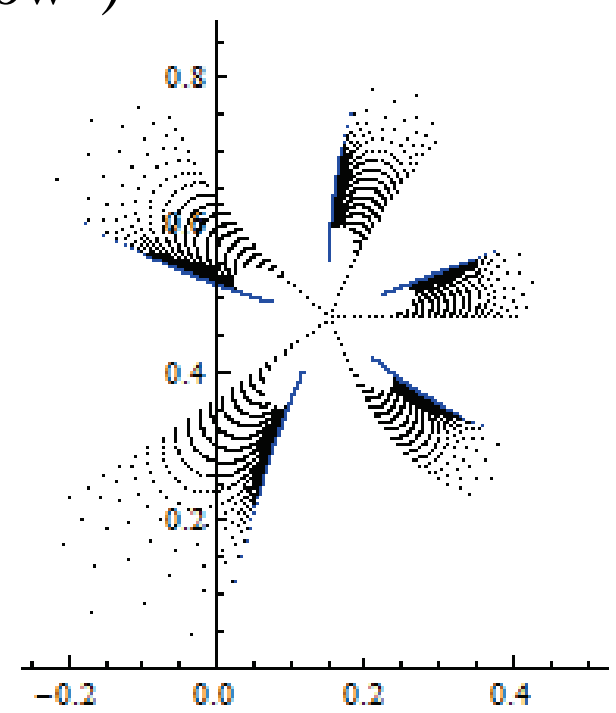

The invariant curves shown on the left are the complex analog of invariant KAM curves. (There are no invariant curves on the right hand plot because it corresponds to gamma = 1/5. This is a rationally indifferent fixed point and there are no stable regions of any size so linearization is impossible. The divergent channels seen here extend all the way to z0. Julia knew this and it is one reason why he was so pessimistic about the irrational case.)

Returning to the main plot, the 2 magenta points show that the symmetric point -z0 maps to z0 so by continuity the neighborhood of -z0 maps to the plot above. In both plots the blue orbit is the orbit of z = 0. This is always a 'critical' point of *f*, since the derivative vanishes there. On the left we can see that the blue orbit forms the boundary of the invariant region. The interior of this region is called the Siegel Disk of z0. The boundary curve is a continuum but it is never (real) analytic. However this does not prevent it from possibly being smooth. (No examples of this are known for quadratic polynomials.) As expected, the Siegel Disk is in the Fatou set of *f* and the boundary points are in the Julia set.

Fatou and Julia knew that for rational complex functions, the orbit of the critical point can be used to characterize the dynamics.

**Definition**: The Julia set, Jc of $f(z)= z^2+ c$ is the closure of the set of repelling periodic points. (A point p of period k is repelling if $|f^{k}\,'(p)| > 1$.)

Thus $J_c$ consists of all repelling periodic points and their limit points. An equivalent characterization is: Jc is the boundary of the set of points which diverge under $f^{k}$.

**Theorem** (P. Fatou & G. Julia, 1918) Let Ω denote the set of critical points for a polynomial function *f* and let K denote the set of points which do not diverge to ∞
(i) $\Omega \subseteq K \Leftrightarrow$ J is connnected
(ii) $\Omega \cap K = \varnothing \Leftrightarrow$ J is a Cantor set

A polynomial such as $f(z) = z^2 + c$ will always have z = 0 as a critical point and z = ∞ is another critical point - but it is not very interesting because it is always attracting. This leaves z = 0 as the only interesting point. If it does not diverge than J is connected

**Definition**:The Mandelbrot set $M = \{c \in \mathbb{C} : Jc \text{ is connected}\} = \{c \in \mathbb{C} : f^k(0) \text{ does not diverge}\}$

It appears from the blue orbit above, that z = 0 does not diverge so $c_0$ is in M, and it must be on the boundary because the interior of M consists of c values where the fixed point is attracting. If c corresponded to a rationally indifferent fixed point, such as the period 5 case above, it would be the cusp of a period 5 bud on the boundary of M. (These are also called parabolic points). The only other possibility for the Golden Mean case is a boundary point which is not a cusp.

Below are the two Julia sets (complements of Peitgen & Richter [PR]). Both are (marginally) connected, but in the rational case the chaotic dynamics reach all the way into the fixed point z0. For the Golden Mean case on left, the invariant curves around z0 are a sign of stability - just as they are for the KAM Theorem.These invariant curves are shown here just for reference- they are not part of the Julia set - which is the chaotic boundary.

.

The Siegel disk boundary forms part of the Jilia set and the rest are preimages. A typical point in one of the smaller buds would progress from smaller buds to larger buds and finally it would be mapped to one of the invariant curves. For polynomial functions, M. Herman proved that the Julia set always contains the critical point so z = 0 is in the Julia set and if we zoom into this region it has a self-similar fractal boundary

**Theorem** (C. Siegel, A. Brjuno and J. Yoccoz). A quadratic polynomial *f* is linearizable near 0 if and only if $\sum_{n-1}^{\infty} \frac{\log q_{n+1}}{q_n} < \infty$ where $p_n/q_n$ are the rational convergents of θ.

The work of Siegel, Brjuno and Yoccoz was concerned entirely with the local dynamics at z0 and very little is known about the global dynamics in the region of a Siegel disk D. The boundary forms a non-analytic continuum but is it a Jordan curve ? It is possible that for some Siegel disks, the boundary might be an indecomposable continuum. On the other extreme there are rotation numbers θ which yield smooth boundaries - but no examples are known for quadratic polynomials.

**Return to DKHO**

The DKHO map from Section 7 can be written as $z_{k+1} = e^{-i\omega}(z_k - i\,\text{sgn}(\text{Re}(z_k)))$
where z = x + iy. In Mathematica: **F[z_]:=Exp[-I*w]*(z - I*Sign[Re[z]])**. To get webs we will scan the y axis, but to make the plots match those of DKHO, we will swap Real and Imaginary.

```
 w=N[ 2*Pi/7]; H0=Table[ x*I, {x, -4, 4, .005}];
Web=Flatten[Table[NestList[F, H0[[k]], 500],{k,1,Length[H0]}],1];
Graphics[{AbsolutePointSize[1.0],Blue,Point[{Im[#],Re[#]}&/@Web]},Axes->True]
```

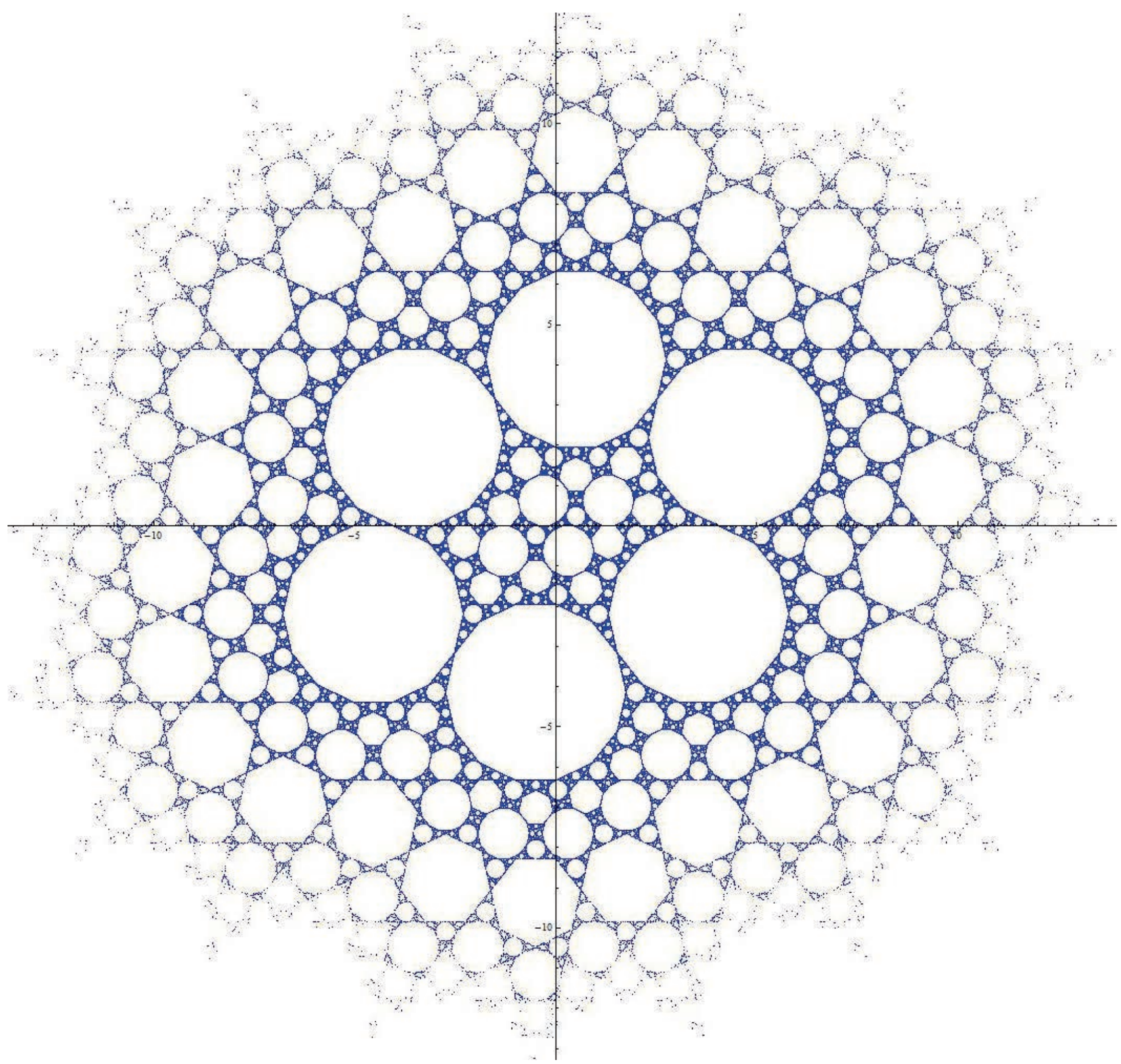

This matches the DKHO web for the regular heptagon, N = 7. We remarked earlier that outside the first ring of D's, the web is locally faithful to the Tangent Map web. But inside the ring of D's, the dynamics are highly irregular as can be seen above. The problem is that the kicks are not synchronized with a regular heptagon. This is easy to rectify - just generate the kicks with y instead of x: **F[z_]:=Exp[-I*w]*(z - Sign[Im[z]]).**

This is clearly not the same as swapping x and y, but that step is no longer necessary because the new web is already aligned with the x axis. It is no longer necessary to swap Real and Imaginary.

**H0=Table[ x , {x, -4, 4, .005}];** (*scan the x axis*) **w=N[2*Pi/7];**
**Web=Flatten[Table[NestList[F,H0[[k]],500],{k,1,Length[H0]}],1];**
**Graphics[{AbsolutePointSize[1.0],Blue,Point[{Re[#],Im[#]}&/@Web]},Axes->True]**
(*On the left below*)

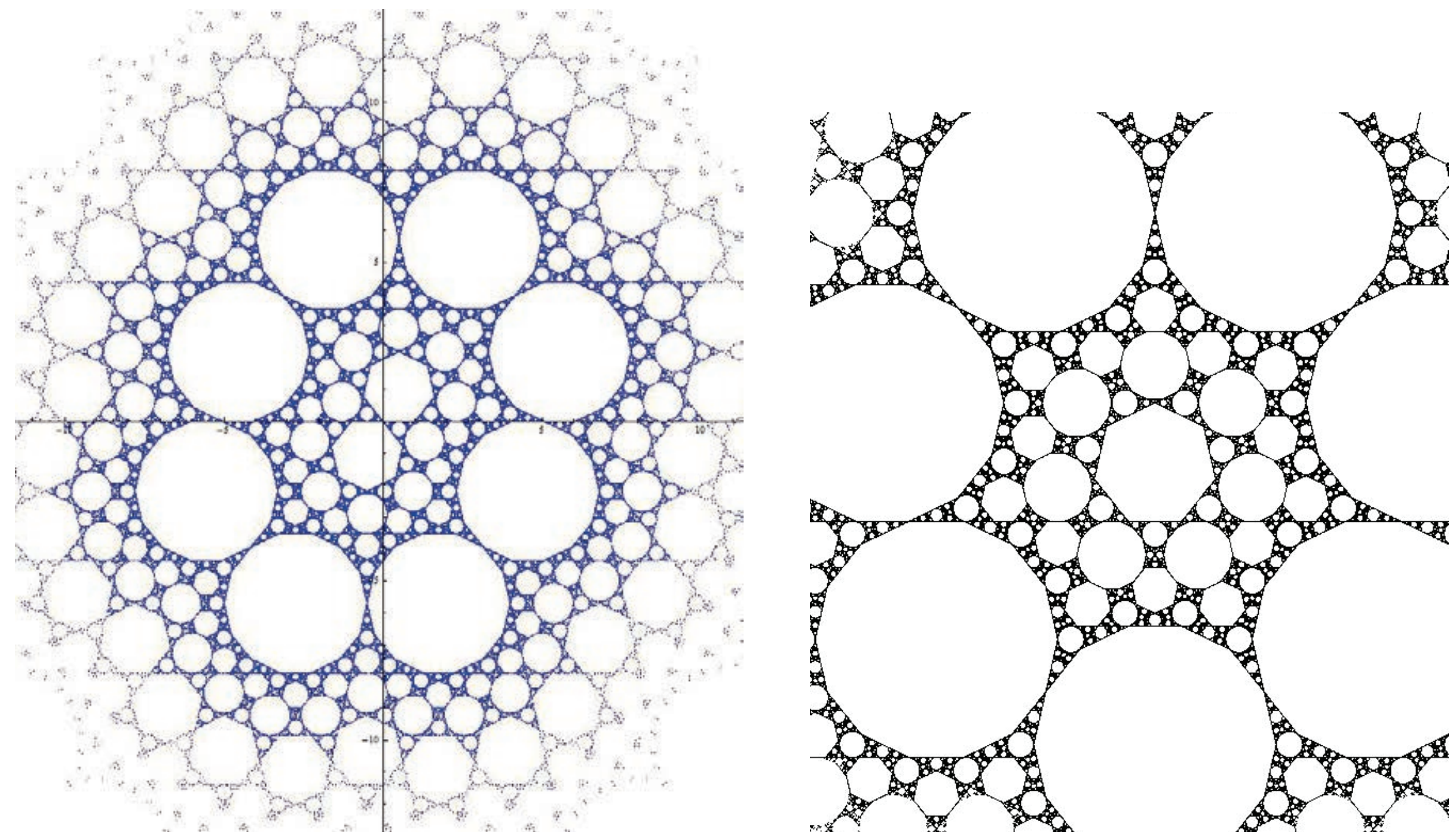

This is a perfect Tangent Map web for y ≥ 0 and a perfect reflection below. For the sake of comparison we have reproduced the actual N = 7 web on the right. This new map appears to work for any regular N-gon. We call it the 'Y' version of DKHO:

$$\begin{bmatrix} x_{k+1} \\ y_{k+1} \end{bmatrix} = \text{DKHOY}\begin{bmatrix} x_k \\ y_k \end{bmatrix} = \begin{bmatrix} \cos\omega & \sin\omega \\ -\sin\omega & \cos\omega \end{bmatrix}\begin{bmatrix} x_k - \operatorname{sgn} y_k \\ y_k \end{bmatrix}$$

**DKHOY[{x_,y_}]:= {(x - Sign[y])*Cos[w]+ y*Sin[w], (Sign[y]-x)*Sin[w]+y*Cos[w]}**;

So the kicks are now dependent on the 'momentum' but they alter the 'position'. This puts them in synchronization with the N-gon (or the other way around). Below are the first 3 stages in the new web:

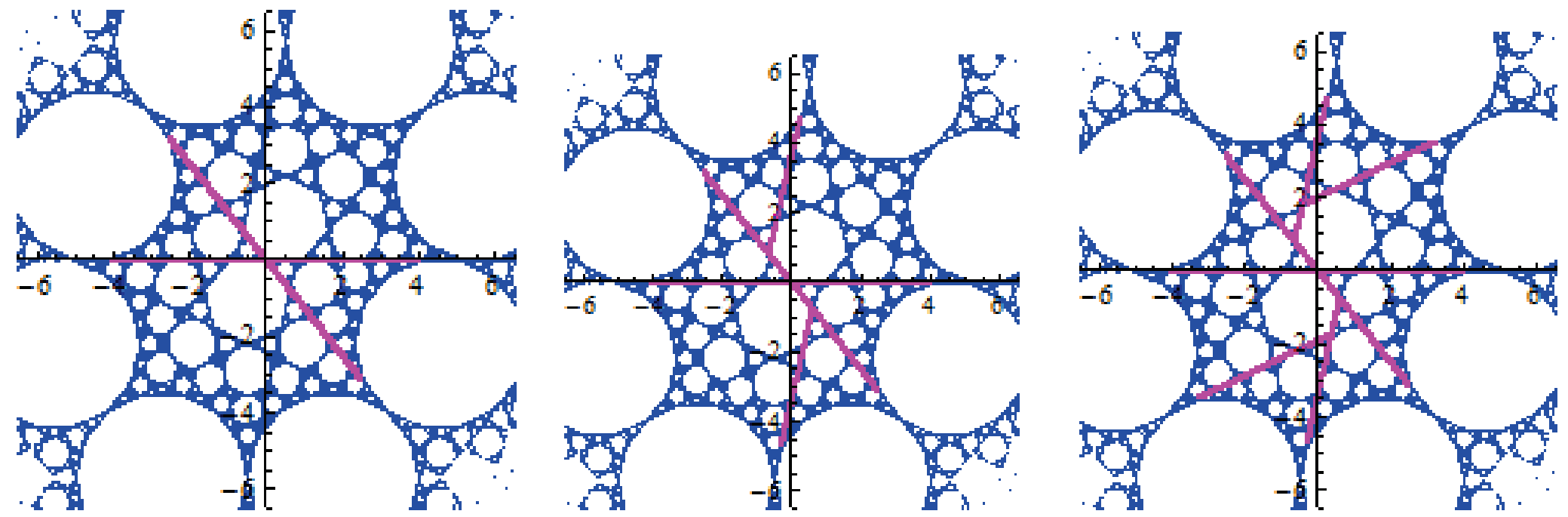

These webs are clones of the Tangent Map webs, but they are generated in a much more efficient manner, from just a single line of discontinuity.

There are broad classes of related mappings based on discontinuous kicked oscillators, and many of these exhibit some form of chaotic 'transport' similar to the plots above for the Harper Kicked Map.

The general form for such maps is $\begin{bmatrix} x_{k+1} \\ y_{k+1} \end{bmatrix} = \begin{bmatrix} \cos\omega & \sin\omega \\ -\sin\omega & \cos\omega \end{bmatrix} \begin{bmatrix} x_k \\ y_k + f(x_k) \end{bmatrix}$

The Tangent Map is in this class and it exhibits a form of boundary transport for many non-regular polygons. But these mappings are not chaotic in the normal sense of the term because one key ingredient of a chaotic map is at least one positive Lyapunov exponent. The Tangent map has zero Lyapunov exponents, yet it often has the fractal web structure which is associated with chaotic systems. For this reason some authors refer to this as 'pseudo-chaos'.

Peter Ashwin noted that the complex form of DKHO is conjugate to the 'dual' center maps of Arek Goetz which have the form:

$$\mathrm{T}(z) = \begin{cases} e^{i\theta}(z - z_0) + z_0 \ \text{ if Re(z)} < 0 \\ e^{i\theta}(z - z_1) + z_1 \ \text{ if Re(z)} \geq 0 \end{cases}$$

Recall the DKHO complex form was $\mathrm{F(z)} = e^{-i\omega}(z - i\,\mathrm{sgn}(\mathrm{Re}(z)))$ which can be written as $\mathrm{F(z)} = e^{-i\omega}(z) + W\,\mathrm{sgn}(\mathrm{Re}(z))$ where $\mathrm{W} = -i/e^{\mathrm{i}\omega}$. This has the same form as T(z) above.